\documentclass[11pt,reqno]{amsart}

\usepackage{amsthm}
\usepackage{amsfonts}
\usepackage{mathrsfs}
\usepackage{amsmath}
\usepackage{amssymb, amsmath}
\usepackage{color}
\usepackage{graphicx}
\usepackage{dsfont}

\newtheorem{Theorem}{Theorem}[section]
\newtheorem{Proposition}{Proposition}[section]
\newtheorem{Lemma}{Lemma}[section]

\newtheorem{Remark}{Remark}[section]
\newtheorem{Definition}{Definition}[section]
\numberwithin{equation}{section}

\newcommand{\non}{\nonumber}
\newcommand{\no}{\noindent}

\def\ef{\hphantom{MM}\hfill\llap{$\square$}\goodbreak}

\def\eqdefa{\buildrel\hbox{\footnotesize def}\over =}

\newcommand{\beq}{\begin{equation}}
\newcommand{\eeq}{\end{equation}}
\newcommand{\ben}{\begin{eqnarray}}
\newcommand{\een}{\end{eqnarray}}
\newcommand{\beno}{\begin{eqnarray*}}
\newcommand{\eeno}{\end{eqnarray*}}
\newcommand{\bali}{\begin{aligned}}
\newcommand{\eali}{\end{aligned}}

\newcommand{\al}{\alpha}

\newcommand{\ve}{\varepsilon}
\newcommand{\ga}{\gamma}

\newcommand{\Om}{\Omega}
\newcommand{\na}{\nabla}

\newcommand{\f}{\frac}

\newcommand{\ud}{\mathrm{d}}
\newcommand{\ue}{\mathrm{e}}
\newcommand{\uu}{\mathbf{u}}
\newcommand{\vv}{\mathbf{v}}
\newcommand{\ii}{\mathbf{i}}
\newcommand{\jj}{\mathbf{j}}
\newcommand{\ww}{\mathbf{w}}

\newcommand{\xx}{\mathbf{x}}
\newcommand{\yy}{\mathbf{y}}
\newcommand{\nn}{\mathbf{n}}
\newcommand{\ee}{\mathbf{e}}
\newcommand{\hh}{\mathbf{h}}
\newcommand{\mm}{\mathbf{m}}
\newcommand{\kk}{\mathbf{k}}

\newcommand{\MM}{\mathbf{M}}
\newcommand{\NN}{\mathbf{N}}
\newcommand{\RR}{\mathbf{R}}
\newcommand{\QQ}{\mathbf{Q}}
\newcommand{\DD}{\mathbf{D}}
\newcommand{\FF}{\mathbf{F}}
\newcommand{\II}{\mathbf{I}}

\newcommand{\CR}{\mathcal{R}}
\newcommand{\CG}{\mathcal{G}}
\newcommand{\CU}{\mathcal{U}}

\newcommand{\CA}{\mathcal{A}}
\newcommand{\CP}{\mathcal{P}}
\newcommand{\CH}{\mathcal{H}}
\newcommand{\BA}{\mathbf{A}}
\newcommand{\BS}{{\mathbb{S}^2}}
\newcommand{\BR}{{\mathbb{R}^3}}
\newcommand{\BOm}{\mathbf{\Omega}}

\newcommand{\Ue}{{U_{\varepsilon}}}
\newcommand{\mue}{{\mu_{\varepsilon}}}
\newcommand{\Be}{{B_{\varepsilon}}}
\newcommand{\sqe}{{\sqrt{\varepsilon}}}
\newcommand{\pa}{\partial}
\newcommand{\pai}{\partial_i}
\newcommand{\Ee}{\mathfrak{E}_\ve}
\newcommand{\Fe}{\mathfrak{F}_\ve}
\newcommand{\tf}{\widetilde{f}}
\newcommand{\htheta}{{\hat{\theta}}}
\newcommand{\hvphi}{{\hat{\varphi}}}
\newcommand{\hmm}{{\hat{\mm}}}
\newcommand{\hm}{{\hat{m}}}
\newcommand{\HMM}{{\hat{\MM}}}
\newcommand{\HM}{{\hat{M}}}
\newcommand{\tv}{\widetilde{\vv}}
\newcommand{\tD}{\widetilde{\DD}}
\newcommand{\Pin}{\mathds{P}_{\mathrm{in}}}
\newcommand{\Pout}{\mathds{P}_{\mathrm{out}}}

\begin{document}

\title[From the Doi-Onsager equation  to the Ericksen-Leslie equation]
{The small Deborah number limit of the Doi-Onsager equation  to the
Ericksen-Leslie equation}

\author{Wei Wang}
\address{School of  Mathematical Sciences, Peking University, Beijing 100871, China}
\email{wangw07@@pku.edu.cn}

\author{Pingwen Zhang}
\address{School of  Mathematical Sciences, Peking University, Beijing 100871, China}
\email{pzhang@@pku.edu.cn}

\author{Zhifei Zhang}
\address{School of  Mathematical Sciences, Peking University, Beijing 100871, China}
\email{zfzhang@@math.pku.edu.cn}

\date{\today}

\vspace{-0.3in}
\begin{abstract}
We present a rigorous derivation of the Ericksen-Leslie equation
starting from the Doi-Onsager equation. As in the fluid dynamic limit of the Boltzmann equation, we
first make the Hilbert expansion for the solution of the  Doi-Onsager equation.
The existence of the Hilbert expansion is connected to an open question whether the energy of the Ericksen-Leslie equation is dissipated.
We show that the energy is dissipated for the Ericksen-Leslie equation
derived from the Doi-Onsager equation. The most difficult step is to prove a uniform bound for the remainder in the Hilbert
expansion. This question is connected to the spectral stability of the linearized Doi-Onsager operator
around a critical point. By introducing two important auxiliary operators,  the detailed spectral information
is obtained for the linearized operator around all critical points. However, these are not enough to justify
the small Deborah number limit for the inhomogeneous Doi-Onsager equation, since the elastic stress in the velocity equation is also strongly
singular. For this, we need to establish a precise lower bound for a bilinear form associated with the linearized operator.
In the bilinear form, the interactions between the part inside the kernel
and the part outside the kernel of the linearized operator are very complicated.
We find a coordinate transform and introduce a five dimensional space called the Maier-Saupe space
such that the interactions between two parts can been seen explicitly by a delicate argument of completing the square.
However, the lower bound is very weak for the part inside the Maier-Saupe space.
In order to apply them to the error estimates,
we have to analyze the structure of the singular terms and  introduce a suitable energy functional.
\end{abstract}

\maketitle

\section{Introduction}

\subsection{The Doi-Onsager theory}

Liquid crystals are a state of matter that have properties between
those of a conventional liquid and those of a solid crystal. One of
the most common liquid crystal phases is the nematic. The nematic
liquid crystals are composed of rod-like molecules with the long
axes of neighboring molecules aligned approximately to one another.
A classic model which predicts isotropic-nematic phase transition is
the hard-rod model proposed by Onsager \cite{On}. Onsager introduced
the notion of orientational distribution function and considered a
mean-field model in which the rod-rod interaction was modeled by the
excluded volume effect. Following Onsager, Maier and Saupe \cite{MS}
proposed a slightly modified interaction potential, now known as the
Maier-Saupe potential. Doi and Edwards \cite{Doi} extended the Onsager
theory for describing the behavior of liquid crystal polymer flows.

We use $\xx\in\Omega\subseteq\BR$ to denote the material point and
$f(\xx,\mm,t)$ to represent the number density for the number of
molecules whose orientation is parallel to $\mm$ at point $\xx$ and
time $t$. For the spatially homogeneous liquid crystal flow,
the Doi-Onsager equation \cite{Doi} takes
\begin{eqnarray}\label{eq:Doi-Onsager}
\frac{\partial{f}}{\partial{t}}=\f 1 {De}\CR\cdot(\CR{f}+f\CR
U)-\CR\cdot\big(\mm\times\kappa\cdot\mm{f}\big),
\end{eqnarray}
where $De$ is the Deborah number, $\CR$ is the rotational gradient
operator(see Section 3), $\kappa$ is a constant velocity gradient, and $U$ is the
mean-field interaction potential. Onsager \cite{On} considered the
potential \beno U=\CU
f(\mm,t)=\alpha\int_{\BS}|\mm\times\mm'|f(\mm',t)d\mm', \eeno where
$\al$ is a parameter that measures the potential intensity. In this
paper, we will use the Maier-Saupe potential \cite{MS} defined by
\beno\label{eq:MS potential}
U=\alpha\int_{\BS}|\mm\times\mm'|^2f(\mm',t)d\mm'. \eeno This model
has a free energy
\begin{equation}\label{eq:free energy-homo}
A[f]=\int_{\BS}\big(f(\mm,t)\ln{f(\mm,t)}+\frac12f(\mm,t)U(\mm,t)\big)\ud\mm
\end{equation}
as its Lyapunov functional. The chemical potential is given by \beno
\mu=\frac {\delta A}{\delta f}=\ln f+U. \eeno The equation
(\ref{eq:Doi-Onsager}) can be written as \beno
\frac{\partial{f}}{\partial{t}}=\f 1
{De}\CR\cdot\big(f\CR\mu\big)-\CR\cdot\big(\mm\times\kappa\cdot\mm{f}\big).
\eeno The stress tensor is given by
\begin{eqnarray}\label{eq:tensor-hom}
\sigma^{De}=\frac{1}{2}\DD:\langle\mm\mm\mm\mm\rangle_f-\frac{1}{De}\langle\mm\mm\times\CR\mu\rangle_f,
\end{eqnarray}
where $\DD=\f1 2(\kappa+\kappa^T)$ is the symmetric part of
$\kappa$, and \beno \langle(\cdot)\rangle_f\eqdefa \int_{\BS}(\cdot)
f(\mm,t)\ud{\mm}. \eeno

The homogeneous Doi-Onsager equation has been very successful in describing the
properties of liquid crystal polymers in a solvent. This model takes
into account the effects of hydrodynamic flow, Brownian motion and
intermolecular forces on the molecular orientation distribution.
However, it does not include effects such as distortional
elasticity. Therefore it is valid only in the limit of spatially
homogeneous flows.

The inhomogeneous flows were first studied by Marrucci and Greco
\cite{MG}, and subsequently by many people \cite{FLS, Wang}. Instead
of using the distribution as the sole order parameter, they used a
combination of the tensorial order parameter and the distribution
function, and used the spatial gradients of the tensorial order
parameter to describe the spatial variations. This is a departure
from the original motivation that led us to the kinetic theory.
Wang, E, Liu and Zhang \cite{WELZ} set up a formalism in which the
interaction between molecules is treated more directly using the
position-orientation distribution function via interaction
potentials. They extend the free energy (\ref{eq:free energy-homo})
to include the effects of nonlocal intermolecular interaction
through an interaction potential as follows:
\begin{eqnarray}\label{eq:free energy-inhom}
A[f]=k_BT\int_{\Omega}\int_{\BS}f(\xx,\mm,t)(\ln{f(\xx,\mm,t)}-1)+\frac{1}{2k_BT}U(\xx,\mm,t)f(\xx,\mm,t)\ud\mm\ud\xx,
\end{eqnarray}
where $k_B$ is the Boltzmann constant, $T$ is the absolute
temperature, and the mean-field interaction potential $U$ is defined
by \beno
U(\xx,\mm,t)=k_BT\int_{\Omega}\int_{\BS}B(\xx,\xx',\mm,\mm')f(\xx',\mm',t)\ud\mm'\ud\xx'.
\eeno Here $B(\xx,\xx';\mm,\mm')$ is the interaction kernel between
the two polymers in the configurations $(\xx,\mm)$ and
$(\xx',\mm')$. It should be symmetric with respect to the
interchange of $\mm$ and $\mm '$, $\xx$ and $\xx '$. $B$ is often
translation invariant and hence it can be written in the form \beno
B(\xx-\xx';\mm,\mm'). \eeno In this paper, we take the following
form  as in \cite{EZ, YZ}:
$$B(\xx,\xx',\mm,\mm')=\alpha|\mm\times\mm'|^2\frac{1}{L^3}g\big(\frac{\xx-\xx'}{L}\big),$$
where $L$ is the lenth of the rods, and $g(\xx)$ is a radial Schwartz
function with $\int_\BR g(\xx)\ud{\xx}=1.$
This potential neglects the interaction between orientation and
position. But is sufficient in many cases. The chemical potential is
given by
$$\mu=\frac{\delta{A[f]}}{\delta{f}}=k_BT\ln{f(\xx,\mm,t)+U(\xx,\mm,t)}.$$

The inhomogeneous Doi-Onsager equation  takes the form
\begin{align*}
\frac{\pa{f}}{\pa{t}}+\vv\cdot\nabla{f}=&\frac{1}{k_BT}\nabla\cdot\big\{
\big(D_{\|}\mm\mm+D_{\bot}(\II-\mm\mm)\big)\cdot(\nabla\mu)f\big\}\nonumber\\
&+\frac{D_r}{k_BT}\CR\cdot(f\CR\mu)
-\CR\cdot(\mm\times\kappa\cdot\mm{f}),\\\nonumber
\frac{\pa{\vv}}{\pa{t}}+\vv\cdot\nabla\vv=&-\nabla{p}+\nabla\cdot\tau+\textbf{F}^e,\quad \nabla\cdot\vv=0.
\end{align*}
Here $D_{\|}$ and $D_{\bot}$ are respectively the translational diffusion coefficients parallel and
normal to the orientation of the LCP molecule, $D_r=\frac {k_BT}
{\xi_r}$ is the rotary diffusivity, $\nabla$ is the gradient
operator with respect to the spatial variable $\xx$. The total
stress $\tau$ is the sum of the viscous stress $\tau^s$ and the
elastic stress $\tau^e$. There are two contributions to the viscous
stress, one from the solvent and the other from the constraint force
arising from the rigidity of the rod (derived in \cite{Doi}), \beno
\tau^s=2\eta_s\DD+\frac12\xi_r\DD:\langle\mm\mm\mm\mm\rangle_f,
\eeno where $\DD=\f12\big(\kappa+\kappa^T\big)$, $\kappa=(\na
\vv)^T$ is the velocity gradient tensor, $\eta_s$ is the solvent
viscosity. The elastic stress $\tau^e$ and body force $\FF^e$ are
given by \beno \tau^e=-\langle\mm\mm\times\CR\mu\rangle_f,\quad
\FF^e=-\langle\nabla\mu\rangle_f. \eeno

Let $L_0$ be the typical size of the flow region, $V_0$ be the
typical velocity scale, $T_0=\f{L_0}{V_0}$ be a typical convective
time scale. Another important time scale is the relaxational time
scale due to orientation diffusion: $T_r=\f{\xi_r}{k_BT}$. The ratio
of these two time scales is an important parameter called the
Deborah number \beno De=\f {T_r} {T_0}=\f {\xi_rV_0} {k_BTL_0}.
\eeno Let $\eta_p=\xi_r, \eta=\eta_s+\eta_p, \ga=\eta_s/\eta,$ and
$Re=\f {V_0L_0} {\eta}$ be the Reynolds number. We denote \beno
&&\Ue(\xx,\mm,t)=\int_{\Omega}\int_{\BS}\Be(\xx,\xx',\mm,\mm')f(\xx',\mm',t)\ud\mm'\ud\xx',\\
&&\tau^e_\ve=-\langle\mm\mm\times\CR\mue\rangle_f,\quad
\FF^e_\ve=-\langle\nabla\mue\rangle_f, \eeno where the small
parameter $\sqrt{\ve}=\f {L} {L_0}$ represents the typical
interaction distance and \beno
&&\Be(\xx,\xx',\mm,\mm')=\alpha|\mm\times\mm'|^2\frac{1}{\sqe^{3}}g\big(\frac{\xx-\xx'}{\sqe}\big),\\
&&\mue=\ln{f(\xx,\mm,t)+U_\ve(\xx,\mm,t)}. \eeno We set \beno
f'(\xx,\mm,t)=f(L_0\xx, \mm, T_0t ),\quad \vv'(\xx,t)=\vv(L_0\xx,
T_0t)/V_0. \eeno
Then the non-dimensional Doi-Onsager equation takes the following form (drop
the prime for the simplicity):
\begin{eqnarray}\label{eq:LCP-non}
&&\frac{\pa{f}}{\pa{t}}+\vv\cdot\nabla{f}=\frac{\ve}{De}\nabla\cdot\big\{
\big(\gamma_{\|}\mm\mm+\gamma_{\bot}(\II-\mm\mm)\big)\cdot(\nabla{f}+f\nabla{\Ue})\big\}\nonumber\\\label{Doi-f}
&&\qquad\qquad+\frac{1}{De}\CR\cdot(\CR{f}+f\CR{\Ue})
-\CR\cdot(\mm\times\kappa\cdot\mm{f}),\\\nonumber
&&\frac{\pa{\vv}}{\pa{t}}+\vv\cdot\nabla\vv=-\nabla{p}+\frac{\gamma}{Re}\Delta\vv
+\frac{1-\gamma}{2Re}\nabla\cdot(\DD:\langle\mm\mm\mm\mm\rangle_f)+\frac{1-\gamma}{De
Re}(\nabla\cdot\tau^e_\ve+\FF^e_\ve),
\end{eqnarray}
where \beno \gamma_{\|}=\f {L_0De} {V_0L^2}D_\|,\quad
\gamma_{\bot}=\f {L_0De} {V_0L^2}D_\bot. \eeno
The system (\ref{eq:LCP-non}) has the following energy dissipation relation:
\begin{eqnarray}
&&-\frac{\ud}{\ud{t}}\Big(\int_{\Omega}\frac{1}{2}|\vv|^2\ud\xx+\frac{1-\gamma}{DeRe}A_\ve[f]\Big)\non\\
&&\quad=\int_{\Omega}\frac{\gamma}{Re}\DD:\DD+\frac{1-\gamma}{2Re}\langle(\mm\mm:\DD)^2\rangle
+\frac{1-\gamma}{De^2Re}\big\langle\CR\mue\cdot\CR\mue\big\rangle\nonumber\\
&&\quad\qquad+\frac{\ve}{De^2Re}\big\langle\nabla\mue\cdot(\gamma_{\|}\mm\mm+\gamma_\bot(\II-\mm\mm))\cdot\nabla\mue\big\rangle
\ud\xx,
\end{eqnarray}
where \beno
A_\ve[f]=\int_{\Omega}\int_{\BS}f(\xx,\mm,t)(\ln{f(\xx,\mm,t)}-1)+\frac{1}{2}U_\ve(\xx,\mm,t)f(\xx,\mm,t)\ud\mm\ud\xx.
\eeno
We refer to \cite{YZ, ZZ-SIAM} for the numerical study and the well-posedness of the system (\ref{eq:LCP-non}).

\subsection{The Ericksen-Leslie theory}
Ericksen-Leslie theory \cite{Eri,Les} is an elastic continuum
theory. The liquid crystal material is treated as a continuum and
molecular details are entirely ignored, and this theory considers
perturbations to a presumed oriented sample. Elastic continuum
theory is a very powerful tool for modeling liquid crystal devices.

The configuration of the liquid crystals is described by a director
field $\nn(\xx, t)$. The hydrodynamic equation takes the form
\begin{eqnarray}\label{eq:EL-v}
\frac{\pa{\vv}}{\pa{t}}+\vv\cdot\nabla\vv=-\nabla{p}+\frac{\gamma}{Re}\Delta\vv
+\frac{1-\gamma}{Re}\nabla\cdot\sigma,
\end{eqnarray}
where the stress $\sigma$ is modeled by the phenomenological
constitutive relation: \beno \sigma=\sigma^L+\sigma^E. \eeno Here
$\sigma^L$ is the viscous (Leslie) stress
\begin{eqnarray}\label{eq:Leslie stress}
\sigma^L=\alpha_1(\nn\nn\cdot\DD)\nn\nn+\alpha_2\nn\NN+\alpha_3\NN\nn+\alpha_4\DD
+\alpha_5\nn\nn\cdot\DD+\alpha_6\DD\cdot\nn\nn \een
with \beno
\NN=\frac{\pa\nn}{\pa{t}}+\vv\cdot\nabla\nn+\BOm\cdot\nn,\quad\BOm=\frac{1}{2}(\kappa^T-\kappa).
\eeno The six constants $\al_1, \cdots, \al_6$ are called the Leslie
coefficients. Parodi's relation \cite{Pa} gives a constraint for
Leslie coefficients: $\alpha_2+\alpha_3=\alpha_6-\alpha_5$. While,
$\sigma^E$ is the elastic (Ericksen) stress
\begin{eqnarray}\label{eq:Ericksen}
\sigma^E=-\frac{\partial{E_F}}{\partial(\nabla\nn)}\cdot(\nabla\nn)^T,
\end{eqnarray}
where $E_F=E_F(\nn,\nabla\nn)$ is the Frank energy. The dynamic
equation for the director field is given by
\begin{eqnarray}\label{eq:EL-n}
&&\nn\times\big(\hh-\gamma_1\NN-\gamma_2\DD\cdot\nn\big)=0,
\end{eqnarray}
where $\gamma_1=\alpha_3-\alpha_2, \gamma_2=\alpha_6-\alpha_5$, and
$\hh$ is the molecular field \beno
\hh=-\frac{\delta{E_F}}{\delta{\nn}}=
\nabla\cdot\frac{\partial{E_F}}{\partial(\nabla\nn)}-\frac{\partial{E_F}}{\partial\nn}.
\eeno In this paper, we will consider
$E_F=\frac{k}{2}\int_{\Om}|\nabla\nn(\xx)|^2\ud\xx$. In this case,
we have \ben
\hh=k\Delta\nn,\quad\sigma^E=-k\nabla\nn\odot\nabla\nn=-k(\nabla_i{n_k}\nabla_jn_k)_{3\times3}.
\een

The energy dissipation for Ericksen-Leslie equation is given by
\begin{eqnarray}
&&-\frac{\ud}{\ud{t}}\Big(\int_{\Omega}\frac{Re}{2(1-\gamma)}|\vv|^2\ud\xx+E_F\Big)\non\\
&&\quad=\int_{\Omega}\Big(\frac{\gamma}{1-\gamma}|\nabla\vv|^2+(\alpha_1+\frac{\gamma_2^2}{\gamma_1})|\DD:\nn\nn|^2
+\alpha_4\DD:\DD\nonumber\\
&&\qquad\qquad+(\alpha_5+\alpha_6-\frac{\gamma_2^2}{\gamma_1})|\DD\cdot\nn|^2
+\frac{1}{\gamma_1}|\nn\times\hh|^2\Big)\ud\xx.\quad\label{EL_energy_law}
\end{eqnarray}
We refer to \cite{EZ} for a derivation of (\ref{EL_energy_law}). Concerning the mathematical study of the simplified
Ericksen-Leslie equation, we refer to \cite{Lin, Lin1, Lin2, Lin3} and references therein.

\subsection{From the Doi-Onsager theory to the Ericksen-Leslie theory}

Two kinds of theories were put forward to investigate the liquid
crystalline polymers  from the different points of view. The
Ericksen-Leslie theory is phenomenological in nature, and will be
become invalid near defects where the director cannot be defined.
The Ericksen-Leslie equation contain six unknown parameters called
the Leslie coefficients, which are difficult to determine by using
experimental results. Especially, whether the energy defined in
(\ref{EL_energy_law}) is dissipated remains unknown in Physics.
Hence, it is very important to establish the relationship between
two theories.

Kuzuu and Doi \cite{KD} formally derive the Ericksen-Leslie equation
from the Doi-Onsager equation (\ref{eq:Doi-Onsager}), and determine
the Leslie coefficients. However, the Ericksen stress is missed in
the homogeneous case. E and Zhang \cite{EZ} extend Kuzuu and Doi's
formal derivation to the inhomogeneous case. To recover the Ericksen
stress, they find that the Deborah number $De$ and the interaction
distance $\sqrt{\ve}$ should satisfy  $De\sim \ve$.

Roughly speaking, Kuzuu and Doi shows that when the Deborah number
is small, the solution $f$ of (\ref{eq:Doi-Onsager}) has the formal
expansion \beno f=f_0(\mm\cdot\nn)+\ve f_1+\cdots, \eeno where
$f_0(\mm\cdot\nn)$ denotes the equilibrium distribution function
satisfying
\beno \CR\cdot(\CR{f_0}+f_0\CR\CU f_0)=0,
\eeno
and $\nn$ is determined by (\ref{eq:El-home-4}). E and Zhang shows that the
solution $(f,\vv)$ of (\ref{eq:LCP-non}) has the formal expansion
\beno
&&f=f_0(\mm\cdot\nn)+\ve f_1+\cdots,\\
&&\vv=\vv_0+\ve\vv_1+\cdots, \eeno
where $(\vv_0,\nn)$ is determined by (\ref{eq:EL-v}) and (\ref{eq:EL-n}).

The main goal of this paper is to give a rigorous derivation of the
Ericksen-Leslie equation from the Doi-Onsager equation. This is a singular small Deborah number
limit problem. To justify this limit, we first make the Hilbert
expansion for the solution of the Doi-Onsager equation, then
show that the error term is small in a suitable Sobolev space.
The existence of the Hilbert expansion is connected to the question whether the energy of the Ericksen-Leslie equation is dissipated.
We will show that the energy is dissipated for the Ericksen-Leslie equation
derived from the Doi-Onsager equation. The error
estimates rely heavily on the  spectral analysis of the
linearized Doi-Onsager operator around the critical point, which includes

\begin{itemize}

\item[1.] Give a complete classification for all critical points $h$ of $A[f]$, which satisfies
\beno \CR\cdot(\CR{h}+h\CR\CU h)=0. \eeno

\item[2.] The spectral analysis of the linearized Doi-Onsager  operator $\CG_h$ around a critical point $h$ defined by
\beno \CG_hf=\CR\cdot\big(\CR{f}+h\CR\CU{f}+f\CR\CU{h}\big). \eeno

\item[3.] Establish a precise lower bound for the bilinear form $\big\langle \CG_{h}^\ve f,\CH^\ve_h f\big\rangle$ with
\beno \CG_h^\ve f=\CR\cdot\big(\CR{f}+h\CR \CU_\ve f+f\CR\CU
h\big),\quad\CH^\ve_h f=\f f {h}+\CU_\ve f. \eeno
\end{itemize}
The first point has been given by the second author and coworkers
\cite{LZZ}. The second point and the third point are completely new.
To prove the second point, we introduce two important auxiliary operators
$\CA_h$ and $\CH_h$ defined by \beno \CA_h{f}=-\CR\cdot(h\CR
f),\quad \CH_h=\f f {h}+\CU f. \eeno It is easy to see that
$\CG_h=-\CA_h\CH_h$ and $\CH_h$ is self-adjoint. Then we reduce the
spectral analysis of $\CG_h$ to that of $\CH_h$. The proof of the
third point is very subtle. Since the orthogonal structure is
destroyed when $\ve \neq 0$, the interactions between the part
inside the kernel of $\CG_h^\ve$ and the part outside the kernel
become very complicated. To prove a lower bound, we find a
coordinate transform and introduce a generalized kernel space of
$\CG_h^\ve$(this is a five dimensional space called the Maier-Saupe
space) such that the interactions between two parts can be seen
explicitly by a delicate argument of completing the square.

With the above preparations, it is still not enough to complete the
error estimates in the inhomogeneous case. When $\ve\neq 0$,
we can only get a strong lower bound of
$\big\langle \CG_{h}^\ve f,\CH_h^\ve f\big\rangle$ for the part
outside the Maier-Saupe space, and a weak lower bound for the part
inside the Maier-Saupe space. In order to apply them to the error
estimates, we have to analyze the nonlinear interactions between two
parts for the singular term like $\f1\ve\big\langle f_R,\pa_t(\f1 {f_0})f_R\big\rangle$ and introduce a
suitable energy functional. Since we have no decay in $\ve$ for the part of $f_R$ inside the kernel,
the term $\f1\ve\big\langle f_R,\pa_t(\f1 {f_0})f_R\big\rangle$ seems to have an order $\f 1 \ve$(\textbf{Very singular}).
Surprisingly, we will show that it is bounded.

We believe that the spectral information of the linearized operator
will be very important to study the other problems like the
nonlinear stability and instability of the critical points. These
will be left to the future work.

\section{Presentation of main results}

\subsection{The homogeneous case}

We consider the homogeneous Doi-Onsager equation
\begin{align}\label{eq:Doi-Onsager-4}
\frac{\partial{f^\ve}}{\partial{t}}
=&\frac{1}{\ve}\CR\cdot\big(\CR{f^\ve}+f^\ve\CR{\CU}f^\ve\big)
-\CR\big(\mm\times(\DD-\BOm)\cdot\mm{f^\ve}\big).
\end{align}
Here $\DD=\f12(\kappa+\kappa^T), \Om=\f12(\kappa^T-\kappa),$ and
$\ve$ is the Deborah number. The corresponding stress tensor
$\sigma^\ve$ is given by
\begin{eqnarray}\label{eq:tensor-hom-4}
\sigma^{\ve}=\frac{1}{2}\DD:\langle\mm\mm\mm\mm\rangle_{f^\ve}-\frac{1}{\ve}\langle\mm\mm\times\CR\mu^\ve\rangle_{f^\ve}
\end{eqnarray}
with $\mu^\ve=\ln f^\ve+\CU f^\ve$ and  $\CU
f=\alpha\int_{\BS}|\mm\times\mm'|^2f(\mm',t)\ud\mm'$.

In the homogeneous case, the Ericksen-Leslie equation is reduced to
\begin{eqnarray}\label{eq:El-home-4}
\nn\times\big(\frac{\pa\nn}{\pa{t}}+\BOm\cdot\nn-\lambda\DD\cdot\nn\big)=0,
\end{eqnarray}
together with the stress $\sigma^L$ given by
\begin{eqnarray}\label{eq:ELhome-stress-4}
\sigma^L=\alpha_1(\nn\nn\cdot\DD)\nn\nn+\alpha_2\nn\NN+\alpha_3\NN\nn+\alpha_4\DD
+\alpha_5\nn\nn\cdot\DD+\alpha_6\DD\cdot\nn\nn.
\end{eqnarray}

Our main results are stated as follows.

\begin{Theorem}\label{thm:Home}
Let $h_{\eta, \nn}$ be a stable critical point of $A[f]$, and
$\nn(t)$ be a solution of (\ref{eq:El-home-4}) with the initial data
$\nn_0\in \BS$ and $\lambda$ given by \ben\label{eq:lambda}
\lambda(\alpha)=\frac{\big\langle3(\mm\cdot\nn)^2-1\big\rangle_{h_{\eta,
\nn}}}
{\big\langle{g_0\frac{\ud{u_0}}{\ud\theta}}\big\rangle_{h_{\eta,
\nn}}}, \quad u_0=\CU{h_{{\eta, \nn}}}, \een and $g_0$ is a solution
of (\ref{eq:g0}). Assume that the initial data $f^{\ve}_0(\mm)\in
H^1(\BS)$ with $\int_{\BS}f^\ve_0(\mm)\ud\mm=1$ takes the form \beno
f^{\ve}_0(\mm)=h_{\eta,\nn_0}(\mm)+\sum_{k=1}^3\ve^kf_{k}(\mm,0)+\ve^2f_{R,0}^\ve(\mm),
\eeno where $f_k(\mm,t)(k=1,2,3)$ is determined by Proposition
\ref{prop:Hilbert-home}, and $f_{R,0}^\ve(\mm)$ satisfies
$\|f_{R,0}^\ve\|_{H^{-1}(\BS)}\le C$. Then for any $T>0$, there
exists an $\ve_0>0$ such that for each $0<\ve<\ve_0$, the solution
$f^\ve(\mm,t)$ of (\ref{eq:Doi-Onsager-4}) takes the form \beno
f^\ve(\mm,t)=h_{\eta, \nn(t)}(\mm)+\sum_{k=1}^3\ve^kf_k(\mm,
t)+\ve^2f_R^\ve(\mm,t), \eeno where $f_R^\ve(\mm,t)$ satisfies \beno
\|f_R^\ve(t)\|_{H^{-1}(\BS)}\le C \quad\textrm{ for any }\quad t\in
[0,T]. \eeno
\end{Theorem}

\begin{Remark}
The Ericksen-Leslie equation (\ref{eq:El-home-4}) is equivalent to
\beno
\frac{\pa\nn}{\pa{t}}+\BOm\cdot\nn-\lambda(\II-\nn\nn)\DD\cdot\nn=0.
\eeno It is easy to show that it has a unique global solution.
\end{Remark}

Let $S_2=\langle{P}_2(\mm\cdot\nn)\rangle_{h_{\eta,\nn}}$ and
$S_4=\langle{P}_4(\mm\cdot\nn)\rangle_{h_{\eta,\nn}}$, where
$P_k(x)$ is the $k$-th Legendre polynomial. We take the Leslie
coefficients $\al_1,\cdots,\al_6$ in the definition of $\sigma^L$ as
follows \ben
&&\alpha_1=-\frac{S_4}{2},\quad\alpha_2=-\frac{1}{2}\big(1+\frac{1}{\lambda}\big)S_2,\quad\alpha_3=-\frac{1}{2}\big(1-\frac{1}{\lambda}\big)S_2,
\label{Leslie cofficients1}\\
&&\alpha_4=\frac{4}{15}-\frac{5}{21}S_2-\frac{1}{35}S_4,\quad
\alpha_5=\frac{1}{7}S_4+\frac{6}{7}S_2,\quad
\alpha_6=\frac{1}{7}S_4-\frac{1}{7}S_2.\label{Leslie cofficients2}
\een Then we have
\begin{Theorem}\label{thm:stress-home}
Let $p(t)=-\frac{S_2}{14}\DD:\nn\nn$. For any $T>0$, there exist an
$\ve_0>0$ such that for each $0<\ve<\ve_0$,
there holds \beno |\sigma^\ve(t)-\sigma^L(t)-p(t)\II|\le {C\ve}\quad
\textrm{for}\quad t\in [0,T]. \eeno

\end{Theorem}

\subsection{The inhomogeneous case}
In order to derive the Ericksen-Leslie equation with the Ericksen
stress, we have to consider the system (\ref{eq:LCP-non}) with $De=\ve$.
For the simplicity of presentation, we will consider the case
when the translational diffusion coefficients vanish. Then the non-dimensional
Doi-Onsager equations takes
\begin{eqnarray}\label{eq:LCP-nonL-f}
&&\frac{\pa{f^\ve}}{\pa{t}}+\vv^\ve\cdot\nabla{f^\ve}=\frac{1}{\ve}\CR\cdot(\CR{f^\ve}+f^\ve
\CR\CU_\ve f^\ve)
-\CR\cdot(\mm\times\kappa^\ve\cdot\mm{f^\ve}),\\
&&\frac{\pa{\vv^\ve}}{\pa{t}}+\vv^\ve\cdot\nabla\vv^\ve=-\nabla{p^\ve}+\frac{\gamma}{Re}\Delta\vv^\ve
+\frac{1-\gamma}{2Re}\nabla\cdot\big(\DD^\ve:\langle\mm\mm\mm\mm\rangle_{f^\ve}\big)\non\\
&&\qquad\qquad\qquad\qquad+\frac{1-\gamma}{\ve
Re}(\nabla\cdot\tau^e_\ve+\FF^e_\ve),\label{eq:LCP-nonL-v}
\end{eqnarray}
where $\kappa^\ve=(\na v^\ve)^T,
\DD^\ve=\f12\big(\kappa^\ve+(\kappa^\ve)^T\big)$, and \beno
&&\tau^e_\ve=-\langle\mm\mm\times\CR\mue\rangle_{f^\ve},\quad
\FF^e_\ve=-\langle\nabla\mue\rangle_{f^\ve},\quad \mue=\ln f^\ve+\CU_\ve f,\\
&&\CU_\ve
f=\int_{\Omega}\int_{\BS}\alpha|\mm\times\mm'|^2\frac{1}{\sqe^{3}}g\big(\frac{\xx-\xx'}{\sqe}\big)f(\xx',\mm',t)\ud\mm'\ud\xx'.
\eeno
We also require that the Fourier transform of $g$ satisfies
\beno
0\le \hat g(\xi)<1 \quad\textrm{ for }\quad \xi\neq 0,\quad \hat g''(0)<0.
\eeno

Now we can derive the full Ericksen-Leslie equation
\begin{align}\label{eq:EL-nh}
&\nn\times\big(\hh-\gamma_1\NN-\gamma_2\DD\cdot\nn\big)=0,\\
&\frac{\pa{\vv}}{\pa{t}}+\vv\cdot\nabla\vv=-\nabla{p}+\frac{\gamma}{Re}\Delta\vv
+\frac{1-\gamma}{Re}\nabla\cdot\sigma,\label{eq:EL-vh}
\end{align}
where $\gamma_1=\alpha_3-\alpha_2$, $\gamma_2=\alpha_6-\alpha_5$, $\hh=k\Delta\nn$,
$\NN=\frac{\pa\nn}{\pa{t}}+\vv\cdot\nabla\nn+\BOm\cdot\nn$,
$\DD=\f12\big(\na\vv+(\na\vv)^T\big),
\BOm=\f12\big(\na\vv-(\na\vv)^T\big)$, and $\sigma=\sigma^L+\sigma^E$
with \beno
&&\sigma^L=\alpha_1(\nn\nn\cdot\DD)\nn\nn+\alpha_2\nn\NN+\alpha_3\NN\nn+\alpha_4\DD
+\alpha_5\nn\nn\cdot\DD+\alpha_6\DD\cdot\nn\nn,\\
&&\sigma^E=-k\nabla\nn\odot\nabla\nn. \eeno

Our main result is stated as follows.

\begin{Theorem}\label{thm:deborah-inhom-L}
Let $\Omega=\mathbb{R}^3$ and $h_{\eta, \nn}$ be a stable critical
point of $A[f]$, and let $(\nn,\vv_0)\in C\big([0,T]; H^{20}(\Om)\big)$
be a solution of (\ref{eq:EL-nh})-(\ref{eq:EL-vh}) on $[0,T]$ for some
$T>0$ with the initial data $(\nn_0, \vv_{0,0})$, $\lambda$ given by
(\ref{eq:lambda}), and Leslie coefficients defined by (\ref{Leslie
cofficients1})-(\ref{Leslie cofficients2}). Assume also that there
exist constant vector $\textbf{c}\in \BS$ and constant $c_0\in
(0,1)$ such that \ben\label{ass:degenerate}
\big|\nn(\xx,t)\times\textbf{c}\big|\ge c_0 \quad\textrm{ for
any}\quad (x,t)\in \Om\times [0,T]. \een Assume that the initial
data $\big(f^{\ve}_0(\xx,\mm), \vv_0^\ve\big) $ with
$\int_{\BS}f^\ve_0(\xx, \mm)\ud\mm=1$ takes the form \beno
&&f^{\ve}_0(\xx,\mm)=h_{\eta,\nn_0(\xx)}(\mm)+\sum_{k=1}^3\ve^kf_{k}(\xx,\mm,0)+\ve^3f_{R,0}^\ve(\xx,\mm),\\
&&\vv^\ve_0(\xx)=\sum_{k=0}^2\ve^k\vv_{k,0}(\xx)+\ve^3v_{R,0}^\ve(\xx),
\eeno
where $\big(f_1,f_2, f_3, \vv_1,\vv_2\big)$ is determined by
Proposition \ref{prop:Hilbert-inhome}, and
$\big(f_{R,0}^\ve(\xx,\mm),\vv_{R,0}^\ve(\xx)\big)$ satisfies \beno
\|f_{R,0}^\ve\|_{H^2(\Om\times\BS)}+\|\vv_{R,0}^\ve\|_{H^2(\Omega)}\le
C<+\infty,\quad \|f_{R,0}^\ve\|_{L^2(\Om\times\BS)}\le C\ve. \eeno
Then there exist $\ve_0>0$ such that for each
$0<\ve<\ve_0$, the system (\ref{eq:LCP-nonL-f})-
(\ref{eq:LCP-nonL-v}) has a unique solution $\big(f^\ve(\xx,\mm,t),
\vv^\ve(\xx,\mm,t)\big)$ on $[0,T]$ which takes the form
\begin{align*}
&f^{\ve}(\xx,\mm,t)=h_{\eta,\nn(\xx,t)}(\mm)+\sum_{k=1}^3\ve^kf_{k}(\xx,\mm, t)+\ve^3f_{R}^\ve(\xx,\mm, t),\\
&\vv^{\ve}(\xx,t)=\sum_{k=0}^2\ve^k\vv_{k}(\xx,t)+\ve^3\vv_{R}^\ve(\xx,t),
\end{align*}
where $\big(f_R^\ve, \vv_R^\ve\big)$ satisfies \beno
\big\|\big(f_{R}^\ve, \ve^{1/2}\nabla{f}_{R}^\ve,
\ve^{3/2}\Delta{f}_{R}^\ve\big)(t)\big\|_{L^2(\Om\times\BS)}
+\big\|\big(\vv_{R}^\ve, \ve\nabla{\vv}_{R}^\ve,
\ve^2\Delta{\vv}_{R}^\ve\big)(t)\big\|_{L^2(\Omega)}\le C \eeno for
any $t\in [0,T]$.
\end{Theorem}

\begin{Remark}
The non-degenerate assumption (\ref{ass:degenerate}) allows us to
construct a global coordinate transformation, which is the key to
establish a lower bound of a bilinear form associated with the linearized
operator in Section 5.
\end{Remark}

\begin{Remark}
We will study the existence of the solution for the full Ericksen-Leslie equation in a separate paper. We refer
to \cite{Lin1} for the simplified Ericksen-Leslie equation.
\end{Remark}

\section{Classification and stability of critical points of energy functional}

We consider the homogeneous energy functional $A[f]$ defined by
\beno A[f]=\int_{\BS}\big(f(\mm)\ln{f(\mm)}+\frac12f(\mm)\CU
f(\mm)\big)\ud\mm \eeno for $f\in L^2(\BS)$. We define \beno
\CP_0(\BS)=\Big\{\varphi\in
L^2(\BS):\int_\BS\varphi(\mm)\ud\mm=0\Big\}. \eeno

We are concerned with the local minimizer of $A[f]$. That is, we
find all $h\in L^2(\BS)$ such that
$$A[h+\epsilon{\phi}]\ge{A[h]}$$
for all $\phi\in\CP_0(\BS)$ when $\epsilon$ is small enough. Taking
a formal expansion, we find that
$$A[h+\epsilon{\phi}]=A[h]+\epsilon\big\langle\ln{h}+\CU{h},
\phi\big\rangle+\epsilon^2\big\langle\frac{\phi}{h}+\CU{\phi},\phi\big\rangle+O(\epsilon^3).
$$
This motivates us to introduce the following definition.
\begin{Definition}
We say that $h\in L^2(\BS)$ is a critical point of the energy
functional $A[f]$ if \beno
\frac{\delta{A[f]}}{\delta{f}}\big|_{f=h}=\ln{h}+\CU{h}=\mathrm{const}.
\eeno A critical point $h$ is said to be stable if for any
$\phi\in\CP_0(\BS)$, there holds \beno
\big\langle\frac{\phi}{h}+\CU{\phi},\phi\big\rangle \ge0. \eeno
\end{Definition}

It is easy to see that if $h$ is a critical point of $A[f]$, then
$h$ is a solution of stationary Doi-Onsager equation
\ben\label{eq:EL}
\CR\cdot\big(\CR{h}+h\CR\CU h\big)=0.
\een

A complete classification for all critical points of $A[f]$ was
given by Liu, Zhang and Zhang \cite{LZZ}; see also \cite{CKT, FS,
ZWFW}.

\begin{Proposition}\label{prop:critical point}
All the critical points of $A[f]$ take the form
$$h_{\eta,\nn}(\mm)=\frac{\ue^{\eta(\mm\cdot\nn)^2}}{\int_\BS\ue^{\eta(\mm\cdot\nn)^2}\ud\sigma}, $$
where $\nn$ is an arbitrary unit vector, and $\eta=\eta(\alpha)$ is
determined by the equation
\begin{eqnarray}\label{eta-alpha}
\frac{3\ue^\eta}{\int_0^1\ue^{\eta{z^2}}\ud{z}}=3+2\eta+\frac{\eta^2}{\alpha}.
\end{eqnarray}
Furthermore, we have
\begin{itemize}
\item For all $\alpha>0$, $\eta=0$(i.e. $h=\frac{1}{4\pi}$) is always a solution.

\item For $\alpha<\alpha^*\approx6.731393$, $\eta=0$ is the only
solution. While for $\alpha=\alpha^*$, there is another solution
$\eta=\eta^*$.

\item For $\alpha>\alpha^*$, besides $\eta=0$, there
are exactly two solutions $\eta=\eta_1(\alpha),\eta_2(\alpha)$
satisfying

\begin{itemize}

\item $\eta_1(\alpha)>\eta^*>\eta_2(\alpha)$, $\lim_{\alpha\to\alpha^*}\eta_1(\alpha)=\lim_{\alpha\to\alpha^*}\eta_2(\alpha)=\eta^*$;

\item $\eta_1(\alpha)$ is an increasing function of $\alpha$, while $\eta_2(\alpha)$ is a decreasing function;

\item $\eta_2(7.5)=0$.

\end{itemize}
\end{itemize}
\end{Proposition}

Except the monotonicity of $\eta(\alpha)$, the others  have been
proved in \cite{LZZ}. The monotonicity will be proved in Lemma
\ref{lem:integral}. 


Concerning the stability of the critical point, Zhang and Zhang
\cite{ZZ} showed

\begin{Proposition}\label{prop:energy stability}
$h=\frac{1}{4\pi}$ is a stable critical point of $A[f]$ if and only
if $\alpha<7.5$; If $\alpha>\alpha^*$, $h_{\eta_1,\nn}$ is stable,
while $h_{\eta_2,\nn}$ is unstable.
\end{Proposition}

Let us conclude this section by collecting some properties of the rotational operator,
which will be used throughout the paper.
Let $\mm\in \BS$ and $\nabla_\mm$ be the gradient operator on the
unit sphere $\BS$. The rotational gradient operator $\CR$ is defined
by
$$\CR=\mm\times\nabla_{\mm}.$$
Let $(\theta, \phi)$ be the sphere coordinate on $\BS$. Then $\CR$ can be written as
\begin{align*}
\CR=&(-\sin\phi\mathbf{i}+\cos\phi\mathbf{j})\partial_\theta
-(\cos\theta\cos\phi\mathbf{i}+\cos\theta\sin\phi\mathbf{j}-\sin\theta\mathbf{k})
\frac{1}{\sin\theta}\partial_\phi\\
\eqdefa&\ii{\CR}_1+\jj{\CR_2}+\kk{\CR_3}.
\end{align*}

The following properties can be easily verified.

\begin{itemize}

\item[1.] $\CR\cdot\CR=\Delta_\BS$;

\item[2.] $\CR_im_j=-\epsilon^{ijk}m_k$, where $\mm=(m_1, m_2, m_3)$.  If $\uu$ is a constant vector, then
\beno
\CR(\mm\cdot\uu)=\mm\times\uu,\quad
\CR\cdot(\mm\times\uu)=-2\mm\cdot\uu; \eeno

\item[3.] $[\CR_j,~\CR_k]=\epsilon^{ijk}\CR_i$;

\item[4.] $\int_{\BS}\CR f_1 f_2\ud\mm=-\int_{\BS}f_1\CR f_2\ud\mm$;

\item[5.] $[\CR,\CU]=0$.
\end{itemize}
\no Here $\Delta_\BS$ is the Laplace-Beltrami operator on $\BS$,
$\epsilon^{ijk}$ is the Levi-Civita symbol.

\section{Spectral analysis of the  linearized operator}

We linearize the Doi-Onsager equation $\CR\cdot(\CR f+f\CR\CU f)=0$ around
a critical point $h$. The linearized Doi-Onsager operator $\CG_h$ is given by
\begin{eqnarray}
\CG_hf\eqdefa \CR\cdot\big(\CR{f}+h\CR\CU{f}+f\CR\CU{h}\big).
\end{eqnarray}
We denote by $H^k(\BS)$ the Sobolev space on $\BS$, and
$H_0^k(\BS)=H^2(\BS)\cap \CP_0$. $\CG_h$ is a bounded operator from
$H^2(\BS)$ to $L^2(\BS)$, and has the discrete spectra.This section
is devoted to studying the kernel and spectra of the linearized
operator $\CG_h$ . These information will play a vital role in the
study of small Deborah limit, and will be very important in the
study of nonlinear stability and instability of the critical point.

When $h$ is a trivial critical point $h_0=\frac{1}{4\pi}$, the
linearized operator $\CG_h$ is reduced to
\beno
\CG_{h_0}f=\Delta_\BS\big(f+\frac{1}{4\pi}\CU f\big).
\eeno

\begin{Proposition}\label{prop:G-h0}
The eigenvalues of $\CG_{h_0}$ are $\lambda_k=-k(k+1)$ (for
$k\neq2$, $k\ge1$) and $-6+\frac{4\alpha}{5}$, and the corresponding
eigenfunction is the spherical harmonics $Y_{k,\ell}$ of degree $k$.
Specifically, $\CG_{h_0}$ has a positive eigenvalue if and only if
$\alpha>7.5$.
\end{Proposition}

\begin{Remark}
The critical value $7.5$ is consistent with that in Proposition
\ref{prop:energy stability} deduced from the energy stability
analysis.
\end{Remark}

\no{\bf Proof.}\,Let $\psi$ be an eigenfunction of $\CG_{h_0}$
associated with the eigenvalue $\lambda$, that is,
\beno
\CG_{h_0}\psi=\lambda\psi.
\eeno
We choose the spherical harmonics
$\{Y_{2,\ell}\}_{1\le \ell\le 5}$ of degree 2 as \beno
Y_1=(m_1^2-m_2^2),~~Y_2=m_3^2-\frac{1}{3},~~Y_3=m_1m_2,~~Y_4=m_1m_3,~~Y_5=m_2m_3.
\eeno
Then we make a spherical harmonics expansion for $\psi$:
\ben\label{eq:3.2} \psi=\sum_{i=1}^5\mu_iY_i+\sum_{k\neq
2,\ell}\mu_{k,\ell}Y_{k,\ell}. \een
We have
\beno
\Delta_\BS \CU \psi&=&-\alpha\CR\cdot\CR \int_{\BS}(\mm\cdot\mm')^2\psi(\mm')\ud\mm'\\
&=&-2\alpha\CR\cdot\int_{\BS}(\mm\cdot\mm')(\mm\times \mm')\psi(\mm')\ud\mm'\\
&=& 6\alpha\int_{\BS}(\mm\cdot\mm')^2\psi(\mm')\ud\mm'. \eeno Hence,
\begin{eqnarray}\label{eq:3.3}
\lambda{\psi}&=&\Delta_\BS \psi+\frac{3\alpha}{2\pi}\int_{\BS}(\mm\cdot\mm')^2\psi(\mm')\ud\mm'\non\\
&=& \Delta_\BS \psi+\frac{3\alpha}{2\pi}m_im_jM_{ij}
\end{eqnarray}
with $M_{ij}=\int_{\BS}\mm_i\mm_j\psi(\mm)\ud\mm$. Noting that
\beno \Delta_\BS Y_{k,\ell}=-k(k+1)Y_{k,\ell}, \eeno and plugging
(\ref{eq:3.2}) into (\ref{eq:3.3}), we find that \beno
\lambda\mu_{k,\ell}=-k(k+1)\mu_{k,\ell}. \eeno This implies that
\beno \lambda=-k(k+1)\quad \textrm{or} \quad \mu_{k,\ell}=0. \eeno
Hence, if $\lambda>0$, then $\mu_{k,\ell}=0$ and we have
$$
\psi(\mm)=\sum_{i=1}^5 \mu_iY_i.
$$
Then it follows from (\ref{eq:3.3}) that
\begin{eqnarray}
&&(\lambda+6)\sum_{i=1}^5 \mu_iY_i\non\\
&&=\frac{3\alpha}{2\pi}\int_\BS\big(\frac{1}{2}Y_1Y_1'+\frac{3}{2}Y_2Y_2'
+2Y_3Y_3'+2Y_4Y_4'+2Y_5Y_5'+\frac{1}{3}\big)\sum_{i=1}^5 \mu_iY_i\ud\mm\nonumber\\
&&=\frac{3\alpha}{2\pi}\int_\BS\big(\frac{1}{2}\mu_1Y_1Y_1'^2+\frac{3}{2}\mu_2Y_2Y_2'^2
+2\mu_3Y_3Y_3'^2+2\mu_4Y_4Y_4'^2+2\mu_5Y_5Y_5'^2\big)\ud\mm.\non
\end{eqnarray}
A direct computation shows that
$$\int_\BS Y_1^2\ud\mm=\frac{16\pi}{15},~~\int_\BS Y_2^2\ud\mm=\frac{16\pi}{45},~~\int_\BS Y_3^2\ud\mm=\frac{4\pi}{15},$$
which implies that \beno
\big(\lambda+6-\frac{4\alpha}{5}\big)\sum_{i=1}^5 \mu_iY_i=0. \eeno
Hence, $\lambda=\frac{4\alpha}{5}-6$. Specifically, $\CG_{h_0}$ has
a positive eigenvalue if and only if $\alpha>7.5$.\ef \vspace{0.1cm}

When $h=h_{\eta_i,\nn} (i=1,2)$, the problem becomes more
complicated. Kuzzu and Doi \cite{KD} conjectured that all the
eigenvalues of $\CG_h$ are non-positive, and
$\mathrm{Ker}~\CG_h=\big\{\mathbf{\Theta}\cdot\CR{h},\mathbf{\Theta}\in\BR\big\}$.
Here we will give a rigorous proof of Kuzzu and Doi's conjecture
when $h$ is a stable critical point. Let us introduce an important
operator $\CA_h$ defined by \beno
\CA_h{\phi}\eqdefa-\CR\cdot(h\CR\phi) \eeno The operator $\CA_h$ has
the following properties:

\begin{Lemma}\label{lem:A-h}
The operator $\CA_h$ is a one-one mapping from $H_0^2(\BS)$ to
$\CP_0(\BS)$. We denote by $\CA_h^{-1}$ its inverse. Then it holds
that \beno \CA_h=\CA_h^*,\quad \big\langle\CA_h\phi,
\phi\big\rangle\ge0,\quad \big\langle\CA_h^{-1}\phi,
\phi\big\rangle\ge0. \eeno
\end{Lemma}

Now we introduce another important operator $\CH_h$ defined by \beno
\CH_h{f}\eqdefa\frac{f}{h}+\CU{f}. \eeno We have the following
important relation: \ben\label{eq:relation} \CG_h f=-\CA_h\CH_h f.
\een Basically, we can reduce the spectral analysis of $\CG_h$ to
that of $\CH_h$.

\begin{Proposition}\label{prop:G-symm}
If $h$ is a critical point of $A[f]$, then $\CG_h\CA_h$ is a
symmetric operator and \beno
\big\langle\CG_h\psi,\CA^{-1}_h\phi\big\rangle=-\big\langle\CH_h\psi,\phi\big\rangle
=\big\langle\CG_h\phi,\CA^{-1}_h\psi\big\rangle. \eeno Moreover, if
$h$ is a stable critical point of $A[f]$, then $\CG_h\CA_h$ is a
non-positive operator, and $\CG_h$ has only non-positive
eigenvalues.
\end{Proposition}

\no{\bf Proof.}\,The identity follows from (\ref{eq:relation}).
Since $h$ is a critical point, we have by (\ref{eq:EL}) that \beno
\CA_h\phi=-h\CR\cdot\CR\phi-\CR{h}\cdot\CR\phi
=-h\CR\cdot\CR\phi+h\CR(\CU{h})\cdot\CR\phi. \eeno Then for any
$\psi, \phi\in H^2(\BS)$, we have
\begin{align}
\big\langle\CG_h\CA_h\psi,\phi\big\rangle
&=-\big\langle\CR\CA_h\psi+h\CR\CU\CA_h\psi+\CA_h\psi\CR\CU{h}, \CR\phi\big\rangle\nonumber\\
&=-\big\langle\CR\CA_h\psi+\CA_h\psi\CR\CU{h}, \CR\phi\rangle+\langle\CU\CA_h\psi,\CR\cdot(h\CR\phi)\big\rangle\nonumber\\
&=\big\langle\CA_h\psi, \CR\cdot\CR\phi-\CR(\CU{h})\cdot\CR\phi\big\rangle-\big\langle\CU\CA_h\psi,\CA_h\phi\big\rangle\nonumber\\
&=-\big\langle\CA_h\psi,\frac{\CA_h\phi}{h}\big\rangle-\big\langle\CU\CA_h\psi,\CA_h\phi\big\rangle.\nonumber
\end{align}
Specifically, \beno \big\langle\CG_h\CA_h\phi,\phi\big\rangle
=-\big\langle\CA_h\phi,\frac{\CA_h\phi}{h}+\CU\CA_h\phi\big\rangle.
\eeno This means that $\CG_h\CA_h$ is a non-positive operator if $h$
is a stable critical point. Furthermore, if $\phi$ is an
eigenfunction of $\CG_{h}$ associated with the eigenvalue $\lambda$,
then we have \beno 0\ge
\big\langle\CG_h\CA_h\CA_h^{-1}\phi,\CA_h^{-1}\phi\big\rangle=
\big\langle\CG_h\phi,\CA_h^{-1}\phi\big\rangle=\lambda\big\langle
\phi,\CA_h^{-1}\phi\big\rangle. \eeno Hence, $\lambda\le 0$.\ef

\vspace{0.1cm}

Now we establish a lower bound of the operator $\CH_h$.

\begin{Proposition}\label{prop:H-positive}
Let $h_1=h_{\eta_1,\nn}$. For any $f\in{\CP}_0$, there holds \beno
\big\langle\CH_{h_1}{f}, f\big\rangle\ge 0, \eeno and the
equality holds if and only if
$f\in\mathrm{span}\big\{\CR_i{h_1},i=1,2,3\big\}$.
Moreover, there exists a positive constant $c_0$ depending only on
$\eta_1$ such that if $f$ satisfies
\beno
\int_\BS f(\mm)\CA^{-1}_{h_1}\CR{h_1}\ud\mm=0, \eeno then we have a
lower bound
\beno
\big\langle\CH_{h_1}{f}, f\big\rangle\ge{c_0}\langle{f}, f\rangle.
\eeno
\end{Proposition}

We need the following key lemma.
\begin{Lemma}\label{lem:integral}
Let $\eta=\eta(\al)$ be determined by (\ref{eta-alpha}), and define
$A_k(\eta)=\int_{-1}^1z^k\ue^{\eta{z^2}}\ud{z}$. Then there hold
\beno A_{k+2}=\frac{\ue^\eta}{\eta}-(k+1)\frac{A_k}{2\eta},\quad
A_0=\alpha(A_2-A_4). \eeno Moreover,
$\frac{\partial\alpha(\eta)}{\partial\eta}>0$ when $\eta>\eta^*$;
$\frac{\partial\alpha(\eta)}{\partial\eta}<0$ when $\eta<\eta^*$.
\end{Lemma}

\no{\bf Proof.}\,The first equality can be easily verified by
integrating by parts. While, the relation (\ref{eta-alpha}) is
equivalent to
\begin{eqnarray*}
6\alpha\ue^{\eta}-(3+2\eta)\alpha{A_0}=4\eta^2A_0\Longleftrightarrow
A_0=\alpha(A_2-A_4).
\end{eqnarray*}

In order to prove the second statement, it suffices to show that the
equation $\frac{\partial\alpha(\eta)}{\partial\eta}=0$ has only one
root, since $\frac{\partial\alpha(\eta^*)}{\partial\eta}=0$. We have
\beno
&&\frac{\partial}{\partial\eta}\Big(\ue^{-\eta}\big({A}_0(A_4-A_6)-A_2(A_2-A_4)\big)\Big)\\
&&=\frac{1}{2}\frac{\partial}{\partial\eta}\int_{-1}^1\int_{-1}^1
(x^2y^4+x^4y^2-x^6-y^6+x^4+y^4-2x^2y^2)\ue^{\eta(x^2+y^2-1)}\ud{x}\ud{y}\\
&&=\frac{1}{2}\int_{-1}^1\int_{-1}^1
-(x^2-y^2)^2(1-x^2-y^2)^2\ue^{\eta(x^2+y^2-1)}\ud{x}\ud{y}<0. \eeno
Hence, ${A}_0(A_4-A_6)-A_2(A_2-A_4)=0$ has only one root. Then from
the fact that \beno
\frac{\partial\alpha(\eta)}{\partial\eta}=\Big(\frac{A_0}{A_2-A_4}\Big)'
=\frac{A_2(A_2-A_4)-A_0(A_4-A_6)}{(A_2-A_4)^2}, \eeno we know that
$\frac{\partial\alpha(\eta)}{\partial\eta}=0$ has only one root.\ef
\vspace{0.1cm}

\no{\bf Proof of Proposition \ref{prop:H-positive}.}\,Without loss
of generality, we may assume $\nn=(0,0,1)$. Introduce the sphere
coordinates $(\theta,\phi)\in [0,\pi]\times [0,2\pi]$ with
$\mm=(\sin\theta\cos\phi,\sin\theta\sin\phi, \cos\theta)$. Hence,
\beno
h_{\eta,\nn}(\mm)=\frac{\ue^{\eta(\cos\theta)^2}}{\int_\BS\ue^{\eta(\cos\theta)^2}\ud\sigma}.
\eeno We make a Fourier expansion for $f$ with respect to the
variable $\phi$: \beno
f=a_0(\theta)+\sum_{k\ge1}\big(a_k(\theta)\cos(k\phi)+b_k(\theta)\sin(k\phi)\big).
\eeno Noting that the area element
$\ud\mm=\sin\theta\ud\theta\ud\phi$, we make a change of variable
$z=\cos\theta$ to get
\begin{align*}
\big\langle\frac{f}{h_{\eta,\nn}},
f\big\rangle&=\int_\BS\ue^{\eta(\cos\theta)^2}\ud\mm\int_\BS
\ue^{-\eta\cos^2\theta}\Big(a_0^2+\frac{1}{2}\sum_{k\ge1}(a_k^2+b_k^2)\Big)\ud\mm\\
&=2\pi^2\Big(\int_{-1}^1\ue^{\eta{z^2}}\ud{z}\Big)\cdot\Big(\int_{-1}^1\ue^{-\eta{z^2}}\big\{2a_0^2
+\sum_{k\ge1}(a_k^2+b_k^2)\big\}\ud{z}\Big).
\end{align*}
Routine computations show that
\begin{align*}
\big\langle\CU{f}, f\big\rangle=&-\alpha\int_\BS\int_\BS(\mm\cdot\mm')^2f(\mm)f(\mm')\ud\mm\ud\mm'\\
=&-\alpha\sum_{i,j=1}^3\Big(\int_\BS{m_im_j}f(\mm)\ud\mm\Big)^2\\
=&-2\alpha\pi^2\Big\{\big(\int_{-1}^{1}(1-z^2)a_0\ud{z}\big)^2+
2\big(\int_{-1}^{1}z^2a_0\ud{z}\big)^2+\big(\int_{-1}^{1}z\sqrt{1-z^2}a_1\ud{z}\big)^2\\
&+\big(\int_{-1}^{1}z\sqrt{1-z^2}b_1\ud{z}\big)^2+\frac{1}{4}\big(\int_{-1}^{1}(1-z^2)a_2\ud{z}\big)^2+
\frac{1}{4}\big(\int_{-1}^{1}(1-z^2)b_2\ud{z}\big)^2\Big\}.
\end{align*}

We use Lemma \ref{lem:integral} and Cauchy-Schwartz
inequality to get
\begin{eqnarray*}
&&2\pi^2\Big(\int_{-1}^1\ue^{\eta{z^2}}\ud{z}\Big)\cdot\Big(\int_{-1}^1\ue^{-\eta{z^2}}a^2_1\ud{z}\Big)
-2\alpha\pi^2\big(\int_{-1}^{1}z\sqrt{1-z^2}a_1\ud{z}\big)^2\\
&&=2\alpha\pi^2\Big(\int_0^1\ue^{\eta{z^2}}z^2(1-z^2)\ud{z}\Big)\cdot\Big(\int_{-1}^1\ue^{-\eta{z^2}}a^2_1\ud{z}\Big)
-2\alpha\pi^2\Big(\int_{-1}^{1}z\sqrt{1-z^2}a_1\ud{z}\Big)^2\ge 0.
\end{eqnarray*}
Moreover, the equality holds if and only if
$a_1(z)=C\ue^{\eta{z^2}}z\sqrt{1-z^2}$ or $0$ for some constant $C$.
If $a_1(z)$ satisfies
$\int_{-1}^{1}\ue^{\eta{z^2}}z\sqrt{1-z^2}a_1(z)\ud{z}=0$, then
\begin{eqnarray}\label{eq:a-3.5}
&&\Big(\int_0^1\ue^{\eta{z^2}}z^2(1-z^2)\ud{z}\Big)\cdot\Big(\int_{-1}^1\ue^{-\eta{z^2}}a^2_1\ud{z}\Big)
-\Big(\int_{-1}^{1}z\sqrt{1-z^2}a_1\ud{z}\Big)^2\non\\
&&=\Big(\int_0^1\ue^{\eta{z^2}}z^2(1-z^2)\ud{z}\Big)\cdot\Big(\int_{-1}^1\ue^{-\eta{z^2}}a^2_1\ud{z}\Big)
-\Big(\int_{-1}^{1}z\sqrt{1-z^2}(1-\zeta\ue^{\eta{z^2}})a_1\ud{z}\Big)^2\non\\
&&\ge \Big(\int_0^1\ue^{\eta{z^2}}\big(z^2(1-z^2)-z^2(1-z^2)(1-\zeta\ue^{\eta{z^2}})^2\big)\ud{z}\Big)
\cdot\Big(\int_{-1}^1\ue^{-\eta{z^2}}a^2_1\ud{z}\Big)\non\\
&&=\Big(\int_0^1\ue^{2\eta{z^2}}z^2(1-z^2)\big(2\zeta-\zeta^2\ue^{\eta{z^2}}\big)\ud{z}\Big)
\cdot\Big(\int_{-1}^1\ue^{-\eta{z^2}}a^2_1\ud{z}\Big)\non\\
&&\ge c_0(\eta)\int_{-1}^1a^2_1\ud{z},
\end{eqnarray}
if we take $\zeta$ small enough. We denote \beno
W_1(\eta)=\int_{-1}^1\ue^{\eta{z^2}}(1-z^2)(5z^2-1)\ud{z}, \eeno
which is positive for $\eta>0$ by noting that \beno
W_1(\eta)=\int_{-1}^1(\ue^{\eta{z^2}}-\ue^{\eta\frac{1}{5}})(1-z^2)(5z^2-1)\ud{z}>0.
\eeno By Lemma \ref{lem:integral} and Cauchy-Schwartz inequality, we
get
\begin{eqnarray*}
&&2\pi^2\Big(\int_{-1}^1\ue^{\eta{z^2}}\ud{z}\Big)\cdot\Big(\int_{-1}^1\ue^{-\eta{z^2}}a^2_2\ud{z}\Big)
-\frac{1}{2}\alpha\pi^2\Big(\int_{-1}^{1}(1-z^2)a_2\ud{z}\Big)^2\\
&&=2\alpha\pi^2\Big(\int_0^1\ue^{\eta{z^2}}z^2(1-z^2)\ud{z}\Big)\cdot\Big(\int_{-1}^1\ue^{-\eta{z^2}}a^2_2\ud{z}\Big)
-\frac{1}{2}\alpha\pi^2\Big(\int_{-1}^{1}(1-z^2)a_2\ud{z}\Big)^2\\
&&=\frac{1}{2}\alpha\pi^2\Big(\int_0^1\ue^{\eta{z^2}}(1-z^2)^2\ud{z}\Big)\cdot\Big(\int_{-1}^1\ue^{-\eta{z^2}}a^2_2\ud{z}\Big)
-\frac{1}{2}\alpha\pi^2\Big(\int_{-1}^{1}(1-z^2)a_2\ud{z}\Big)^2\\
&&\quad+\frac{1}{2}\alpha\pi^2W_1(\eta)\int_{-1}^1\ue^{-\eta{z^2}}a^2_2\ud{z}\ge
c_0(\eta)\int_{-1}^1a^2_2\ud{z}
\end{eqnarray*}
for some $c_0(\eta)>0$.

In the following, we take $\eta=\eta_1(\al)$ for $\al>\al^*$. By
Lemma \ref{lem:integral}, we have \beno
0<\frac{\partial\alpha(\eta)}{\partial\eta}
=\frac{A_2(A_2-A_4)-A_0(A_4-A_6)}{(A_2-A_4)^2}=\frac{3A_2^2+2A_0A_2-5A_0A_4}{2\eta(A_2-A_4)^2},
\eeno which implies that \beno 3(A_0A_4-A_2^2)<2A_0(A_2-A_4). \eeno
Then using the fact $\int_{-1}^1a_0\ud{z}=0$ and Cauchy-Schwartz
inequality, we infer that
\begin{eqnarray*}
&&\Big(\int_{-1}^{1}(1-z^2)a_0\ud{z}\Big)^2+2\Big(\int_{-1}^{1}z^2a_0\ud{z}\Big)^2\\
&&=3\Big(\int_{-1}^{1}(\frac{A_2}{A_0}-z^2)a_0\ud{z}\Big)^2\\
&&\le
3\Big(\int_0^1\ue^{\eta{z^2}}(\frac{A_2}{A_0}-z^2)^2\ud{z}\Big)\cdot
\Big(\int_{-1}^1\ue^{-\eta{z^2}}a_0^2\ud{z}\Big)\\
&&=3\Big(A_4-\frac{A_2^2}{A_0}\Big)\int_{-1}^1\ue^{-\eta{z^2}}a_0^2\ud{z}\\
&&<2(A_2-A_4)\int_{-1}^1\ue^{-\eta{z^2}}a_0^2\ud{z}=\frac{2}{\alpha}
\Big(\int_{-1}^1\ue^{\eta{z^2}}\ud{z}\Big)\cdot
\Big(\int_{-1}^1\ue^{-\eta{z^2}}a_0^2\ud{z}\Big),
\end{eqnarray*}
which implies that
\begin{eqnarray*}
&&\frac{2}{\alpha} \Big(\int_{-1}^1\ue^{\eta{z^2}}\ud{z}\Big)\cdot
\Big(\int_{-1}^1\ue^{-\eta{z^2}}a_0^2\ud{z}\Big)
-\Big(\int_{-1}^{1}(1-z^2)a_0\ud{z}\Big)^2-2\Big(\int_{-1}^{1}z^2a_0\ud{z}\Big)^2\\
&&{\ge}c_0(\eta)\int_{-1}^1a_0^2\ud{z}\quad\textrm{ for some }
c_0(\eta)>0.
\end{eqnarray*}

Summing up all the above estimates, we conclude that if
$\eta=\eta_1(\al)$, then \ben\label{eq:3.5}
\big\langle\CH_{h_1}f, f\big\rangle=\big\langle\frac{f}{h_{\eta_1,\nn}}, f\big\rangle+\big\langle\CU{f},
f\big\rangle\ge c_0(\eta)\sum_{k\neq 1}\int_{-1}^1a_k^2\ud{z}, \een
and hence,
\begin{eqnarray*}
\big\langle\frac{f}{h_{\eta_1,\nn}}, f\big\rangle+\big\langle\CU{f},
f\big\rangle=0
\end{eqnarray*}
if and only if $f$ is of the form
$C_1z\sqrt{1-z^2}e^{\eta{z^2}}\cos\phi+C_2z\sqrt{1-z^2}e^{\eta{z^2}}\sin\phi$,
which belongs to $\text{ span}\big\{\CR_ih_{\eta,\nn},
i=1,2,3\big\}$, since we have
\begin{align*}
\CR h_{\eta,\nn}=&2\eta(\mm\times \nn)(\mm\cdot\nn)e^{\eta(\mm\cdot\nn)^2}
=2\eta\big(\sin\theta\sin\phi,-\sin\theta\cos\phi,0\big)\cos\theta e^{\eta\cos^2\theta}\\
=& 2\eta\big(\sin\phi,-\cos\phi,0\big)z\sqrt{1-z^2}e^{\eta z^2}.
\end{align*}
This proves the first statement of Proposition
\ref{prop:H-positive}. To obtain a lower bound of $\CH_{h_1}$, we
decompose $g$ into \beno f=f_1+f_2,\quad f_1\in \text{
span}\big\{\CR_ih_1, i=1,2,3\big\}. \eeno If $g$ satisfy
$\int_\BS f(\mm)\CA^{-1}_{h_1}\CR{h_1}\ud\mm=0$, then we have
\beno \int_\BS f_1(\mm)\CA^{-1}_{h_1}f_1\ud\mm=-\int_\BS
f_2(\mm)\CA^{-1}_{h_1}g_1\ud\mm, \eeno which implies that \beno
\langle{f_1}, f_1\rangle\le C\langle{f_2}, f_2\rangle, \eeno since
$\langle{f_1}, f_1\rangle\sim \big\langle{f_1},
\CA^{-1}_{h_1}f_1\big\rangle$.
This together with (\ref{eq:a-3.5}) and (\ref{eq:3.5}) gives \beno
\big\langle\CH_{h_1}{f}, f\big\rangle
=\big\langle\CH_{h_1}{f_2}, f_2\big\rangle
\ge{c_0}\langle{f_2}, f_2\rangle\ge c_0\langle{f}, f\rangle. \eeno
The proof is finished.\ef \vspace{0.1cm}

We define $\mathrm{Ker}\CG_{h_{\eta,\nn}}\eqdefa \big\{\phi\in
H_0^2(\BS): \CG_{h_{\eta,\nn}}\phi=0\big\}$.

\begin{Theorem}\label{thm:G-kernel}
Let $h_i=h_{\eta_i,\nn}, i=1,2$. For $\al>\alpha^*$, it holds that
\begin{itemize}

\item[1.] $\CG_{h_1}$ has no positive eigenvalues, while $\CG_{h_2}$ has at least one positive eigenvalue;

\item[2.] $\mathrm{Ker}\CG_{h_1}=\big\{\mathbf{{\Theta}}\cdot
\CR{h_1}; \mathbf{{\Theta}}\in \RR^3\big\}$ is a two dimensional
space;

\item[3.] If $\phi\in \mathrm{Ker}\CG_{h_1}$, then $\CH_{h_1}\phi=0$;

\end{itemize}

\end{Theorem}

\no{\bf Proof.}\,Since $h_1$ is a stable critical point
of $A[f]$, $\CG_{h_1}$ has no positive eigenvalues by
Proposition \ref{prop:G-symm}. From the proof of Proposition
\ref{prop:H-positive}, we know that there exists $g\in\CP_0(\BS)$
such that \ben\label{eq:H-3.6}
\big\langle\CH_{h_2}g,g\big\rangle<0. \een Assume that
all eigenvalues $\{\lambda_k\}$ of $\CG_{h_2}$ are
non-positive. We denote by $E_k$ the eigen-subspaces of
$\CG_{h_2}$ corresponding to $\lambda_k$. Then for
$\psi_k\in E_k, \psi_\ell \in E_\ell (k\neq{\ell})$,  we have \beno
\lambda_k\big\langle\psi_k, \CA^{-1}_{h_2}\psi_\ell\big\rangle=
\big\langle\CG\psi_k,\CA^{-1}_{h_2}\psi_\ell\big\rangle=\big\langle\CG\psi_\ell,\CA^{-1}_{h_2}\psi_k\big\rangle
=\lambda_\ell\big\langle\psi_k, \CA^{-1}_{h_2}\psi_\ell\big\rangle. \eeno
Hence, $\langle\psi_k,~\CA^{-1}_{h_2}\psi_\ell\rangle=0$ for
$k\neq{\ell}$. We write $g=\sum_k\psi_k$ with $\psi_k\in E_k$. Then
\beno \big\langle\CH_{h_2}g,g\big\rangle=-\big\langle\CG{g},
\CA^{-1}_{h_2}g\big\rangle=-\sum_k\lambda_k\big\langle\psi_k,
\CA^{-1}_{h_2}\psi_k\big\rangle\ge0, \eeno which leads to a contradiction
with (\ref{eq:H-3.6}). Thus, $\CG_{h_2}$ has at least one
positive eigenvalue.

If $\phi\in \mathrm{Ker}\CG_{h_1}$, then
$\CH_{h_1}\phi=\textrm{constant}$. Hence,
$\langle\CH_{h_1}\phi, \phi\rangle=0$, and then $\phi\in
\mathrm{span}\big\{\CR_i{h_1},i=1,2,3\big\}$ by
Proposition \ref{prop:H-positive}. On the other hand, if $\phi=\CR
h_1$, then we find
\beno \CH_{h_1}\phi=\CR\big(\ln h_1+\CU h_1\big)=0.
\eeno This proves
$\mathrm{Ker}\CG_{h_1}=\big\{\mathbf{{\Theta}}\cdot
\CR{h_1}; \mathbf{{\Theta}}\in \RR^3\big\}$ and the third point. Due
to $\nn\cdot\CR{h_1}=0$,
$\mathrm{Ker}\CG_{h_1}$ is a two dimensional space. \ef
\vspace{0.1cm}

Finally let us give a characterization of the functions in
$\mathrm{Ker} \CG_{h_{\eta_1,\nn}}^*$, see also \cite{KD}.

\begin{Proposition}\label{prop:kernel of Gstar}
If $\psi_0\in \mathrm{Ker}\CG_{h_{\eta_1,\nn}}^*$, then $\psi_0$
takes the form $(\theta, \phi)$ \beno
\psi_0(\theta,\phi)=\mathbf\Theta\cdot\ee_\phi{g_0}(\cos\theta),
\eeno in the spherical coordinate, where $\ee_\phi=(-\sin\phi,
\cos\phi, 0)$ and $g_0$ satisfies
\begin{equation}\label{eq:g0}
\frac{1}{\sin\theta}\frac{\ud}{\ud\theta}\Big(\sin\theta\frac{\ud{g_0}}{\ud\theta}\Big)
-\frac{g_0}{\sin^2\theta}
-\frac{\ud{u_0}}{\ud\theta}\frac{\ud{g_0}}{\ud\theta}=-\frac{\ud{u_0}}{\ud\theta}.
\end{equation}
\end{Proposition}

\no{\bf Proof.}\,Note that
$\mathrm{Ker}\CG_{h_{\eta_1,\nn}}^*=\CA_{h_{\eta_1,\nn}}^{-1}\mathrm{Ker}\CG_{h_{\eta_1,\nn}}$.
Hence, $\psi_0\in \mathrm{Ker}\CG_{h_{\eta_1,\nn}}^*$ if and only if
there exits a vector $\mathbf{\Theta}$ such that \beno
\CR\cdot(h_{h_{\eta_1,\nn}}\CR{\psi_0})=\mathbf{\Theta}\cdot\CR{h_{h_{\eta_1,\nn}}},
\eeno which is equivalent to
\begin{equation}\label{eq:3.7}
\CR\cdot\CR{\psi_0}-\CR{u_0}\cdot\CR\psi_0=-\mathbf{\Theta}\cdot\CR{u_0},
\end{equation}
where $u_0=\CU{h_{\eta_1,\nn}}$ is a function of $\mm\cdot\nn$. We
take $\theta$ be the angle between $\mm$ and $\nn$, and rewrite
(\ref{eq:3.7}) in terms of in the spherical coordinate $(\theta,
\phi)$ as
\beno
\frac{1}{\sin\theta}\frac{\partial}{\partial\theta}\Big(\sin\theta\frac{\partial\psi_0}{\partial\theta}\Big)
+\frac{1}{\sin^2\theta}\frac{\partial^2{\psi_0}}{\partial\phi^2}
-\frac{\ud{u_0}}{\ud\theta}\frac{\partial\psi_0}{\partial\theta}=-\mathbf\Theta\cdot\ee_\phi\frac{\ud{u_0}}{\ud\theta}.
\eeno We rewrite $\psi_0(\theta,\phi)$ as \beno
\psi_0(\theta,\phi)=\mathbf\Theta\cdot\ee_\phi{g_0}. \eeno Then it
easy to find that $g_0$ satisfies (\ref{eq:g0}).\ef

\section{Lower bound of a bilinear form for the linearized operator}

In the inhomogeneous case, the linearized operator $\CG^\ve_h$ around $h$ is given by
\beno \CG_h^\ve f=\CR\cdot\big(\CR{f}+h\CR
\CU_\ve f+f\CR\CU h\big).
\eeno
To justify the small Deborah number limit for the inhomogeneous system,
the main difficulty is that the elastic stress in the
velocity equation is strongly singular(a loss of $\f 1 \ve$). To
overcome it, we have to establish a precise lower bound for the following
bilinear form:
\beno \big\langle
\CG_{h}^\ve f,\CH^\ve_h f\big\rangle,\quad \CH^\ve_h f=\f f {h}+\CU_\ve
f.
\eeno
When $\ve \neq 0$, the orthogonal structure is destroyed
such that the interactions between the part inside the kernel and
the part outside the kernel of $\CG_h$ become very complicated. We find a coordinate transformation and introduce a
generalized kernel space of $\CG_h^\ve$ such that the interactions
between two parts can been seen explicitly, then a lower bound is
obtained by very subtle calculations.

\subsection{New coordinates frame}
At each point $\xx$, we choose a right hand cartesian coordinate
frame $(\kk_1(\xx), \kk_2(\xx), \kk_3(\xx))$ such that
$\kk_3(\xx)=\nn(\xx)$, and $\kk_1(\xx)$, $\kk_2(\xx)$ depend on
$\nn(\xx)$ smoothly. For instance, under the assumption that
$|n_1(\xx)|<1-c_0$ for all $\xx\in\Omega$, we can take
\begin{align*}
&\kk_2(\xx)=\frac{\nn\times(1,0,0)^T}{|\nn\times(1,0,0)^T|}=\Big(0,
\frac{n_3}{(n_2^2+n_3^2)^{1/2}},-\frac{n_2}{(n_2^2+n_3^2)^{1/2}}\Big)^T, \\
&\kk_1=\kk_2\times\nn=\Big((n_2^2+n_3^2)^{1/2},-\frac{n_1n_2}{(n_2^2+n_3^2)^{1/2}},
-\frac{n_1n_3}{(n_2^2+n_3^2)^{1/2}}\Big)^T.
\end{align*}
At each point $\xx$, let $(\htheta, \hvphi)$ be the sphere
coordinate on the unit sphere $\BS$, that is,
\begin{align*}
\mm&=\sin\htheta\cos\hvphi\kk_1(\xx)+\sin\htheta\sin\hvphi\kk_2(\xx)+\cos\htheta\kk_3(\xx)\\
&=\mathbf{A}(\xx)\cdot(\sin\htheta\cos\hvphi,\sin\htheta\sin\hvphi,
\cos\htheta)^T,
\end{align*}
where the matrix $\mathbf{A}=[\kk_1~\kk_2~\kk_3]$. We set \beno
\ee_{\hvphi}=-\sin\hvphi\kk_1+\cos\hvphi\kk_2,\quad
\ee_\htheta=-(\cos\htheta\cos\hvphi\kk_1
+\cos\htheta\sin\hvphi\kk_2-\sin\htheta\kk_3). \eeno We denote
$\hat{\mm}=(\sin\htheta\cos\hvphi,\sin\htheta\sin\hvphi,
\cos\htheta)^T$, hence $\mm=\BA(\xx)\cdot\hat{\mm}$.

In this coordinate, the rotational gradient operator
$\CR=\mm\times\nabla_{\mm}$ can be written as
\begin{align}\label{eq:R operator}
\CR&=\big(-\sin\hvphi\kk_1+\cos\hvphi\kk_2)\frac{\partial}{\partial\htheta}
-(\cos\htheta\cos\hvphi\kk_1+\cos\htheta\sin\hvphi\kk_2-\sin\htheta\kk_3\big)
\frac{1}{\sin\htheta}\frac{\partial}{\partial\hvphi}\non\\
&=\BA\cdot\big(-\sin\hvphi\frac{\partial}{\partial\htheta}-\cos\hvphi\frac{\cos\htheta}{\sin\htheta}
\frac{\partial}{\partial\hvphi},\cos\hvphi\frac{\partial}{\partial\htheta}
-\sin\hvphi\frac{\cos\htheta}{\sin\htheta}\frac{\partial}{\partial\hvphi}, \frac{\partial}{\partial\hvphi}\big)^T.
\end{align}
We also have
\begin{align}\label{eq:R operator-p}
\CR{f}\cdot\CR{g}&=\big(\ee_\hvphi\frac{\partial{f}}{\partial\htheta}+\ee_\htheta\frac{1}{\sin\htheta}
\frac{\partial{f}}{\partial\hvphi}\big)\cdot\big(\ee_\hvphi\frac{\partial{g}}{\partial\htheta}
+\ee_\htheta\frac{1}{\sin\htheta}\frac{\partial{g}}{\partial\hvphi}\big)\non\\
&=\frac{\partial{f}}{\partial\htheta}\frac{\partial{g}}{\partial\htheta}
+\frac{1}{\sin^2\htheta}\frac{\partial{f}}{\partial\hvphi}\frac{\partial{g}}{\partial\hvphi}.
\end{align}

\subsection{The Maier-Saupe space and the lower bound inequality}
For any $f\in L^2(\Om\times \BS)$ with $\int_{\BS}f(\xx,\mm)d\mm=0$,
we decompose it as
\begin{align}\label{eq:f-decompo}
f=a_0(\xx,\htheta)+\sum_{k\ge1}\big(a_k(\xx,\htheta)\cos{k}\hvphi+b_k(\xx,\htheta)\sin{k}\hvphi\big),
\end{align}
with $a_k(\xx,0)=a_k(\xx,\pi)=b_k(\xx,0)=b_k(\xx,\pi)=0$ for $k\ge
1$. We set
\beno
A_k=\int_{0}^\pi\ue^{\eta\cos^2\theta}(\cos\theta)^k\sin\theta\ud\theta,\quad
f_0(\htheta)=\frac{\ue^{\eta\cos^2\htheta}}{2\pi{A}_0}. \eeno We
further decompose the coefficients $a_0, a_1, a_2, b_1, b_2$ as
follows \beno
&&a_0(\xx,\htheta)=\zeta_{0}(\xx)f_0(\htheta)(\cos^2\htheta-A_2/A_0)+\gamma_0(\xx,\htheta),\\
&&a_1(\xx,\htheta)=\zeta_{a,1}(\xx)f_0(\htheta)\sin\htheta\cos\htheta+\gamma_{a,1}(\xx,\htheta),\\
&&b_1(\xx,\htheta)=\zeta_{b,1}(\xx)f_0(\htheta)\sin\htheta\cos\htheta+\gamma_{b,1}(\xx,\htheta),\\
&&a_2(\xx,\htheta)=\zeta_{a,2}(\xx)f_0(\htheta)\sin^2\htheta+\gamma_{a,2}(\xx,\htheta),\\
&&b_2(\xx,\htheta)=\zeta_{b,2}(\xx)f_0(\htheta)\sin^2\htheta+\gamma_{b,2}(\xx,\htheta),
\eeno
where the functions $\ga_0, \ga_{a,1},\cdots,\ga_{b,2}$ satisfy
\beno
&&\int_0^\pi\gamma_0(\xx,\htheta)\sin\htheta\ud\htheta=0, \quad \int_0^\pi(3\cos^2\htheta-1)\gamma_0(\xx,\htheta)\sin\htheta\ud\htheta=0,\\
&&\int_0^\pi\sin\htheta\cos\htheta\gamma_{a,1}(\xx,\htheta)\sin\htheta\ud\htheta=0,\quad
\int_0^\pi\sin\htheta\cos\htheta\gamma_{b,1}(\xx,\htheta)\sin\htheta\ud\htheta=0,\\
&&\int_0^\pi\sin^2\htheta\gamma_{a,2}(\xx,\htheta)\sin\htheta\ud\htheta=0,\quad
\int_0^\pi\sin^2\htheta\gamma_{b,2}(\xx,\htheta)\sin\htheta\ud\htheta=0.
\eeno Noting that
\begin{align*}
\int_0^\pi{f}_0(\cos^2\htheta-A_2/A_0)(3\cos^2\htheta-1)\sin\htheta\ud\htheta=
\frac{3}{A_0}\big(A_4-\frac{A_2^2}{A_0}\big)>0,
\end{align*}
hence $\zeta_0$ is uniquely determined. Obviously,
$\zeta_{a,1},\cdots,\zeta_{b,2}$ are also uniquely determined. Thus,
the above decompositions make sense.

The space spanned by the following five functions \beno
&f_0(\htheta)(\cos^2\htheta-A_2/A_0),\quad
f_0(\htheta)\sin\htheta\cos\htheta\cos\hvphi,\quad
f_0(\htheta)\sin\htheta\cos\htheta\sin\hvphi,\\
&f_0(\htheta)\sin^2\htheta\cos(2\hvphi),\quad \textrm{and} \quad
f_0(\htheta)\sin^2\htheta\sin(2\hvphi) \eeno can be viewed as a
generalized kernel of the operator $\CG^\ve_h$, and will be called
as \textbf{the Maier-Saupe space}. The kernel of $\CG_h$ is spanned  by the second and the third
function.

We denote
\beno
&&\HMM(\xx)=\int_\BS\hmm\hmm{f}\ud\hmm,\quad \MM(\xx)=\int_\BS\mm\mm{f}\ud\mm=\int_\BS(\BA\cdot\hmm)(\BA\cdot\hmm)f\ud\hmm,\\
&&N_{ij}(\xx)=A_{ki}A_{lj} \big({g}_\ve*M_{kl}\big),\quad
N_0(\xx)=2N_{33}-N_{11}-N_{22},\quad N_2(\xx)=N_{11}-N_{22}. \eeno

\begin{Proposition}(Lower bound inequality)\label{prop:lower bound}
Let $h_{\nn(x),\eta}$ be a  stable critical point of $A[f]$. Then
there exits $c>0$ such that any $f\in H^1(\Om\times \BS)$ with
$\int_{\BS}f(\xx,\mm)d\mm=0$, there holds
\beno &&\big\langle
\CG_{h_{\nn,\eta}}^\ve f,\CH^\ve_{h_{\nn,\eta}} f\big\rangle\ge
c\Big\{\int_\Omega\int_0^\pi\big(\gamma_0^2+\gamma_{a,1}^2+\gamma_{b,1}^2
+\gamma_{a,2}^2+\gamma_{b,2}^2+\sum_{k\ge3}(a_k^2+b_k^2)\big)\sin\htheta\ud{\htheta}\ud\xx\\
&&+\int_\Omega\big((2{\zeta}_{0}-\alpha{N}_0)^2+\big({\zeta}_{a,1}-2\alpha{N_{13}}\big)^2
+\big({\zeta}_{b,1}-2\alpha{N_{23}}\big)^2+\big({\zeta}_{a,2}-\alpha{N_{2}}\big)^2
+\big({\zeta}_{b,2}-\alpha{N_{12}}\big)^2\big)\ud\xx\Big\}.
\eeno
\end{Proposition}
\begin{Remark}
This inequality gives a good bound for the part outside the
Maier-Saupe space, and a weak bound for the part inside the
Maier-Saupe space.
\end{Remark}

\begin{Remark}
By letting $\ve$ tend to zero, the lower bound inequality implies that
\beno
\big\langle\CG_{h_{\nn,\eta}} f,\CH_0 f\big\rangle\ge c\big\langle\CH_0 f, f\big\rangle,\quad \CH_0=\CH_{h_{\nn,\eta}},
\eeno
which can also be deduced the following simple argument. Noting that $\CH(\mathrm{Ker}\CG_{h_{\nn,\eta}})=0$,
we may assume that $f\in(\mathrm{Ker}\CG^*_{h_{\nn,\eta}})^\bot$.
Then by Poincar\'{e} inequality and Proposition \ref{prop:H-positive}, we have
\begin{eqnarray*}
\big\langle\CH_0{f}, \CA\CH_0{f}\big\rangle\ge\big\langle\CR\CH_0{f},\CR\CH_0{f}\big\rangle
\ge\big\langle\CH_0{f}-\overline{\CH_0{f}}, \CH_0{f}-\overline{\CH_0{f}}\big\rangle\\
\ge\frac{\big\langle\CH_0{f}-\overline{\CH_0{f}},f\big\rangle^2}{\langle{f},f\rangle}
\ge{c}\big\langle\CH_0{f},f\big\rangle\ge{c}\langle{f},f\rangle,
\end{eqnarray*}
where $\overline{\CH_0{f}}=\frac{1}{4\pi}\int_\BS\CH_0{f}\ud\mm$.
\end{Remark}

\subsection{Proof of the lower bound inequality}
First of all, we can get by direct calculations that
\begin{align*}
\CU_\ve{f}&=\alpha\int_\Omega\int_\BS{g}_\ve(\xx-\xx')(1-\mm\mm:\mm'\mm')f(\xx',\mm')\ud\mm'\ud\xx'\\
&=-\alpha\mm\mm:\int_\Omega
{g}_\ve(\xx-\xx')\MM(\xx')\ud\xx'\\
&=-\alpha(\BA\cdot\hmm)(\BA\cdot\hmm):\int_\Omega
{g}_\ve(\xx-\xx')\MM(\xx')\ud\xx'\\
&=-\alpha\hm_i\hm_jA_{ki}A_{lj}
\big({g}_\ve*(A_{ki'}A_{lj'}\HM_{i'j'})\big).
\end{align*}
For the simplicity, we denote $f_0=h_{\eta,\nn}$ in what follows.
Using (\ref{eq:R operator})-(\ref{eq:f-decompo}), we get by very
tedious calculations of competing the square that
\begin{align*}
&\big\langle \CG_{f_0}^\ve f,\CH^\ve_{f_0} f\big\rangle
=\big\langle{f_0}\CR(\frac{f}{f_0}+\CU_\ve{f}),\CR(\frac{f}{f_0}+\CU_\ve{f})\big\rangle\\
&=\big\langle{f}_0\CR(\frac{f}{f_0}-\alpha\hm_i\hm_jN_{ij}), \CR(\frac{f}{f_0}-\alpha\hm_i\hm_jN_{ij})\big\rangle\\
&=\int_\Omega\int_\BS{f}_0\big|\partial_\htheta(\frac{f}{f_0}-\alpha\hm_i\hm_jN_{ij})\big|^2\ud\hmm\ud\xx
+\int_\Omega\int_\BS\frac{{f}_0}{\sin^2\htheta}\big|\partial_\hvphi
(\frac{f}{f_0}-\alpha\hm_i\hm_jN_{ij})\big|^2\ud\hmm\ud\xx\\
&=\int_\Omega\int_\BS{f}_0\big|\partial_\htheta(\frac{a_0}{f_0})
+\alpha(2N_{33}-N_{11}-N_{22})\sin\htheta\cos\htheta\big|^2\ud\hmm\ud\xx\\
&\quad+\frac{1}{2}\int_\Omega\int_\BS{f}_0\big|\partial_\htheta(\frac{a_1}{f_0})
-2\alpha{N}_{13}\cos2\htheta\big|^2+\frac{{f}_0}{\sin^2\htheta}\big(\frac{a_1}{f_0}
-2\alpha{N}_{13}\sin\htheta\cos\htheta\big)^2\ud\hmm\ud\xx\\
&\quad+\frac{1}{2}\int_\Omega\int_\BS{f}_0\big|\partial_\htheta(\frac{b_1}{f_0})
-2\alpha{N}_{23}\cos2\htheta\big|^2+\frac{{f}_0}{\sin^2\htheta}\big(\frac{b_1}{f_0}
-2\alpha{N}_{23}\sin\htheta\cos\htheta\big)^2\ud\hmm\ud\xx\\
&\quad+\frac{1}{2}\int_\Omega\int_\BS{f}_0\big|\partial_\htheta(\frac{a_2}{f_0})
-\alpha({N}_{11}-N_{22})\sin\htheta\cos\htheta\big|^2+
\frac{{f}_0}{\sin^2\htheta}\big(\frac{2a_2}{f_0}-\alpha({N}_{11}-N_{22})\sin^2\htheta\big)^2\ud\hmm\ud\xx\\
&\quad+\frac{1}{2}\int_\Omega\int_\BS{f}_0\big|\partial_\htheta(\frac{b_2}{f_0})
-2\alpha{N}_{12}\sin\htheta\cos\htheta\big|^2+\frac{4{f}_0}{\sin^2\htheta}
\big(\frac{b_2}{f_0}-\alpha{N}_{12}\sin^2\htheta\big)^2\ud\hmm\ud\xx\\
&\quad+\frac{1}{2}\sum_{k\ge3}\int_\Omega\int_\BS{f}_0\big|\partial_\htheta(\frac{a_k}{f_0})\big|^2
+{f}_0\big|\partial_\htheta(\frac{b_k}{f_0})\big|^2+
\frac{k^2}{\sin^2\htheta}\frac{a_k^2+b^2_k}{f_0}\ud\hmm\ud\xx.
\end{align*}
Making a change of variable $z=\cos\htheta$, we obtain
\begin{align*}
&\big\langle \CG_{f_0}^\ve f,\CH^\ve_{f_0} f\big\rangle\\
&=\int_\Omega\int_{-1}^1{f}_0{(1-z^2)}\big(\partial_z(\frac{a_0}{f_0})
-\alpha(2N_{33}-N_{11}-N_{22})z\big)^2\ud{z}\ud\xx\\
&\quad+\frac{1}{2}\int_\Omega\int_{-1}^1{f}_0\big|\partial_z(\frac{a_1}{f_0})\sqrt{1-z^2}
+2\alpha{N}_{13}(2z^2-1)\big|^2+\frac{{f}_0}{1-z^2}\big(\frac{a_1}{f_0}-2\alpha{N}_{13}z\sqrt{1-z^2}\big)^2\ud{z}\ud\xx\\
&\quad+\frac{1}{2}\int_\Omega\int_{-1}^1{f}_0\big|\partial_z(\frac{b_1}{f_0})\sqrt{1-z^2}
+2\alpha{N}_{23}(2z^2-1)\big|^2+\frac{{f}_0}{1-z^2}\big(\frac{b_1}{f_0}-2\alpha{N}_{23}z\sqrt{1-z^2}\big)^2\ud{z}\ud\xx\\
&\quad+\frac{1}{2}\int_\Omega\int_{-1}^1{f}_0{(1-z^2)}\big(\partial_z(\frac{a_2}{f_0})
+\alpha({N}_{11}-N_{22})z\big)^2+
\frac{{f}_0}{1-z^2}\big(\frac{2a_2}{f_0}-\alpha({N}_{11}-N_{22})(1-z^2)\big)^2\ud{z}\ud\xx\\
&\quad+\frac{1}{2}\int_\Omega\int_{-1}^1{f}_0{(1-z^2)}\big(\partial_z(\frac{b_2}{f_0})
+2\alpha{N}_{12}z\big)^2+\frac{4{f}_0}{1-z^2}\big(\frac{b_2}{f_0}-\alpha{N}_{12}(1-z^2)\big)^2\ud{z}\ud\xx\\
&\quad+\frac{1}{2}\sum_{k\ge3}\int_\Omega\int_\BS{f}_0\big|\partial_\htheta(\frac{a_k}{f_0})\big|^2
+{f}_0\big|\partial_\htheta(\frac{b_k}{f_0})\big|^2+
\frac{k^2}{\sin^2\htheta}\frac{a_k^2+b^2_k}{f_0}\ud\hmm\ud\xx .
\end{align*}

\no{\bf$\bullet$ Lower bound for the term including
$a_0$}\vspace{0.1cm}

To deal with the cross term, we need to introduce a slightly different decomposition
\begin{align*}
a_0(\xx,z)=\hat{\zeta}_{0}(\xx)f_0(z)\big(z^2-\frac{A_2}{A_0}\big)+\hat{\gamma}_0(\xx,z),
\end{align*}
where $\hat{\gamma}_0$ satisfies
$\int_{-1}^1\hat{\gamma}_0(\xx,z)\big(3z^2-1-2\eta{z}^2(1-z^2)\big)\ud{z}=0$.
From the fact that
\begin{align*}
&\int_{-1}^1f_0(A_0z^2-A_2)\big(3z^2-1-2\eta{z}^2(1-z^2)\big)\ud{z}\\
&=3(A_0A_4-A_2^2)+2\eta(A_2(A_2-A_4)-A_0(A_4-A_6))>0,
\end{align*}
we know that $\hat{\zeta}_{0}$ are uniquely determined. Then we have
\begin{align*}
&\int_\Omega\int_{-1}^1{f}_0(1-z^2)\big|\partial_z(\frac{a_0}{f_0})
-\alpha{N}_0(\xx)z\big|^2\ud{z}\ud\xx\\
&=\int_\Omega\int_{-1}^1{f}_0(1-z^2)\big(\partial_z(\frac{\hat{\gamma}_0(\xx,z)}{f_0})\big)^2\ud{z}\ud\xx
+\Big(\int_\Omega\big(2\hat{\zeta}_{0}-\alpha{N_{0}}\big)^2\ud\xx\Big)\Big(\int_{-1}^1(1-z^2)z^2f_0\ud{z}\Big)\\
&\quad-2\int_\Omega\big(2A_0\hat{\zeta}_{0}(\xx)-\alpha{N}_{0}(\xx)\big)\int_{-1}^1
\partial_z(\frac{\hat{\gamma}_0(\xx,z)}{f_0})f_0(z)z{(1-z^2)}\ud{z}\ud\xx\\
&=\int_\Omega\int_{-1}^1{f}_0(1-z^2)\big(\partial_z(\frac{\hat{\gamma}_0(\xx,z)}{f_0})\big)^2\ud{z}\ud\xx
+\Big(\int_\Omega\big(2\hat{\zeta}_{0}-\alpha{N_{0}}\big)^2\ud\xx\Big)\Big(\int_{-1}^1(1-z^2)z^2f_0\ud{z}\Big)\\
&\quad+2\int_\Omega\big(2A_0\hat{\zeta}_{0}(\xx)-\alpha{N}_{0}(\xx)\big)\int_{-1}^1
\hat{\gamma}_0(\xx,z)\big(1-3z^2+2\eta{z^2}{(1-z^2)}\big)\ud{z}\ud\xx\\
&=\int_\Omega\int_{-1}^1{f}_0(1-z^2)\big(\partial_z(\frac{\hat{\gamma}_0(\xx,z)}{f_0})\big)^2\ud{z}\ud\xx
+\Big(\int_\Omega\big(2\hat{\zeta}_{0}-\alpha{N_{0}}\big)^2\ud\xx\Big)\Big(\int_{-1}^1(1-z^2)z^2f_0\ud{z}\Big).
\end{align*}
Recall that we have decomposed $a_0(x,z)$ as
\begin{align*}
a_0(\xx,z)=\zeta_{0}(\xx)f_0(z)(z^2-\frac{A_2}{A_0})+\gamma_0(\xx,z),
\end{align*}
thus
$(\zeta_0-\hat{\zeta}_0)f_0(z)(z^2-\frac{A_2}{A_0})=\hat{\gamma}_0-\gamma_0$,
which implies that \beno
\zeta_{0}(\xx)-\hat{\zeta}_{0}(\xx)=\int_{-1}^1\hat{\gamma}_0(3z^2-1)\ud{z}/(3A_0A_4-3A_2^2).
\eeno Thus we have
\begin{align*}
\int_\Omega\int_{-1}^1{f}_0(1-z^2)\big(\partial_z(\frac{{\gamma}_0}{f_0})\big)^2\ud{z}\ud\xx
&=\int_\Omega\int_{-1}^1{f}_0(1-z^2)\Big(\partial_z(\frac{\hat{\gamma}_0}{f_0})-2A_0(\zeta_{0}
-\hat{\zeta}_{0})z\Big)^2\ud{z}\ud\xx\\
&\le{C}\int_\Omega\int_{-1}^1{f}_0(1-z^2)\Big(\partial_z(\frac{\hat{\gamma}_0}{f_0})\Big)^2
+{\hat{\gamma}_0^2}\ud{z}\ud\xx.
\end{align*}
Then we infer that
\begin{align}\label{eq:a0}
&\int_\Omega\int_{-1}^1{f}_0(1-z^2)\big|\partial_z(\frac{a_0}{f_0})
-\alpha{N}_0(\xx)z\big|^2\ud{z}\ud\xx\nonumber\\
&\ge{c}\Big\{\int_\Omega\int_{-1}^1{f}_0(1-z^2)\big(\partial_z(\frac{{\gamma}_0(\xx,z)}{f_0})\big)^2\ud{z}\ud\xx
+\int_\Omega\big(2{\zeta}_{0}-\alpha{N_{0}}\big)^2\ud\xx\int_{-1}^1(1-z^2)z^2f_0\ud{z}\Big\}\non\\
&\ge{c}\Big\{\int_\Omega\int_{-1}^1{\gamma}_0^2+(\partial_z{\gamma}_0)^2\ud{z}\ud\xx
+\int_\Omega\big(2{\zeta}_{0}(\xx)-\alpha{N_{0}}\big)^2\ud\xx\Big\},
\end{align}
where we used the following Poincar\'{e} type
inequality in the last inequality:
\begin{Lemma}
There exists $c>0$ such that if $\ga(z)$ satisfies \beno
\int_{-1}^1{\gamma}\ud{z}=0,\quad
\int_{-1}^1(3z^2-1){\gamma}\ud{z}=0, \eeno then there holds
$$\int_{-1}^1{f}_0(1-z^2)\big(\partial_z(\frac{{\gamma}}{f_0})\big)^2\ud{z}\ge{c}\int_{-1}^1{\gamma}^2\ud{z}.$$
\end{Lemma}

\no{\bf Proof.}\, We define
$\bar{\gamma}(\mm)={\gamma}(\mm\cdot\nn)$. By the assumption, we
know \beno \int_\BS\mm\mm\bar{\gamma}(\mm)\ud\mm=0,\quad
\CH_0\bar{\gamma}=\frac{\bar\gamma}{f_0}. \eeno Hence, $\bar\ga\in
(\textrm{Ker} \CG_{f_0})^\perp$ and we have
$$\int_\BS{f_0}|\CR\CH_{f_0}\bar\gamma|^2\ud\mm
=2\pi\int_{-1}^1{f}_0(1-z^2)\big(\partial_z(\frac{{\gamma}}{f_0})\big)^2\ud{z}.$$
Set $\bar{C}=\frac{1}{4\pi}\int\CH_{f_0}\bar\gamma\ud\mm$. It follows
from Poincar\'{e} inequality and Proposition \ref{prop:H-positive}
that
\begin{align*}
\int_\BS{f_0}|\CR\CH_{f_0}\bar\gamma|^2\ud\mm&\ge{c}\int_\BS(\CH_{f_0}\bar\gamma-\bar{C})^2\ud\mm
\ge{c}\frac{\big(\int_\BS(\CH_{f_0}\bar\gamma-\bar{C})\bar\gamma\ud\mm\big)^2}{\int_\BS\bar\gamma^2\ud\mm}\\
&\ge{c}\frac{\big(\int_\BS(\CH_{f_0}\bar\gamma)\bar\gamma\ud\mm\big)^2}{\int_\BS\bar\gamma^2\ud\mm}
\ge{c}\int_\BS\bar\gamma^2\ud\mm=c4\pi^2\int_{-1}^1{\gamma}^2\ud{z},
\end{align*}
which completes the proof.\ef\vspace{0.1cm}

\no{\bf$\bullet$ Lower bound for the terms including $a_1, b_1$}\vspace{0.1cm}

Recall that we have a decomposition for $a_1(x,z)$ as
\begin{align*}
a_1(\xx,z)={\zeta}_{a,1}(\xx)f_0(z)z\sqrt{1-z^2}+{\gamma}_{a,1}(\xx,z),
\end{align*}
where $\int_{-1}^1z\sqrt{1-z^2}{\gamma}_{a,1}\ud{z}=0$. Thus, we can
get
\begin{align}
&\int_\Omega\int_{-1}^1\frac{{f}_0}{1-z^2}\big(\frac{a_1}{f_0}-2\alpha{N}_{13}z\sqrt{1-z^2}\big)^2\ud{z}\ud\xx\nonumber\\
&\ge\int_\Omega\int_{-1}^1{{f}_0}\Big(\frac{{\gamma}_{a,1}(\xx,z)}{f_0}
+\big({\zeta}_{a,1}(\xx)-2\alpha{N}_{13}(\xx)\big)z\sqrt{1-z^2}\Big)^2\ud{z}\ud\xx\nonumber\\\label{eq:a1}
&=\int_\Omega\int_{-1}^1\frac{{\gamma}_{a,1}^2(\xx,z)}{f_0}\ud{x}\ud\xx
+\int_\Omega\big({\zeta}_{a,1}(\xx)-2\alpha{N}_{13}(\xx)\big)^2\ud\xx\int_{-1}^1f_0z^2(1-z^2)\ud{z}.
\end{align}
A lower bound for the terms including $b_1$ can be obtained in the same way. \vspace{0.1cm}

\no{\bf$\bullet$ Lower bound for the terms including $a_2, b_2$}\vspace{0.1cm}

We have decomposed $a_2(x,z)$ as
\begin{align*}
a_2(\xx,z)={\zeta}_{a,2}(\xx)f_0(z)(1-z^2)+{\gamma}_{a,2}(\xx,z),
\end{align*}
where $\int_{-1}^1{\gamma}_{a,2}(\xx,z)(1-z^2)\ud{z}=0$. Then we
have
\begin{align}
&\int_\Omega\int_{-1}^1\frac{4{f}_0}{1-z^2}\big(\frac{a_2}{f_0}-\alpha{N}_{2}(1-z^2)\big)^2\ud{z}\ud\xx\nonumber\\
&\ge\int_\Omega\int_{-1}^1{4{f}_0}\Big(\frac{{\gamma}_{a,2}(\xx,z)}{f_0}
+\big({\zeta}_{a,2}(\xx)-\alpha{N}_{2}(\xx)\big)(1-z^2)\Big)^2\ud{z}\ud\xx\nonumber\\\label{eq:a2}
&=\int_\Omega\int_{-1}^1\frac{{\gamma}_{a,2}^2(\xx,z)}{f_0}\ud{z}\ud\xx
+\int_\Omega\big({\zeta}_{a,2}(\xx)-\alpha{N}_{2}(\xx)\big)^2\ud\xx\int_{-1}^1f_0(1-z^2)^2\ud{z}.
\end{align}
We can obtain a similar bound for the terms including $b_2$.

Finally, the lower bound inequality follows from (\ref{eq:a0}),
(\ref{eq:a1}) and (\ref{eq:a2}).\ef

\section{Small Deborah number limit for the homogeneous system}

This section is devoted to  justifying the small Deborah number limit
for the homogeneous system (\ref{eq:Doi-Onsager-4}). For the
simplicity of notations, throughout this section we denote \beno
&&h_{\nn}=h_{\eta,\nn},\quad \CG_\nn=\CG_{h_\nn},\quad \CA_\nn=\CA_{h_\nn},\\
&&\|f\|_k=\|f\|_{H^k(\BS)},\quad |f|_p=\|f\|_{L^p(\BS)},\quad
\big\langle f,g\big\rangle=\int_{\BS}f(\mm)g(\mm)\ud\mm. \eeno

\subsection{Hilbert expansion}
As in the fluid dynamic limit of the Boltzmann equation \cite{Caf},
we make the Hilbert expansion for $f^\ve(\mm,t)$:
\begin{eqnarray}\label{eq:Hilbert-home}
f^\ve(\mm,t)=\sum_{k=0}^3\ve^kf_k(\mm,
t)+\ve^2f_R^\ve(\mm,t).
\end{eqnarray}
Plugging it into (\ref{eq:Doi-Onsager-4}) and collecting the terms
with the same order with respect to $\ve$, we find that
$f_k(\mm,t)(k=0,1,2,3)$ satisfies \ben
&&\CR f_0+f_0\CR \CU f_0=0,\quad\textrm{ that is}\quad f_0(\mm,t)= h_{\nn(t)}(\mm),\label{eq:hilbert-0}\\
&&\frac{\pa{f_0}}{\pa{t}}=\CG_{\nn(t)}f_1-\CR\cdot(\mm\times\kappa\cdot\mm{f_0}),\label{eq:hilbert-1}\\
&&\frac{\pa{f_1}}{\pa{t}}=\CG_{\nn(t)}f_2+\CR\cdot(f_1\CR\CU{f_1})-\CR\cdot(\mm\times\kappa\cdot\mm{f_1}),
\label{eq:hilbert-2}\\
&&\frac{\pa{f_2}}{\pa{t}}=\CG_{\nn(t)}f_3+\CR\cdot(f_1\CR\CU{f_2})+\CR\cdot(f_2\CR\CU{f_1})
-\CR\cdot(\mm\times\kappa\cdot\mm{f_2}).\label{eq:hilbert-3} \een

The global in time existence of the Hilbert expansion is nontrivial,
since $f_1$ satisfies a nonlinear equation. However, we find that
the part of $f_1$ inside the kernel of $\CG_{\nn(t)}$ satisfies a
linear equation by Lemma \ref{lem:Hibert-2}.

\begin{Proposition}\label{prop:Hilbert-home}
Let $\nn(t)$ be a solution of (\ref{eq:El-home-4}) on $[0,T]$ with
$\lambda$ given by (\ref{eq:lambda}). We can construct smooth
functions $f_k(\mm,t)(k=0, 1,2,3)\in \CP_0$ defined on $[0,T]$ such
that (\ref{eq:hilbert-0})-(\ref{eq:hilbert-3}) hold on $[0,T]$.
\end{Proposition}

We need the following two lemmas in order to prove it.

\begin{Lemma}\cite{Doi, EZ}\label{lem:Hibert-1}
Let $f_0(\mm,t)= h_{\nn(t)}(\mm)$. Then $f_0(\mm,t)$ satisfies
\ben\label{eq:f_0}
\Big\langle\frac{\pa{f_0}}{\pa{t}}+\CR\cdot(\mm\times\kappa\cdot\mm{f_0}),
\psi\Big\rangle=0 \een for any $\psi\in \mathrm{Ker}\CG^*_{\nn(t)}$ if
and only if $\nn(t)$ is a solution of (\ref{eq:El-home-4}) with
$\lambda=\lambda(\alpha)$.
\end{Lemma}

\begin{Lemma}\label{lem:Hibert-2}
For any $\phi, \tilde\phi\in\mathrm{Ker}\CG_\nn$, there holds \beno
\phi=-h_\nn\CU\phi,\quad
\big\langle(\CU\phi)^2,\tilde\phi\big\rangle=0. \eeno
\end{Lemma}

\no{\bf Proof.}\,For $\phi\in\mathrm{Ker}\CG_\nn$, there exists
$\mathbf{\Theta}\in\BR$ such that
$\phi=\mathbf{\Theta}\cdot\CR{h_\nn}$ by Theorem \ref{thm:G-kernel}.
Due to $\ln{h_\nn}=-\CU{h_\nn}$, we see that \beno
\CR{h_\nn}=-h_\nn\CR\CU{h_\nn}=-h_\nn\CU\CR{h_\nn}. \eeno Hence,
$\phi=-h_\nn\CU\phi$. To prove the second equality, we choose the
sphere coordinates such that $\nn=(0,0,1)$ and
$\mm=(\sin\theta\cos\varphi,\sin\theta\sin\varphi,\cos\theta)$. Let
$\tilde\phi=\tilde{\mathbf{\Theta}}\cdot\CR{h_\nn}$ for some
$\tilde{\mathbf{\Theta}}\in \BR$. Then we have \beno
&&\phi=\eta{h}_\nn\cos\theta\sin\theta\big(\Theta_1\sin\varphi-\Theta_2\cos\varphi\big),\\
&&\tilde\phi=\eta{h}_\nn\cos\theta\sin\theta\big(\tilde\Theta_1\sin\varphi-\tilde\Theta_2\cos\varphi\big).
\eeno It is easy to check that \beno
\int_0^{2\pi}(\Theta_1\sin\varphi-\Theta_2\cos\varphi)^2(\tilde\Theta_1\sin\varphi
-\tilde\Theta_2\cos\varphi)\ud\varphi=0, \eeno which implies that
$\int_\BS\frac{\phi^2}{h_\nn^2}\tilde\phi\ud\mm=0$, or equivalently
$\langle(\CU\phi)^2, \tilde\phi\rangle=0$ by
$\phi=-h_\nn\CU\phi$.\ef

\vspace{0.2cm}

\no{\bf Proof of Proposition \ref{prop:Hilbert-home}.}\,Let us first
solve $f_1$ and write $f_1=\phi(t)+\phi^\bot(t)$, where
$\phi\in\mathrm{Ker}\CG_\nn,
\phi^\bot\in(\mathrm{Ker}\CG^*_\nn)^\bot$. Then $\phi^\bot$ will be
determined by (\ref{eq:hilbert-1}), while $\phi$ will be determined
by (\ref{eq:hilbert-2}). However, (\ref{eq:hilbert-1}) has a
solution $\phi^\bot$ if and only if \beno
\Big\langle\frac{\pa{f_0}}{\pa{t}}+\CR\cdot(\mm\times\kappa\cdot\mm{f_0}),
\psi\Big\rangle=0 \eeno for any $\psi\in \mathrm{Ker}\CG^*_\nn$.
From Lemma \ref{lem:Hibert-1}, this is equivalent to require that
$\nn(t)$ is a solution of (\ref{eq:El-home-4}) with
$\lambda=\lambda(\alpha)$. Given $\phi^\bot$, (\ref{eq:hilbert-2})
implies that $\phi$ satisfies \beno
\Big\langle\frac{\pa{\phi}}{\pa{t}},\psi\Big\rangle=
\big\langle\CR\cdot\big((\phi+\phi^\bot)\CR\CU{(\phi+\phi^\bot)}\big)-
\CR\cdot\big(\mm\times(\kappa\cdot\mm){(\phi+\phi^\bot)}\big),\psi\big\rangle-
\Big\langle\frac{\pa{\phi^\bot}}{\pa{t}}, \psi\Big\rangle \eeno for
any $\psi\in \mathrm{Ker}\CG^*_\nn$. Since $\mathrm{Ker}\CG^*_\nn$
is a two dimensional space, we take $\psi_1, \psi_2$ as a base of
$\mathrm{Ker}\CG^*_\nn$ and write $\phi=a_1(t)\psi_1+a_2(t)\psi_2$.
Then we can get an ODE system for $\big(a_1(t), a_2(t)\big)$. For
any $\psi\in\mathrm{Ker}\CG^*_\nn$, we write
$\psi=\CA^{-1}\tilde\phi$ with $\tilde\phi\in\mathrm{Ker}\CG$. Due
to Lemma \ref{lem:Hibert-2}, we find that
\begin{eqnarray*}
&&\big\langle\CR\cdot(\phi\CR\CU\phi), \psi\big\rangle=
\big\langle\CR\cdot(\phi\CR\CU\phi),
\CA_\nn^{-1}\tilde\phi\big\rangle
=-\big\langle\phi\CR\CU\phi, \CR\CA_\nn^{-1}\tilde\phi\big\rangle\\
&&=\big\langle{h_\nn}\CU\phi\CR\CU\phi,
\CR\CA_\nn^{-1}\tilde\phi\big\rangle
=\frac{1}{2}\big\langle\CR\big[(\CU\phi)^2\big],
h_\nn\CR\CA_\nn^{-1}\tilde\phi\big\rangle
=\frac{1}{2}\big\langle(\CU\phi)^2, \tilde\phi\big\rangle=0.
\end{eqnarray*}
This means that $\big(a_1(t), a_2(t)\big)$ satisfies a linear ODE
system, hence is global in time.

Once $f_1$ is determined, we can get $f_2$ and $f_3$ by solving
(\ref{eq:hilbert-2}) and (\ref{eq:hilbert-3}) in a similar way (note
that the equation for $f_2$ is linear).\ef

\subsection{Error estimates}

This subsection is devoted to proving Theorem \ref{thm:Home}. Due to
the weak nonlinearity of the kinetic equation
(\ref{eq:Doi-Onsager-4}), given the initial data $f_0^\ve\in
H^1(\BS)$, it is easy to show by standard energy method that there
exists a unique global solution $f^{\ve}(\mm,t)$ to
(\ref{eq:Doi-Onsager-4}) such that \beno f^\ve\in C([0,+\infty);
H^1(\BS))\cap L^2(0, T;H^2(\BS))\quad\textrm{ for any }T<+\infty.
\eeno Thanks to Proposition \ref{prop:Hilbert-home} and
(\ref{eq:Doi-Onsager-4}), it is easy to find that $f_R^\ve$
satisfies
\begin{align}\label{eq:remainder-home}
\frac{\pa}{\pa{t}}f_R^\ve(\mm,t)=&\f1
\ve\CG_{\nn}f_R^\ve-\CR\cdot\big(\mm\times\kappa\cdot\mm{f_R^\ve}\big)
+\ve\CR\cdot\big(f_R^\ve\CR\CU{f_R^\ve}\big)\nonumber\\
&+\sum_{i=1}^3\ve^{i-1}\CR\cdot\big(f_i\CR\CU{f_R^\ve}+f_R\CR\CU{f_i}\big)+\ve{A},
\end{align}
where $A$ is given by
$$A=\frac{\pa}{\pa{t}}f_3+\sum_{1\le{i,j}\le3,i+j\ge4}\ve^{i+j-4}\CR\big(f_i\CR\CU{f_j}\big).$$
To complete the proof, it suffices to prove that
\ben\label{eq:remiander-est} \|f_R^\ve(t)\|_{-1}\le C
\quad\textrm{ for all }\quad 0\le t\le T. \een For this purpose, we
need the following lemmas.

\begin{Lemma}\label{lem:remaider-1}
\beno
\big[\frac{\pa}{\pa{t}}, \CA_\nn^{-1}\big]g(\mm,t)=
\CA_\nn^{-1}\CR\cdot\Big(\frac{\pa{f_0}}{\pa{t}}\CR(\CA_\nn^{-1}g)\Big).
\eeno
\end{Lemma}

\no{\bf Proof.}\, Direct calculations give that
\begin{align*}
\frac{\pa}{\pa{t}}(\CA_{\nn}g)=&-\frac{\pa}{\pa{t}}\CR\cdot(h_\nn\CR{g})=
-\CR\cdot\Big(\frac{\partial{h_\nn}}{\partial{t}}\CR{g}\Big)-\CR\cdot\Big(h_\nn\CR{\frac{\pa{}g}{\pa{t}}}\Big)\\
=&-\CR\cdot\Big(\frac{\partial{h_\nn}}{\partial{t}}\CR{g}\Big)+\CA_\nn{\frac{\pa{g}}{\pa{t}}}.
\end{align*}
Replacing $g$ by $\CA_\nn^{-1}g$, the lemma follows.\ef

\begin{Lemma}\label{lem:remaider-2}
For any vector field $V\in C^1(\BS)$, there holds \beno
\big\langle\CR\cdot(V{f}),\CA_\nn^{-1}f\big\rangle\le
C\big(|V|_{\infty}+|\CR{V}|_{\infty}\big)\big\langle{f},
\CA_\nn^{-1}f\big\rangle. \eeno
\end{Lemma}

\no{\bf Proof.}\,Let $V=(V_1,~V_2,~V_3)^T$,
$\CR{f}=(R_1f,~R_2f,~R_3f)^T$, and $g=\CA_\nn^{-1}f$. Recalling $\CA_\nn f=-\CR\cdot(h_\nn{\CR
f})$, we get
\begin{align*}
\big\langle\CR\cdot(V{f}), \CA_\nn^{-1}f\big\rangle=&-\big\langle{R_i}[V_iR_j(h_\nn{R_jg})], g\big\rangle\\
=&-\big\langle{h_\nn}{R_jg},(R_jV_i)R_ig+V_iR_jR_ig\big\rangle\\
=&-\big\langle{h_\nn}{R_jg}, (R_jV_i)R_ig-V_i\epsilon^{kji}R_kg\big\rangle+\big\langle{h_\nn}{R_jg},V_iR_iR_jg\big\rangle\\
=&-\big\langle{h_\nn}{R_jg},
(R_jV_i)R_ig-V_i\epsilon^{kji}R_kg\big\rangle
-\frac{1}{2}\big\langle{R_i(h_\nn}V_i),(R_jg)^2\big\rangle,
\end{align*}
which implies the lemma by using  the fact that \beno
|\CR{g}|_2\le{C}\|f\|_{-1}\le{C}\big\langle{f},
\CA_\nn^{-1}f\big\rangle. \eeno The proof is finished. \ef

Now we are in position to prove (\ref{eq:remiander-est}). With the
help of Lemma \ref{lem:remaider-1}, we make the energy estimate for
(\ref{eq:remainder-home}) to obtain \beno
&&\frac{1}{2}\frac{d}{d{t}}\big\langle{f_R^\ve},
\CA_\nn^{-1}f_R^\ve\big\rangle\\
&&=\big\langle\frac{\pa}{\pa{t}}f_R^\ve,\CA_\nn^{-1}f_R^\ve\big\rangle+\f12\Big\langle
\CA_\nn^{-1}\CR\cdot\Big(\frac{\pa{f_0}}{\pa{t}}\CR(\CA_\nn^{-1}f_R^\ve)\Big), f_R^\ve\Big\rangle\\
&&=\f 1 \ve\big\langle\CG_{\nn}f_R^\ve,
\CA_\nn^{-1}f_R^\ve\big\rangle
-\big\langle\CR\cdot(\mm\times\kappa\cdot\mm{f_R^\ve}),\CA_\nn^{-1}f_R^\ve\big\rangle\\
&&\quad+\big\langle\ve\CR\cdot(f_R^\ve\CR\CU{f_R^\ve}),
\CA_\nn^{-1}f_R^\ve\big\rangle+\sum_{i=1}^3
\ve^{i-1}\big\langle\CR\cdot(f_i\CR\CU{f_R^\ve}+f_R^\ve\CR\CU{f_i}),
\CA_\nn^{-1}f_R^\ve\big\rangle
\nonumber\\&&\quad+\ve\big\langle{A},\CA_\nn^{-1}f_R^\ve\big\rangle
-\frac{1}{2}\Big\langle\frac{\pa{f_0}}{\pa{t}}\CR(\CA_\nn^{-1}f_R^\ve),
\CR\CA_\nn^{-1}f_R^\ve\Big\rangle. \eeno Since $h_\nn$ is a stable
critical point, we know from Proposition \ref{prop:H-positive} that
\beno \big\langle\CG_{\nn}f_R^\ve, \CA_\nn^{-1}f_R^\ve\big\rangle
=-\big\langle \CH_{\nn}f_R^\ve, f_R^\ve\big\rangle\le 0. \eeno We
infer from Lemma \ref{lem:remaider-2} that \beno
&&-\big\langle\CR\cdot(\mm\times\kappa\cdot\mm{f_R^\ve}),\CA_\nn^{-1}f_R^\ve\big\rangle
\le C\big\langle{f_R^\ve}, \CA_\nn^{-1}f_R^\ve\big\rangle,\\
&&\big\langle\ve\CR\cdot(f_R^\ve\CR\CU{f_R^\ve}),
\CA_\nn^{-1}f_R^\ve\big\rangle
\le C\ve\big\langle{f_R^\ve}, \CA_\nn^{-1}f_R^\ve\big\rangle^\f32,\\
&&\big\langle\CR\cdot(f_R^\ve\CR\CU{f_i}),
\CA_\nn^{-1}f_R^\ve\big\rangle \le
C\|f_i\|_{-1}\big\langle{f_R^\ve}, \CA_\nn^{-1}f_R^\ve\big\rangle,
\eeno while the other terms on the right hand side are bounded by
\beno \Big(\sum_{i=1}^3\ve^{i-1}|f_i|_2+|\partial_t
f_0|_2\Big)\big\langle{f_R^\ve}, \CA_\nn^{-1}f_R^\ve\big\rangle
+\ve\|A\|_{-1}\big\langle{f_R^\ve},
\CA_\nn^{-1}f_R^\ve\big\rangle^\f12. \eeno Here we use the fact that
$|\CR^k\CU{f_R^\ve}|_\infty$ for $k\in \NN$ can be bounded by
$\big\langle{f_R^\ve}, \CA_\nn^{-1}f_R^\ve\big\rangle^{\f12}$ and
$\big\langle{f_R^\ve}, \CA_\nn^{-1}f_R^\ve\big\rangle^{\f12}\sim
\|f_R^\ve\|_{-1}$. In conclusion, we obtain \beno
\frac{\ud}{\ud{t}}\|f_R^\ve\|_{-1}^2\le{C}\big(\|f_R^\ve\|_{-1}^2+\ve\|f_R^\ve\|_{-1}^3+\ve\|f_R^\ve\|_{-1}\big).
\eeno This implies (\ref{eq:remiander-est}).\ef

\subsection{The Lesile stress and coefficients}

This subsection is devoted to proving Theorem \ref{thm:stress-home}.
We introduce the 2-order tensor $\QQ_2[f]$ and 4-order tensor
$\QQ_4[f]$ as follows \beno
&&\QQ_2[f]=\big\langle\mm\mm-\frac{1}{3}\II\big\rangle_f,\\
&&\QQ_{4\alpha\beta\gamma\mu}[f]=\Big\langle{m}_{\alpha}{m}_{\beta}{m}_{\gamma}{m}_{\mu}
-\frac{1}{7}(m_{\alpha}m_{\beta}\delta_{\gamma\mu}+
m_{\gamma}m_{\mu}\delta_{\alpha\beta}+m_{\alpha}m_{\gamma}\delta_{\beta\mu}
+m_{\beta}m_{\mu}\delta_{\alpha\gamma}\nonumber\\
&&\qquad+m_{\alpha}m_{\mu}\delta_{\beta\gamma}+
m_{\beta}m_{\gamma}\delta_{\alpha\mu})+\frac{1}{35}(\delta_{\alpha\beta}\delta_{\gamma\mu}
+\delta_{\alpha\gamma}\delta_{\beta\mu}+\delta_{\alpha\mu}\delta_{\beta\gamma})\Big\rangle_f.
\eeno

\begin{Lemma}\label{lem:Q-tensor}
It holds that \beno
&&\QQ_2[h_\nn]=\langle{P}_2(\mm\cdot\nn)\rangle_{h_\nn}(\nn\nn-\frac{1}{3}),\\
&&\QQ_{4\alpha\beta\gamma\mu}[h_\nn]=\langle{P}_4(\mm\cdot\nn)\rangle_{h_\nn}\Big(n_{\alpha}n_{\beta}n_{\gamma}n_{\mu}
-\frac{1}{7}(n_{\alpha}n_{\beta}\delta_{\gamma\mu}+
n_{\gamma}n_{\mu}\delta_{\alpha\beta}+n_{\alpha}n_{\gamma}\delta_{\beta\mu}\nonumber\\
&&\qquad+n_{\beta}n_{\mu}\delta_{\alpha\gamma}+n_{\alpha}n_{\mu}\delta_{\beta\gamma}+
n_{\beta}n_{\gamma}\delta_{\alpha\mu})+\frac{1}{35}(\delta_{\alpha\beta}
\delta_{\gamma\mu}+\delta_{\alpha\gamma}\delta_{\beta\mu}
+\delta_{\alpha\mu}\delta_{\beta\gamma})\Big).~~~~~~~~\label{Q4}
\eeno
\end{Lemma}

\no{\bf Proof.}\,We only prove the first equality, the proof of the
second equality is similar but more complicated. Since both sides of
the first equality are tensors, they are coordinate-independent. So,
we may choose the sphere coordinates such that $\nn=(0,0,1)$ and
$\mm=(\sin\theta\cos\varphi,\sin\theta\sin\varphi,\cos\theta)$.
Since $h_\nn(\mm)$ depends only on $\cos\theta$, it is easy to check
that \beno
&&\int_\BS{m}_i{m}_jh_\nn(\mm)\ud\mm=0\quad \textrm{for}\quad i\neq{j},\\
&&\int_\BS({m}_3^2-\frac{1}{3})h_\nn(\mm)\ud\mm=\frac{2}{3}\int_\BS\frac{1}{2}\big(3\cos^2\theta-1\big)h_\nn(\mm)\ud\mm=
\frac{2}{3}\big\langle{P_2(\cos\theta)}\big\rangle_{h_\nn},\\
&&\int_\BS({m}_1^2-
\frac{1}{3})h_\nn(\mm)\ud\mm=\int_\BS({m}_2^2-\frac{1}{3})h_\nn(\mm)\ud\mm=
-\frac{1}{3}\big\langle{P_2(\cos\theta)}\big\rangle_{h_\nn}, \eeno
which give the first equality. \ef

\begin{Lemma}\label{lem:Diff-mean}
Let $f^{\ve}(t)$ be given by Theorem \ref{thm:Home} and $P(x)$ be a
smooth function on $\RR$. Then we have \beno \big|\big\langle
P(\mm\cdot\nn(t))
\big\rangle_{f^{\ve}(t)}-\big\langle{P}(\mm\cdot\nn(t))
\big\rangle_{h_{\nn(t)}}\big|\le{C\ve}. \eeno
\end{Lemma}

\no{\bf Proof.}\,By the definition, we get \beno \big\langle
P(\mm\cdot\nn)\big\rangle_{f^{\ve}}-\big\langle{P}(\mm\cdot\nn)\big\rangle_{h_\nn}
=\int_\BS{P(\mm\cdot\nn)}\Big(\sum_{k=1}^3\epsilon^kf_k(\mm,
t)+\epsilon^2f_R^\ve(\mm,t)\Big)\ud\mm. \eeno Using the facts that
$f_1, f_2, f_3$ are bounded and \beno
\int_\BS{P(\mm\cdot\nn)}f_R^\ve(\mm,t)\ud\mm&=&\int_\BS\CA_\nn{P(\mm\cdot\nn)}\CA_\nn^{-1}f_R^\ve(\mm,t)\ud\mm\\
&\le&\|\CA_\nn{P(\mm\cdot\nn)}\|_0\|f_R^\ve\|_{-1},\quad~ \eeno the
lemma follows.\ef \vspace{0.2cm} Now we are in position to prove
Theorem \ref{thm:stress-home}. Direct computation shows that
\begin{align*}
\frac{d}{d{t}}\QQ_2[f^\ve]=&\frac{1}{\ve}\int_{\BS}(\mm\mm-\frac{1}{3}\II)\big(\CR\cdot(f^\ve\CR\mu^\ve)
-\ve\CR(\mm\times\kappa\cdot\mm{f^{\ve}})\big)\ud\mm\\
=&\frac{1}{\ve}\big\langle\mm\times\CR\mu^{\ve}\mm+\mm\mm\times\CR\mu^\ve\big\rangle_{f^\ve}
-\big(2\DD:\langle\mm\mm\mm\mm\rangle_{f^\ve}\nonumber\\
&-\DD\cdot\langle\mm\mm\rangle_{f^\ve}+\BOm\cdot\langle\mm\mm\rangle_{f^\ve}
-\langle\mm\mm\rangle_{f^\ve}\cdot(\DD+\BOm)\big).
\end{align*}
So by Lemma \ref{lem:Diff-mean}, the stress $\sigma^\ve$ can be
written as
\begin{align*}
\sigma^\ve=&\frac{1}{2}\DD:\langle\mm\mm\mm\mm\rangle_{f^\ve}-
\frac{1}{\ve}\langle\mm\mm\times\CR\mu[f^\ve]\rangle_{f^\ve}\\
=&\frac{1}{2}\DD:\langle\mm\mm\mm\mm\rangle_{f^\ve}-\frac{1}{2}
\Big(2\DD:\langle\mm\mm\mm\mm\rangle_{f^\ve}-\DD\cdot\langle\mm\mm\rangle_{f^\ve}\\
&+\BOm\cdot\langle\mm\mm\rangle_{f^\ve}
-\langle\mm\mm\rangle_{f^\ve}\cdot(\DD+\BOm)+\frac{\ud}{\ud{t}}\QQ_2[f^\ve]\Big)\\
=&\frac{1}{2}\DD:\langle\mm\mm\mm\mm\rangle_{h_\nn}-\frac{1}{2}
\Big(2\DD:\langle\mm\mm\mm\mm\rangle_{h_\nn}-\DD\cdot\langle\mm\mm\rangle_{h_\nn}\\
&+\BOm\cdot\langle\mm\mm\rangle_{h_\nn}
-\langle\mm\mm\rangle_{h_\nn}\cdot(\DD+\BOm)+\frac{\ud}{\ud{t}}\QQ_2[h_\nn]\Big)\\
&-\frac{1}{2}\QQ_2\Big[\sum_{k=1}^3\epsilon^k\frac{\partial}{\partial{t}}f_k(\mm,
t)+\epsilon^2\frac{\partial}{\partial{t}}f_R^\ve(\mm,t)\Big]+C\ve.
\end{align*}
For any constant vector $U,~V$, we have
\begin{eqnarray*}
&&\ve^2U^T\cdot\QQ_2\Big[\frac{\partial}{\partial{t}}f_R^\ve(\mm,t)\Big]\cdot{V}=
\ve^2\int_\BS(\mm\cdot{U})
(\mm\cdot{V})\Big[\ve^{-1}\CG_{\nn}f^\ve_R-\CR(\mm\times\kappa\cdot\mm{f^\ve_R})\nonumber\\
&&\qquad+\ve\CR\cdot(f^\ve_R\CR\CU{f^\ve_R})+\sum_{i=1}^3\ve^{i-1}\CR\cdot(f^\ve_i\CR\CU{f^\ve_R}
+f^\ve_R\CR\CU{f^\ve_i})+\ve{A}\Big]\ud\mm,
\end{eqnarray*}
which is bounded by $C\ve\|f_R\|_{-1}$ by using the argument of
integration by parts. Then we infer that
\begin{align*}
\sigma^\ve
=&\frac{1}{2}\DD:\langle\mm\mm\mm\mm\rangle_{h_\nn}-\frac{1}{2}
\Big(2\DD:\langle\mm\mm\mm\mm\rangle_{h_\nn}-\DD\cdot\langle\mm\mm\rangle_{h_\nn}\\
&+\BOm\cdot\langle\mm\mm\rangle_{h_\nn}
-\langle\mm\mm\rangle_{h_\nn}\cdot(\DD+\BOm)+\frac{\ud}{\ud{t}}\QQ_2[h_\nn]\Big)+C\ve.
\end{align*}
We see from the definition of $\QQ_4$ that \beno
\DD:\langle\mm\mm\mm\mm\rangle_{h_\nn}=\DD:\QQ_4+\frac{1}{7}\DD:\langle\mm\mm\rangle_{h_\nn}\II+\frac{2}{7}(\DD:\QQ_2+\QQ_2:\DD)+\frac{2}{15}\DD.
\eeno Using Lemma \ref{lem:Q-tensor}, we can calculate
\begin{align*}
\sigma^\ve=&\frac{1}{2}\Big(-S_4(\DD:\nn\nn)\nn\nn-\frac{S_2}{7}(\DD:\nn\nn)\II
-S_2(\nn\NN+\NN\nn)\nonumber\\
&+\big(\frac{8}{15}-\frac{10}{21}S_2-\frac{2}{35}S_4\big)\DD+\big(\frac{5}{7}S_2+\frac{2}{7}S_4\big)(\nn\nn\cdot\DD+\DD\cdot\nn\nn)\Big)+C\ve.
\end{align*}
So, we conclude that
$$\sigma^\ve-\sigma^L-p\II=-\frac{S_2}{2}\Big(\nn(\frac{1}{\lambda}\NN-\DD\cdot\nn)-(\frac{1}{\lambda}\NN-\DD\cdot\nn)\nn\Big)+C\ve=C\ve,
$$
here we used the equation (\ref{eq:El-home-4}). This completes the
proof of Theorem \ref{thm:stress-home}.\ef

\section{Small Deborah number limit for the inhomogeneous system}

This section is devoted to  justifying the small Deborah number limit
for the inhomogeneous system (\ref{eq:LCP-nonL-f})-(\ref{eq:LCP-nonL-v}).
Throughout this section, we will use the following notations.
$\langle,\rangle$ denotes the inner product on $L^2(\Om\times\BS)$ or $L^2(\Om)$.
We also denote $\|f\|_{L^p}=\|f\|_{L^p(\Om\times\BS)}(\|f\|_{L^p(\Om)})$ for $f$ defined on $\Om\times \BS$(for
$f$ defined on $\Om$), and $\|\cdot\|_{H^{0,k}}=\big\|\|\cdot\|_{H^k(\Om)}\big\|_{L^2(\BS)}$.

\subsection{Hilbert expansion}
Due to the  choice of $g$, we have the formal expansion to the operator
$\CU_\ve$:
\begin{align*}
\CU_\ve f-\CU
f=&\int_\BR\int_\BS\alpha|\mm\times\mm'|^2g_\ve(\xx-\xx')
(f(\xx',\mm',t)-f(\xx,\mm',t))\ud\mm'\ud\xx'\\
=&\int_\BR\int_\BS\alpha|\mm\times\mm'|^2g(\yy)
(f(\xx+\sqe\yy,\mm',t)-f(\xx,\mm',t))\ud\mm'\ud\yy\\
=&\int_\BR\int_\BS\alpha|\mm\times\mm'|^2g(\yy)
\Big(\sum_{k\ge1}\frac{\ve^{\frac{k}{2}}}{k!}(\yy\cdot\nabla)^kf(\xx,\mm',t)\Big)\ud\mm'\ud\yy\\
=&\int_\BR\int_\BS\alpha|\mm\times\mm'|^2g(\yy)
\Big(\sum_{k\ge1}\frac{\ve^{{k}}}{(2k)!}(\yy\cdot\nabla)^{2k}f(\xx,\mm',t)\Big)\ud\mm'\ud\yy.
\end{align*}
We denote
\begin{eqnarray*}
U_k[f](\xx,\mm,t)=\frac{1}{(2k)!}\int_\BR\int_\BS\alpha|\mm\times\mm'|^2g(\yy)
(\yy\cdot\nabla)^{2k}f(\xx,\mm',t)\ud\mm'\ud\yy\quad\textrm{
for}\quad k\ge 1.
\end{eqnarray*}
Formally, we have \ben \CU_\ve
f=U_0[f]+\ve{U_1[f]}+\ve^2U_2[f]+\cdots,\quad U_0[f]=\CU f. \een
Then we make a formal expansion for the solution of
(\ref{eq:LCP-nonL-f})-(\ref{eq:LCP-nonL-v}):
\beno
&&f^\ve(\xx,\mm,t)=\sum_{k=0}^3\ve^kf_k(\xx,\mm,t)+\ve^3f_R^\ve(\xx,\mm,t),\\
&&v^\ve(\xx,t)=\sum_{k=0}^2\ve^kv_k(\xx,t)+\ve^3v_R^\ve(\xx,t).
\eeno Plugging them into (\ref{eq:LCP-nonL-f})-(\ref{eq:LCP-nonL-v})
and collecting the terms with the same order with respect to $\ve$,
we find that \ben\label{eq:Hilbert-L-00} \CR{f_0}+f_0\CR{\CU
f_0}=0,\quad \textrm{that is},\quad f_0=h_{\eta, \nn(t,\xx)}(\mm);
\een
and for the terms of order $O(1)$, there hold
\begin{align}
\frac{\pa{f_0}}{\pa{t}}+\vv_0\cdot\nabla{f_0}=&\CG_{f_0}f_1
+\CR\cdot(f_0\CR{U_1}f_0)-\CR\cdot(\mm\times(\nabla\vv_0)^T\cdot\mm{f_0}),\label{eq:Hilbert-L-10}\\
\frac{\pa{\vv_0}}{\pa{t}}+\vv_0\cdot\nabla\vv_0=&-\nabla{p_0}+\frac{\gamma}{Re}\Delta\vv_0
+\frac{1-\gamma}{2Re}\nabla\cdot(\DD_0:\langle\mm\mm\mm\mm\rangle_{f_0})\nonumber\\
&-\frac{1-\gamma}{Re}\Big\{\nabla\cdot\big\langle\mm\mm\times(\CR{f_1}
+\sum_{\scriptscriptstyle{i+j+k=1}}f_i{U}_jf_k)\big\rangle_1+\langle\nabla{U}_1f_0\rangle_{f_0}\Big\};\label{eq:Hilbert-L-11}
\end{align}
and for the terms of order $O(\ve)$, there hold
\begin{align}
&\frac{\pa{f_1}}{\pa{t}}+\vv_0\cdot\nabla{f_1}=\CG_{f_0}f_2+\CR\cdot\Big(\sum_{\scriptscriptstyle{i+j+k=2,j\ge1}}
f_i\CR{U}_jf_k\Big)-\vv_1\cdot\nabla{f_0}\nonumber\\
&\qquad\qquad-\CR\cdot\big(\mm\times(\nabla\vv_0)^T\cdot\mm{f_1}+\mm\times(\nabla\vv_1)^T\cdot\mm{f_0}\big),\label{eq:Hilbert-L-20}\\
&\frac{\pa{\vv_1}}{\pa{t}}+\vv_0\cdot\nabla\vv_1=-\nabla{p_1}+\frac{\gamma}{Re}\Delta\vv_1
+\frac{1-\gamma}{2Re}\nabla\cdot\Big(\sum_{i+j=1}\DD_i:\langle\mm\mm\mm\mm\rangle_{f_j}\Big)-\vv_1\cdot\nabla\vv_0\nonumber\\
&\qquad-\frac{1-\gamma}{Re}\Big\{\nabla\cdot\big\langle\mm\mm\times(\CR{f_2}
+\sum_{\scriptscriptstyle{i+j+k=2}}f_i\CR{U}_jf_k)\big\rangle_1
+\big\langle\sum_{\scriptscriptstyle{i+j+k=2}}f_i\nabla{U}_jf_k\big\rangle_1\Big\};\label{eq:Hilbert-L-21}
\end{align}
and for the terms of $O(\ve^2)$, there hold
\begin{align}
&\frac{\pa{f_2}}{\pa{t}}+\vv_0\cdot\nabla{f_2}
=\CG_{f_0}f_3+\CR\cdot\Big(\sum_{\scriptscriptstyle{i+j+k=3,j\ge1}}f_i\CR{U_j{f_k}}\Big)\non\\
&\qquad\quad-\CR\cdot\Big(\sum_{\scriptscriptstyle{i+j=2}}\mm\times\big((\nabla\vv_i)^T\cdot\mm\big){f_j}\Big)
-\vv_1\cdot\nabla{f_1}-\vv_2\cdot\nabla{f_0},\label{eq:Hilbert-L-30}\\
&\frac{\pa{\vv_2}}{\pa{t}}+\vv_0\cdot\nabla\vv_2
=-\nabla{p_2}+\frac{\gamma}{Re}\Delta\vv_2
+\frac{1-\gamma}{2Re}\nabla\cdot\Big(\sum_{\scriptscriptstyle{i+j=2}}\DD_i:\langle\mm\mm\mm\mm\rangle_{f_j}\Big)\nonumber\\
&\quad-\frac{1-\gamma}{Re}\Big\{\nabla\cdot\big\langle\mm\mm\times(\CR{f_3}+\sum_{\scriptscriptstyle{i+j+k=3}}f_i\CR{U}_j{f}_k)\big\rangle_1
+\langle\sum_{\scriptscriptstyle{i+j+k=3}}f_i\nabla{U}_jf_k\rangle_1\Big\}\nonumber\\
&\qquad\qquad-\vv_1\cdot\nabla\vv_1-\vv_2\cdot\nabla\vv_0.\label{eq:Hilbert-L-31}
\end{align}

\begin{Proposition}\label{prop:Hilbert-inhome}
Let $(\vv_0, \nn)\in C\big([0,T];H^{20}(\Om)\big)$ be a solution of
(\ref{eq:EL-nh})-(\ref{eq:EL-vh}) with $\lambda$ given by
(\ref{eq:lambda}) and Leslie coefficients defined by (\ref{Leslie
cofficients1})-(\ref{Leslie cofficients2}) on $[0,T]$. Then there
exist the functions $f_i\in
C\big([0,T];H^{20-4i}(\Om\times\BS)\big)(i=0,1,2,3)$ and $\vv_i\in
C\big([0,T];H^{20-4i}(\Om)\big)(i=0,1,2)$ such that
(\ref{eq:Hilbert-L-00})-(\ref{eq:Hilbert-L-31}) holds on $[0,T]$.
\end{Proposition}

We need the following lemma.

\begin{Lemma}\cite{EZ}\label{lem:Hibert-11}
Let $f_0= h_{\eta, \nn(\xx, t)}(\mm)$. Then $f_0(\xx,\mm,t)$ satisfies
\beno
\Big\langle\frac{\pa{f_0}}{\pa{t}}+\vv_0\cdot\nabla{f_0}+\CR\cdot(\mm\times(\nabla\vv_0)^T\cdot\mm{f_0})
+\CR\cdot(f_0\CR{U}_1f_0),
\psi\Big\rangle_{L^2(\BS)}=0
\eeno
for any $\psi\in \mathrm{Ker}\CG^*_{f_0}$
if and only if $\nn(\xx, t)$ is a solution of (\ref{eq:EL-nh}) with $\vv=\vv_0$.
\end{Lemma}

\no{\bf Proof of Proposition \ref{prop:Hilbert-inhome}}.
We denote by $\Pin$ the projection operator from $\CP_0(\BS)$ to $\mathrm{Ker}\CG_{f_0}$, and  denote by $\Pout$
the projection operator from $\CP_0(\BS)$ to $(\mathrm{Ker}\CG^*_{f_0})^{\bot}$. Let $\Pin{f}_1=\phi_1,~\Pout{f}_1=\psi_1$.
Now $\psi_1$ will be determined by (\ref{eq:Hilbert-L-10}), whose existence is ensured by Lemma \ref{lem:Hibert-11}.
Once $\psi_1$ is determined, it can be proved that the equation (\ref{eq:Hilbert-L-11}) is equivalent to (\ref{eq:EL-vh}), see \cite{EZ}.
Now we solve $(\phi_1, \vv_1)$. In what follows, we denote by  $L(\phi,\vv)$ the terms which only depend on $(\phi,\vv)$(not their
derivatives) linearly. Let $\phi_1=\nn^\bot\cdot\CR{f}_0$. We have
\begin{align*}
\big(\frac{\pa}{\pa{t}}+\vv_0\cdot\nabla\big)\phi_1
=\big(\frac{\pa}{\pa{t}}+\vv_0\cdot\nabla\big)\nn^\bot\cdot\CR{f}_0
+\nn^\bot\cdot\CR\big((\frac{\pa}{\pa{t}}+\vv_0\cdot\nabla){f}_0\big).
\end{align*}
This means that
\begin{align}
&\Pout\Big((\frac{\pa}{\pa{t}}+\vv_0\cdot\nabla)\phi_1\Big)
=\Pout\Big(\nn^\bot\cdot\CR\big((\frac{\pa}{\pa{t}}+\vv_0\cdot\nabla){f}_0\big)\Big)\triangleq L(\phi_1),\label{eq:out-1}\\
&\Pin\Big(\big(\frac{\pa}{\pa{t}}+\vv_0\cdot\nabla\big)\phi_1\Big)
=\big(\frac{\pa}{\pa{t}}+\vv_0\cdot\nabla\big)\phi_1-L(\phi_1).\label{eq:in-1}
\end{align}
We also have
\begin{align}
&\Pin\big(\CG_{f_0}f_2+\CR\cdot(\phi_1\CR \CU\phi_1)\big)=0.\label{eq:in-2}
\end{align}
For a matrix $\kappa$, we denote
\beno
&&\mathcal{K}(\kappa)=\Pin\big(\CR\cdot\big(\mm\times\kappa\cdot\mm{f_0}\big)\big),\quad
\mathcal{L}(\kappa)=\Pout\big(\CR\cdot\big(\mm\times\kappa\cdot\mm{f_0}\big)\big),\\
&&\mathcal{B}_{in}(\phi_1)=\Pin\big(\CR\cdot(f_0\CR U_1\phi_1)\big),
\quad \mathcal{B}_{out}(\phi_1)=\Pout\big(\CR\cdot(f_0\CR U_1\phi_1)\big).
\eeno
Taking $\Pin$ for the equation of $f_1$, we get by (\ref{eq:in-1}) and (\ref{eq:in-2}) that
\begin{align}
&\big(\frac{\pa}{\pa{t}}+\vv_0\cdot\nabla\big)\phi_1=L(\phi_1)
-\mathcal{K}((\nabla\vv_1)^T)+\mathcal{B}_{in}(\phi_1)\non\\
&\quad+\Pin\big(\CR\cdot(\phi_1\CR\CU\psi_1+\psi_1\CR\CU\phi_1+\psi_1\CR\CU\psi_1)\big)-\vv_1\cdot\nabla{f}_0\non\\
&\quad-\Pin\big(\CR\cdot\big(\mm\times(\nabla\vv_0)^T\cdot\mm{f_1}\big)\big)+
\Pin\big(\CR\cdot(f_0\CR U_1\psi_1)\big)+\Pin\big(\CR\cdot(f_1\CR U_1{f}_0)\big)\non\\
&=-\mathcal{K}((\nabla\vv_1)^T)+\mathcal{B}_{in}(\phi_1)+L(\phi_1,\vv_1).\label{eq:f-in}
\end{align}
Taking $\Pout$ for the  equation of $f_1$, we get by (\ref{eq:out-1}) that
\begin{align*}
&-\big(\CG_{f_0}f_2+\CR\cdot(\phi_1\CR\CU\phi_1)\big)=-L(\phi_1)+\mathcal{L}((\nabla\vv_1)^T)-\mathcal{B}_{out}(\phi_1)\\
&\quad+\Pout\big(\CR\cdot(\phi_1\CR \CU\psi_1+\psi_1\CR\CU\phi_1+\psi_1\CR\CU\psi_1)\big)\\
&\quad-\Pout\big(\CR\cdot\big(\mm\times(\nabla\vv_0)^T\cdot\mm{f_1}\big)\big)+
\Pout\big(\CR\cdot(f_0\CR U_1\psi_1))+\Pout(\CR\cdot(f_1\CR U_1{f}_0\big))\\
&=\mathcal{L}((\nabla\vv_1)^T)-\mathcal{B}_{out}(\phi_1)-L(\phi_1,\vv_1).
\end{align*}
So, the stress in the equation of $\vv_1$ can be rewritten as
\begin{align}
&-\big\langle\mm\mm\times\big(\CR{f_2}+\sum_{\scriptscriptstyle{i+j+k=2}}f_i\CR{U}_jf_k\big)\big\rangle_1\nonumber\\
&=-\frac{1}{2}\big\langle(\mm\mm-\frac{1}{3}\II)(\CG_{f_0}f_2+\CR\cdot(\phi_1\CR\CU\phi_1))\big\rangle_1
-\big\langle\mm\mm\times(f_0\CR U_1\phi_1)\big\rangle_1+L(\phi_1)\non\\\nonumber
&=\frac{1}{2}\big\langle(\mm\mm-\frac{1}{3}\II)(\mathcal{L}((\nabla\vv_1)^T)-\mathcal{B}_{out}(\phi_1))\big\rangle_1
-\big\langle\mm\mm\times(f_0\CR U_1\phi_1)\big\rangle_1+L(\phi_1,\vv_1)\\\nonumber
&\triangleq \sigma_1+\sigma_2+L(\phi_1,\vv_1).
\end{align}
Set $\sigma_3=\DD_1:\langle\mm\mm\mm\mm\rangle_{f_0}$. Then the equation for $\vv_1$ can be rewritten as
\begin{align}
&\frac{\pa{\vv_1}}{\pa{t}}-\frac{\gamma}{Re}\Delta\vv_1+\vv_0\cdot\na \vv_1+\nabla{p_1}=
\frac{1-\gamma}{2Re}\nabla\cdot\big(\sigma_1+\sigma_2+\sigma_3\big)\nonumber\\
&\qquad-\frac{1-\gamma}{Re}\big\langle f_0\nabla{U}_1\phi_1\big\rangle_1+\na L(\phi_1,\vv_1)+L(\phi_1,\vv_1).\label{eq:v-in}
\end{align}
In order to solve (\ref{eq:f-in})-(\ref{eq:v-in}), we introduce the energy functional
\begin{align*}
E(t)=\langle\phi_1, \phi_1\rangle+\langle\phi_1, U_1\phi_1\rangle+\frac{Re}{1-\gamma}\langle\vv_1,\vv_1\rangle.
\end{align*}
Due to the choice of $g$, it is easy to see that
\begin{align*}
\langle\phi_1, \phi_1\rangle+\langle\phi_1, U_1\phi_1\rangle\ge{c}\big(\langle\phi_1, \phi_1\rangle+\langle\nabla\phi_1, \nabla\phi_1\rangle\big).
\end{align*}
Notice that
$
\langle\nabla\cdot\sigma_3,\vv_1\rangle=-\langle\DD_1:\langle\mm\mm\mm\mm\rangle_{f_0}, \DD_1\rangle\le0
$
and by Lemma \ref{Lem:stress-dissipation} and Lemma \ref{Lem:ker-out-proj},
\begin{align*}
&\big\langle-\mathcal{K}((\nabla\vv_1)^T)+\mathcal{B}_{in}(\phi_1),U_1\phi_1\big\rangle-\big\langle\nabla\cdot(\sigma_1+\sigma_2), \vv_1\big\rangle\\
&\le\big\langle-\mathcal{K}((\nabla\vv_1)^T)+\mathcal{B}_{in}(\phi_1),-\CA^{-1}_{f_0}(\mathcal{B}_{in}(\phi_1)+\mathcal{B}_{out}(\phi_1))\big\rangle\\
&\qquad+\frac{1}{2}\big\langle\mathcal{B}_{out}(\phi_1),\mm\cdot\nabla\vv_1\cdot\mm\big\rangle
+\big\langle{U}_1\phi_1,\CR\big(\mm\times(\nabla\vv_1)^T\cdot\mm{f}_0\big)\big\rangle\\
&=\big\langle-\mathcal{K}((\nabla\vv_1)^T)+\mathcal{B}_{in}(\phi_1),-\CA^{-1}_{f_0}(\mathcal{B}_{in}(\phi_1)+\mathcal{B}_{out}(\phi_1))\big\rangle
+\frac{1}{2}\big\langle\mathcal{B}_{out}(\phi_1),-2\CA^{-1}_{f_0}(\mathcal{K}(\DD_1)+\mathcal{L}(\DD_1))\big\rangle\\
&\qquad+\big\langle-\CA^{-1}_{f_0}(\mathcal{B}_{in}(\phi_1)+\mathcal{B}_{out}(\phi_1)),\mathcal{K}((\nabla\vv_1)^T)+\mathcal{L}((\nabla\vv_1)^T)\big\rangle\\
&=-\big\langle\mathcal{B}_{in}(\phi_1), \CA^{-1}_{f_0}\mathcal{B}_{in}(\phi_1)\big\rangle \le 0.
\end{align*}
Then by a simple energy estimate, we can deduce that
\beno
\f d {dt}E(t)\le C\big(1+E(t)\big).
\eeno
The estimate of the higher order derivative for $(\phi_1,\vv_1)$ can be obtained by introducing a similar energy functional.
Once $(f_1,\vv_1)$ is determined, we can get $(f_2,\vv_2)$ and $f_3$ by solving
(\ref{eq:Hilbert-L-30})-(\ref{eq:Hilbert-L-31}) in a similar way.
Here we omit the details.
\ef

\subsection{The remainder equations}

We denote
\begin{align*}
&\tv=\vv_0+\ve\vv_1+\ve^2\vv_2,\quad  \tf=f_1+\ve{f_2}+\ve^2{f_3},
\quad\tD=\DD_0+\ve\DD_1+\ve^2\DD_2,\\
&X_{\mathcal{T}}=(f_1+\ve{f}_2+\ve^2{f_3})\mathcal{T}\CU_{\ve}{f_3}
+f_3\mathcal{T}\CU_{\ve}{(f_1+\ve{f}_2+\ve^2{f_3})}+f_2\mathcal{T}\CU_{\ve}{f_2}+f_3\mathcal{T}\frac{\CU_{\ve}-{U_0}}{\ve}f_0
\\
&\quad+f_1\mathcal{T}\Big(\frac{\CU_{\ve}-{U_0}-\ve{U}_1-\ve^2U_2}{\ve^3}f_0+\frac{\CU_{\ve}-{U_0}-\ve{U}_1}{\ve^2}f_1
+\frac{\CU_{\ve}-{U_0}}{\ve}f_2\Big)\\
&\quad+f_0\mathcal{T}\Big(\frac{\CU_{\ve}-{U_0}-\ve{U}_1-\ve^2U_2-\ve^3U_3}{\ve^4}f_0+
\frac{\CU_{\ve}-{U_0}-\ve{U}_1-\ve^2U_2}{\ve^3}f_1\\
&\quad+\frac{\CU_{\ve}-{U_0}-\ve{U}_1}{\ve^2}f_2+\frac{\CU_{\ve}-{U_0}}{\ve}f_3\Big)
+f_2\mathcal{T}\big(\frac{\CU_{\ve}-{U_0}-\ve{U}_1}{\ve^2}f_0+\frac{\CU_{\ve}-{U_0}}{\ve}f_1\big)
,\quad\text{for }\mathcal{T}=\CR\text{ or }\nabla,
\\
&L_1=-\Big(\frac{\partial{f_3}}{\partial{t}}+\vv_0\cdot\nabla{f_3}+\vv_1\cdot
\nabla({f_2}+\ve{f_3})+\vv_2\cdot\nabla(f_1+\ve{f}_2+\ve^2{f_3})\Big)
\nonumber\\
&\quad+\CR\cdot\Big(X_{\CR}-\mm\times
(\nabla\vv_0)^T\cdot\mm{f_3}-\mm\times(\nabla\vv_1)^T
\cdot\mm(f_2+\ve{f_3})-\mm\times(\nabla\vv_2)^T\cdot\mm{(f_1+\ve{f}_2+\ve^2{f_3})}\Big),\\
&L_2=\frac{1-\gamma}{2Re}\nabla\cdot\Big\{\DD_0:\langle\mm\mm\mm\mm{f_3}\rangle
+\DD_1:\langle\mm\mm\mm\mm({f}_2+\ve{f_3})\rangle
+\DD_2:\langle\mm\mm\mm\mm(f_1+\ve{f}_2+\ve^2{f_3})\rangle\Big\}
\nonumber\\&\quad-\frac{1-\gamma}{Re}\big\{\nabla\cdot\big\langle\mm\mm\times
X_{\CR}\big\rangle_1+\langle{X}_{\nabla}\rangle_1\big\}-\vv_2\cdot\nabla\vv_1-\vv_1\cdot\nabla\vv_2-\ve\vv_2\cdot\nabla\vv_2.
\end{align*}
Then we can deduce the equations of $(f_R^\ve,\vv_R^\ve)$(drop $\ve$
for the simplicity):
\begin{align*}
&\frac{\pa{f_R}}{\pa{t}}+\tv\cdot\nabla{f_R}
+\vv_R\cdot\nabla(f_0+\ve{\tf})+{\frac{1}{\ve}}\CA_{f_0}\CH^{\ve}_{f_0}{f_R}\nonumber\\
&=-\CR\cdot\Big(\mm\times(\nabla\tv)^T\cdot\mm{f_R}+\mm\times(\nabla\vv_R)^T\cdot\mm{(f_0
+\ve{\tf})}+\ve^3\mm\times(\nabla\vv_R)^T\cdot\mm{f_R}\Big)\nonumber\\
&\quad-\CR\cdot\Big(f_R\CR\CU_{\ve}\tf
+f_R\CR\frac{(\CU_\ve-{\CU})f_0}{\ve}+\tf\CR\CU_{\ve}{f_R}+\ve^2{f_R}\CR\CU_{\ve}{f_R}\Big)+L_1,\\
&\frac{\pa{\vv_R}}{\pa{t}}+\vv_R\cdot\nabla\tv+\tv\cdot\nabla\vv_R
+\ve^3\vv_R\cdot\nabla\vv_R+\nabla{p}-\frac{\gamma}{Re}\Delta\vv_R\nonumber\\
&=\frac{1-\gamma}{2Re}\nabla\cdot\Big\{{\tD}:\langle\mm\mm\mm\mm{f_R}\rangle_1
+\DD_R:\langle\mm\mm\mm\mm{(f_0+\ve{\tf})}\rangle_1+\ve^3\DD_R:\langle\mm\mm\mm\mm{f_R}\rangle_1\Big\}\nonumber\\
&\quad-\frac{1-\gamma}{Re}\nabla\cdot\Big\langle\mm\mm\times\Big({\frac{1}{\ve}}
f_0\CR\CH^{\ve}_{f_0}{f_R}+f_R\CR\CU_{\ve}\tf+f_R\CR\frac{(\CU_\ve-{\CU})f_0}{\ve}
+\tf\CR\CU_{\ve}{f_R}+\ve^2{f_R}\CR\CU_{\ve}{f_R}\Big)\Big\rangle_1\nonumber\\
&\quad-\frac{1-\gamma}{Re}\Big\langle\Big({\frac{1}{\ve}}f_0\nabla\CH^{\ve}_{f_0}{f_R}
+f_R\nabla\CU_{\ve}\tf+f_R\nabla\frac{(\CU_\ve-{\CU})f_0}{\ve}
+\tf\nabla\CU_{\ve}{f_R}+\ve^2{f_R}\nabla\CU_{\ve}{f_R}\Big)\Big\rangle_1+L_2.\non
\end{align*}
Here $\DD_R=\f12\big(\na\vv_R+(\na \vv_R)^T\big)$. We denote
$F_R=F_1+\cdots+F_6$ with
\begin{align*}
&F_1=-\vv_R\cdot\nabla(f_0+\ve{\tf}),\\
&F_2=-\CR\cdot\Big(\mm\times(\nabla\tv)^T\cdot\mm{f_R}+f_R\CR\CU_{\ve}\tf
+f_R\CR\frac{(\CU_\ve-\CU_0)f_0}{\ve}+\tf\CR\CU_{\ve}{f_R}\Big),\\
&F_3=-\ve\CR\cdot\big(\mm\times(\nabla\vv_R)^T\cdot\mm{{\tf})}\big),\quad
F_4=-\ve^3\CR\cdot\big(\mm\times(\nabla\vv_R)^T\cdot\mm f_R\big),\\
&F_5=-\ve^2\CR\cdot\big(f_R\CR\CU_\ve f_R\big),\quad
F_6=\ve^3\vv_R\cdot\nabla{f_R},
\end{align*}
and $G_R=G_1+\cdots+G_8$ with
\begin{align*}
&G_1=-\vv_R\cdot\nabla\tv-\tv\cdot\nabla\vv_R,\quad G_2=\frac{1-\gamma}{2Re}\nabla\cdot\big(
\ve\DD_R:\langle\mm\mm\mm\mm{{\tf}}\rangle_1\big),\\
&G_3=-\frac{1-\gamma}{Re}\nabla\cdot\Big\langle\mm\mm\times\big(\widetilde{\DD}:\langle\mm\mm\mm\mm{f_R}\rangle_1+f_R\CR\CU\tf+
f_R\CR\frac{(\CU_\ve-{\CU})f_0}{\ve}+\tf\CR\CU{f_R}\big)\Big\rangle_1,\\
&G_4=\frac{1-\gamma}{Re}\Big\langle{\frac{1}{\ve}}f_0\nabla\CH_{\ve}{f_R}
+f_R\nabla\CU_{\ve}\tf+f_R\nabla\frac{(\CU_\ve-{\CU})f_0}{\ve}+\tf\nabla\CU_{\ve}{f_R}\Big\rangle_1,\\
&G_5=\frac{1-\gamma}{2Re}\ve^3\nabla\cdot\big(
\DD_R:\langle\mm\mm\mm\mm{f_R}\rangle_1\big),\quad G_6=-\ve^3\vv_R\cdot\nabla\vv_R,\\
&G_7=-\ve^2\frac{1-\gamma}{Re}\nabla\cdot\big\langle\mm\mm\times\big({f_R}\CR\CU_\ve{f_R}\big)\big\rangle_1,\quad
G_8=-\ve^2\frac{1-\gamma}{Re}\big\langle
{f_R}\nabla\CU_{\ve}{f_R}\big\rangle_1.
\end{align*}
Then we rewrite the equations for $(f_R,\vv_R)$ as
\begin{align}
&\frac{\pa{f_R}}{\pa{t}}+\tv\cdot\nabla{f_R}
+\ve^3\vv_R\cdot\nabla{f_R}+{\frac{1}{\ve}}\CA_{f_0}\CH^\ve_{f_0}{f_R}\non\\
&\qquad\qquad=-\CR\cdot\big(\mm\times(\nabla\vv_R)^T\cdot\mm{f_0}\big)+F_R+L_1,\label{eq:error-L-f}\\
&\frac{\pa{\vv_R}}{\pa{t}}-\frac{\gamma}{Re}\Delta\vv_R
-\frac{1-\gamma}{2Re}\nabla\cdot\big(\DD_R:\langle\mm\mm\mm\mm{{f_0}}\rangle_1\big)+\nabla{p}\non\\
&\qquad\qquad=-\frac{1-\gamma}{Re}\nabla\cdot\Big\langle\mm\mm\times\big({
\frac{1}{\ve}}f_0\CR\CH^\ve_{f_0}{f_R}\big)\Big\rangle_1+G_R+L_2.
\label{eq:error-L-v}
\end{align}

\subsection{Some key estimates}

In this subsection, we mainly present a control for a singular term in the error estimates.
The proof is based on the lower bound inequality. Since we only have good lower bound for the part outside the Maier-Saupe space,
we have to analyze the nonlinear interactions between the part inside the Maier-Saupe space
and the part outside the Maier-Saupe space. Throughout this section, we will repeatedly use the notations from Section 5.
Due to the assumption (\ref{ass:degenerate}), we can construct a global coordinate transformation so that all results from Section 5
can be applied.

\begin{Proposition}\label{prop:nonest-key}
For any $\delta>0$, there exists $C=C(\delta)$ such that for any
$f\in H^1(\Om\times\BS)$ with $\int_\BS{f}(\xx,\mm)\ud\mm=0$, there
holds
\begin{align*}
\frac{1}{\ve}\big\langle{f},
f\partial_t(\frac{1}{f_0})\big\rangle\le
\frac{\delta}{\ve^2}\big\langle{f_0}\CR\CH^\ve_{f_0}{f},\CR\CH^\ve_{f_0}{f}\big\rangle
+{C}\big(\frac{1}{\ve}\big\langle\CH^\ve_{f_0}{f},
f\big\rangle+\langle{f}, f\rangle\big).
\end{align*}
\end{Proposition}

We need the following lemmas.

\begin{Lemma}\label{lem:M}
It holds that
\begin{align*}
&\HM_{11}-\HM_{22}=\zeta_{a,2}\frac{(A_0-2A_2+A_4)}{2A_0},\\
&2\HM_{33}-\HM_{11}-\HM_{22}=3\zeta_0\frac{A_0A_4-A_2^2}{A_0^2}.
\end{align*}
\end{Lemma}

\no{\bf Proof.}\, Recall that $\HMM(\xx)=\int_\BS\hmm\hmm{f}\ud\hmm$. We get by direct
calculations that
\begin{align*}
\HM_{11}=&\frac{1}{2}\int\sin^2\htheta{a_0}\ud\hmm+\frac{1}{4}\int\sin^2\htheta{a_2}\ud\hmm
=-\zeta_0\frac{A_0A_4-A_2^2}{2A_0^2}+\zeta_{a,2}\frac{(A_0-2A_2+A_4)}{4A_0},\\
\HM_{22}=&\frac{1}{2}\int\sin^2\htheta{a_0}\ud\hmm-\frac{1}{4}\int\sin^2\htheta{a_2}\ud\hmm
=-\zeta_0\frac{A_0A_4-A_2^2}{2A_0^2}-\zeta_{a,2}\frac{(A_0-2A_2+A_4)}{4A_0},\\
\HM_{33}=&\int\cos^2\htheta{a_0}\ud\hmm
=\zeta_0\frac{A_0A_4-A_2^2}{A_0^2},\\
\HM_{12}=&\frac{1}{4}\int\sin^2\htheta{b_2}\ud\hmm
=\zeta_{b,2}\frac{(A_0-2A_2+A_4)}{4A_0},\\
\HM_{13}=&\frac{1}{2}\int\sin\htheta\cos\htheta{a_1}\ud\hmm
=\zeta_{a,1}\frac{A_2-A_4}{A_0},\\
\HM_{23}=&\frac{1}{2}\int\sin\htheta\cos\htheta{b_1}\ud\hmm
=\zeta_{b,1}\frac{A_2-A_4}{A_0}.
\end{align*}
Then the lemma follows.\ef

\begin{Lemma}\label{lem:H-ve}
There exists $c>0$ such that \beno
\big\langle\CH^\ve_{f_0}{f},f\big\rangle\ge{c}\langle f^\bot,f^\bot\rangle
+c\big\langle{M}_{kl}, M_{kl}-g_\ve*M_{kl}\big\rangle. \eeno
\end{Lemma}

\no{\bf Proof.}\,Noting that
\begin{align*}
\big\langle{\CH}^\ve_{f_0}{f}, f\big\rangle=&\langle\CH_{f_0}
f,f\rangle+\alpha\Big(\int_\Omega
\int_\BS\int_\BS(\mm\cdot\mm')^2f(\xx,\mm)f(\xx,\mm')\ud\mm'\ud\mm\ud\xx\\
&-\int_\Omega\int_\Omega\int_\BS\int_\BS(\mm\cdot\mm')^2g_\ve(\xx-\xx')
f_R(\xx,\mm)f_R(\xx',\mm')\ud\mm'\ud\xx'\ud\mm\ud\xx\Big)\\
=&\big\langle\CH_{f_0}f,f\big\rangle+\alpha\int_\Omega\MM(\xx):\big(\MM(\xx)-g_\ve*\MM(\xx)\big)\ud\xx,
\end{align*}
then the lemma follows from Proposition \ref{prop:H-positive}.\ef

\begin{Lemma}\label{lem:communicator}
We have \beno \frac{1}{\ve}\big\langle(f_1-g_\ve*f_1),
f_3f_2\big\rangle
\le{C}\Big(\frac{1}{\ve}\langle(f_1-g_\ve*f_1),f_1\rangle
+\frac{1}{\ve}\langle(f_2-g_\ve*f_2),f_2\rangle+\langle{f_2},f_2\rangle\Big).
\eeno where the constant $C$ depends on $\|f_3\|_{L^\infty}$ and
$\|\na f_3\|_{L^\infty}$.
\end{Lemma}

\no{\bf Proof.}\,We write
$f_1-g_\ve*f_1=\big(1-\chi(\sqrt{\ve}D)\big)^2f_1$, that is $\chi(\xi)=1-\sqrt{1-\hat g(\xi)}$. Hence,
\begin{align*}
\frac{1}{\ve}\big\langle(f_1-g_\ve*f_1), f_3f_2\big\rangle
=&\frac{1}{\ve}\big\langle(1-\chi(\sqrt{\ve}D)\big)^2f_1, f_3f_2\big\rangle\\
=&\frac{1}{\ve}\big\langle(1-\chi(\sqrt{\ve}D)\big)f_1, f_3(1-\chi(\sqrt{\ve}D)f_2\big\rangle\\
&+\frac{1}{\ve}\big\langle(1-\chi(\sqrt{\ve}D)\big)f_1, [f_3,
\chi(\sqrt{\ve}D)]f_2\big\rangle.
\end{align*}
Then the lemma follows from the commutator estimate
\beno
\|[f_3,
\chi(\sqrt{\ve}D)]f_2\|_{L^2}\le C\ve^\f12\|\na
f_3\|_{L^\infty}\|f_2\|_{L^2}.
\eeno
By a scaling argument, it suffices to prove the commutator estimate with $\ve=1$. Let $K_j(\xx)$
be the kernel associated with the Fourier multiplier $i(\pa_j\chi)(D)$
(It is easy to show that $K_j$ is a Calderon-Zygmund kernel). We have
\beno
\big[f_3,\chi(\sqrt{\ve}D)\big]f_2=\sum_{j=1}^3\int_{\mathbb{R}^3}K_j(\xx-\yy)\int_0^1\pa_j f_3(\tau\xx+(1-\tau)\yy)\ud\tau f_2(\yy)\ud\yy,
\eeno
whose $L^2$ norm is bounded by $\|\na f_3\|_{L^\infty}\|f_2\|_{L^2}$; see \cite{Christ} for example.\ef
\vspace{0.1cm}

Now we are ready to prove Proposition \ref{prop:nonest-key}.
With the notations in section 5, we decompose $f$ as
\beno
f=a_0(\xx,\htheta)+\sum_{k\ge1}\big(a_k(\xx,\htheta)\cos{k}\hvphi+b_k(\xx,\htheta)\sin{k}\hvphi\big).
\eeno Then we can get
\begin{align*}
\frac{1}{\ve}\big\langle{f}, f\partial_t(\frac{1}{f_0})\big\rangle
=&\frac{1}{\ve}\int_\Omega\int_\BS\partial_t\big(\frac{1}{f_0}\big)\Big(a_0(\xx,\htheta)
+\sum_{k\ge1}\big(a_k(\xx,\htheta)\cos{k}\hvphi+b_k(\xx,\htheta)\sin{k}\hvphi\big)\Big)^2\ud\hmm\ud\xx\nonumber\\
\le&{C}\int_\Omega\int_\BS\frac{1}{\ve^2}\Big(\gamma_0^2+\gamma_1^2+\beta_1^2
+\gamma_2^2+\beta_2^2+\sum_{k\ge3}(a_k^2+b_k^2)\Big)\ud\hmm\ud\xx\non\\
&+C\int_\Omega\int_\BS\big(\zeta_0^2+\zeta_{a,1}^2+\zeta_{a,2}^2+\zeta_{b,1}^2+\zeta_{b,2}^2\big)\ud\hmm\ud\xx\non\\\non
&+\frac{1}{\ve}\int_\Omega\int_\BS\partial_t\big(\frac{1}{f_0}\big)
\Big(\zeta_0(\cos^2\htheta-\frac{A_2}{A_0})+\sin\htheta\cos\htheta(\zeta_{a,1}\cos\hvphi+\zeta_{b,1}\sin\hvphi)
\\&\qquad\qquad+\sin^2\htheta(\zeta_{a,2}\cos2\hvphi+\zeta_{b,2}\sin2\hvphi)\Big)^2\ud\hmm\ud\xx.
\end{align*}
As $\partial_t\big(1/f_0\big)=-\partial_tf_0/f_0^2$ and $\partial_tf_0\in\mathrm{Ker}~\CG_{f_0}$, we may assume
that \beno
\partial_t(\frac{1}{f_0})=\frac{1}{f_0}(w_1(\xx)\cos\hvphi+w_2(\xx)\sin\hvphi)\sin\htheta\cos\htheta.
\eeno
We have
\begin{align}\nonumber
&\frac{1}{\ve}\int_\Omega\int_\BS\partial_t\big(\frac{1}{f_0}\big)\Big(\zeta_0(\cos^2\htheta-\frac{A_2}{A_0})
+\sin\htheta\cos\htheta(\zeta_{a,1}\cos\hvphi+\zeta_{b,1}\sin\hvphi)\non\\
&\qquad\qquad\qquad+\sin^2\htheta(\zeta_{a,2}\cos2\hvphi+\zeta_{b,2}\sin2\hvphi)\Big)^2\ud\hmm\ud\xx\non\\\nonumber
&=\frac{2}{\ve}\int_\Omega\int_\BS\partial_t\big(\frac{1}{f_0}\big)
\sin\htheta\cos\htheta\big(\zeta_{a,1}\cos\hvphi+\zeta_{b,1}\sin\hvphi\big)
\big(\zeta_0(\cos^2\htheta-\frac{A_2}{A_0})\non\\
&\qquad\qquad+\sin^2\htheta(\zeta_{a,2}\cos2\hvphi+\zeta_{b,2}\sin2\hvphi)\big)\ud\hmm\ud\xx\non\\\nonumber
&\quad+\frac{1}{\ve}\int_\Omega\int_\BS\partial_t\big(\frac{1}{f_0}\big)\big(\zeta_0(\cos^2\htheta-\frac{A_2}{A_0})
+\sin^2\htheta(\zeta_{a,2}\cos2\hvphi+\zeta_{b,2}\sin2\hvphi)\big)^2\ud\hmm\ud\xx\non\\\nonumber
&\le\frac{2}{\ve}\int_\Omega\int_\BS\partial_t\big(\frac{1}{f_0}\big)
\sin\htheta\cos\htheta\big(\zeta_{a,1}\cos\hvphi+\zeta_{b,1}\sin\hvphi\big)\big(\zeta_0(\cos^2\htheta-\frac{A_2}{A_0})\non\\
&\qquad\quad+\sin^2\htheta(\zeta_{a,2}\cos2\hvphi+\zeta_{b,2}\sin2\hvphi)\big)\ud\hmm\ud\xx
+\frac{C}{\ve}\int_\Omega(\zeta_0^2+\zeta_{a,2}^2+\zeta_{b,2}^2)\ud\xx,\label{eq:cross
term}
\end{align}
where we have used the fact
\begin{align*}
\int_{\BS}\frac{1}{f_0}(w_1(\xx)\cos\hvphi+w_2(\xx)\sin\hvphi)\sin\htheta\cos\htheta
\big(\sin\htheta\cos\htheta(\zeta_{a,1}\cos\hvphi+\zeta_{b,1}\sin\hvphi)\big)^2\ud\hmm=0.
\end{align*}
We get by Lemma \ref{lem:M} that
\begin{align*}
&2{\zeta}_{0}(\xx)-\alpha{N_{0}}-\alpha(2\HM_{33}-\HM_{22}-\HM_{11}-N_0)\\
&=2{\zeta}_{0}(\xx)-3\alpha\zeta_0(\xx)\frac{A_0A_4-A_2^2}{A_0^2}
=\frac{3A_2^2+2A_0A_2-5A_0A_4}{A_0(A_2-A_4)}\zeta_0(\xx).
\end{align*}
Note that the coefficient is positive in the front of $\zeta_0(\xx)$
if $f_0$ is a stable critical point. This implies that
\begin{align}\nonumber
\frac{1}{\ve}\int_\Omega{w}_1(\xx)\zeta_{a,1}(\xx)\zeta_0(\xx)\ud\xx
\le&\int_\Omega\frac{\delta}{\ve^2}\big(2{\zeta}_{0}(\xx)-\alpha{N}_0(\xx))^2
+C(\delta,\|w_1\|_{L^\infty})\zeta_{a,1}^2(\xx)\ud\xx\\
&+\frac{C}{\ve}\big\langle2\HM_{33}(\xx)-\HM_{11}(\xx)-\HM_{22}(\xx)-N_0(\xx),\zeta_{a,1}(\xx)w_1(\xx)\big\rangle.\non
\end{align}
On the other hand, we have \beno
&&\HM_{ij}-N_{ij}=A_{ki}A_{lj}M_{kl}-A_{ki}A_{lj}(g_\ve*M_{kl})
=A_{ki}A_{lj}(M_{kl}-g_\ve*M_{kl}),\\
&&\zeta_{a,1}=\frac{2}{\alpha}\HM_{13}=\frac{2}{\alpha}A_{k'1}A_{l'3}M_{k'l'},
\eeno from which,  Lemma \ref{lem:communicator} and Lemma
\ref{lem:H-ve}, it follows that
\begin{align*}
\frac{1}{\ve}\big\langle\HM_{ij}(\xx)-N_{ij}(\xx),
\zeta_{a,1}(\xx)w_1(\xx)\big\rangle
&=\frac{1}{\ve}\big\langle{M}_{kl}-g_\ve*M_{kl},\frac{2}{\alpha}A_{ki}A_{lj}A_{k'1}A_{l'3}M_{k'l'}(\xx)w_1(\xx)\big\rangle\\
&\le{C(\BA,w_1)}\big(\frac{1}{\ve}\langle{M}_{kl}-g_\ve*M_{kl},M_{kl}\rangle+\langle{M}_{kl}, M_{kl}\rangle\big)\\
&\le{C(\BA,w_1)}\big(\frac{1}{\ve}\langle\CH^\ve_{f_0}{f},f\rangle+\langle{f},
f\rangle\big).
\end{align*}
Thus, it follows from Proposition \ref{prop:lower bound} that
\begin{align}\nonumber
\frac{1}{\ve}\int_\Omega{w}_1(\xx)\zeta_{a,1}(\xx)\zeta_0(\xx)\ud\xx
&\le\int_\Omega\frac{\delta}{\ve^2}\big(2{\zeta}_{0}(\xx)-\alpha{N}_0(\xx))^2
+C(\delta,\|w_1\|_{L^\infty})\zeta_{a,1}^2(\xx)\ud\xx\\
&\quad+{C(\BA,w_1)}\big(\frac{1}{\ve}\langle\CH^\ve_{f_0}{f},f\rangle+\langle{f}, f\rangle\big)\non\\
&\le
C\big(\frac{\delta}{\ve^2}\big\langle{f_0}\CR\CH^\ve_{f_0}{f},\CR\CH^\ve_{f_0}{f}\big\rangle
+\frac{1}{\ve}\langle\CH^\ve_{f_0}{f},f\rangle+\langle{f},
f\rangle\big).\non
\end{align}
This gives the desired estimate for the term in (\ref{eq:cross
term}): \beno
\frac{2}{\ve}\int_\Omega\int_\BS\partial_t\big(\frac{1}{f_0}\big)
\sin\htheta\cos\htheta\zeta_{a,1}\cos\hvphi\zeta_0(\cos^2\htheta-\frac{A_2}{A_0})\ud\mm\ud\xx.
\eeno The other terms in (\ref{eq:cross term}) can be treated
similarly. We omit the details. \ef

The following lemma is used to control the other singular terms like
$\f 1 {\ve}\big\langle \widetilde{\vv}\cdot\na f_R, \CH^\ve_{f_0}f_R\big\rangle$
in the error estimates.

\begin{Lemma}\label{lem:refined-estimate}
We have
\begin{align}
\|\na f\|_{L^2}^2\le{C}\big(\big\langle\CH^\ve_{f_0}\na f,
\na f\big\rangle+\frac{1}{\ve}\big\langle\CH^\ve_{f_0}{f}, f\big\rangle\big).
\end{align}
\end{Lemma}

\no{\bf Proof.}\,Let us first claim that \ben\label{claim}
\|f\|_{L^2}^2\le{C}\big(\big\langle\CH_\ve{f},
f\big\rangle+\big\langle\MM[f], \MM[f]\big\rangle\big),
\een
where $\MM[f]=\int_\BS\mm\mm{f}\ud\mm$.
Due to the choice of $g$, we have \beno
|\xi\hat{f}(\xi)|^2\le{C}\big((1-\hat{g}(\ve\xi))|\xi\hat{f}(\xi)|^2+\frac{1}{\ve}(1-\hat{g}(\ve\xi))|\hat{f}(\xi)|^2\big).
\eeno This implies that \beno
\|\na f\|_{L^2}^2\le{C}\big(\langle f-g_\ve\ast\na {f},
\na f\rangle+\frac{1}{\ve}\langle f-g_\ve\ast f, f\rangle\big),
\eeno which along with (\ref{claim}) gives
\begin{align*}
\|\na f\|_{L^2}^2\le& {C}\big(\big\langle\CH^\ve_{f_0}\na {f},
\partial_if\big\rangle
+\langle\na \MM[f], \na \MM[f]\rangle\big)\\
\le &{C}\big(\big\langle\CH^\ve_{f_0}\na {f},
\na f\big\rangle+\frac{1}{\ve}\big\langle\CH^\ve_{f_0},f\big\rangle\big),
\end{align*}
where we have used
\begin{align*}
\big\langle\CH^\ve_{f_0}{f},f\big\rangle\ge&\alpha\int_\Omega\int_\BS\int_\BS(\mm\cdot\mm')^2f(\xx,\mm)\Big(f(\xx,\mm')
-\int_\Omega{g_\ve}(\xx-\xx')f(\xx',\mm')\ud\xx'\Big)\ud\mm'\ud\mm\ud\xx\\
=&\alpha\big\langle(1-g_\ve)\ast\MM[f],\MM[f]\big\rangle.
\end{align*}
To complete the proof, it remains to prove the claim. We write $f=f^\bot+f^\top$ with
$f^\top\in \mathrm{Ker}\CG_{f_0}$ and $f^\bot\in (\mathrm{Ker}\CG^*_{f_0})^\bot$.
By Proposition \ref{prop:H-positive}, we have
$$\|f^\bot\|_{L^2}^2\le C\big\langle\CH_{f_0}f,f\big\rangle\le C\big\langle\CH^\ve_{f_0}f,f\big\rangle.$$
While, from the proof of Lemma \ref{lem:M}, we know that
\beno
\|f^\top\|_{L^2}^2\le C\big(\big\langle \hat M_{13}, \hat M_{13}\big\rangle
+\big\langle \hat M_{23}, \hat M_{23}\big\rangle\big)\le
C\big\langle \hat\MM[f^\top], \hat \MM[f^\top]\big\rangle=C\big\langle \MM[f^\top], \MM[f^\top]\big\rangle.
\eeno
This implies (\ref{claim}). \ef

In the nonlinear estimates, we will frequently use the following basic lemmas.

\begin{Lemma}\label{lem:nonest-1}
It holds that
\begin{align*}
\int_\Omega\int_\BS{M_1}(\xx){M_2(\xx,\mm)M_3(\xx,\mm)}\ud\mm\ud\xx&\le {C}
\|M_1\|_{L^2}\|M_2\|_{H^{0,2}}\|M_3\|_{L^2},\\
\int_\Omega\int_\BS{M_1}(\xx){M_2(\xx,\mm)M_3(\xx,\mm)}\ud\mm\ud\xx&\le {C}
\|M_1\|_{H^1}\|M_2\|_{H^{0,1}}\|M_3\|_{L^2},\\
\int_\Omega\int_\BS{M_1}(\xx){M_2(\xx,\mm)M_3(\xx,\mm)}\ud\mm\ud\xx&\le {C}
\|\na M_1\|_{H^1}\|M_2\|_{L^2}\|M_3\|_{L^2}.
\end{align*}
\end{Lemma}

\no{\bf Proof.}\,By H\"{o}lder inequality and Sobolev embedding, we get
\begin{align*}
&\int_\Omega\int_\BS{M_1}(\xx){M_2(\xx,\mm)M_3(\xx,\mm)}\ud\mm\ud\xx\\
&\le\|M_1\|_{L^2}\|M_2\|_{L_x^{\infty}L_m^2}\|M_3(\xx,\mm)\|_{L^2}
\le\|M_1\|_{L^2}\|M_2\|_{L_m^2L_x^{\infty}}\|M_3\|_{L^2}\\
&\le{C}\|M_1\|_{L^2}\|M_2\|_{H^{0,2}}\|M_3\|_{L^2}.
\end{align*}
The other two inequalities can be proved similarly. \ef

The following Bernstein type lemma is a direct consequence of Young's inequality.
\begin{Lemma}\label{lem:convolution}
Let $k\ge 0$ be an integer and $p\ge 2$. Then it holds that
\begin{align*}
&\|\na^k\CU_\ve{f}\|_{L^{p}(\Omega\times\BS)}\le{C}\ve^{-3/2(1/2-1/p)-k/2}\|{f}\|_{L^{2}(\Omega\times\BS)}.
\end{align*}
\end{Lemma}

\subsection{Error estimates}
Let us first explain how to choose a suitable energy functional.
It's helpful to look at the following toy model for $(f_R,\vv_R)$:
\begin{align*}
&\frac{\pa{f_R}}{\pa{t}}+{\frac{1}{\ve}}\CA_{f_0}\CH^\ve_{f_0}{f_R}
=-\CR\cdot\big(\mm\times(\nabla\vv_R)^T\cdot\mm{f_0}\big),\\
&\frac{\pa{\vv_R}}{\pa{t}}-\frac{\gamma}{Re}\Delta\vv_R+\nabla{p}
=-\frac{1-\gamma}{Re}\nabla\cdot\big\langle\mm\mm\times\big({
\frac{1}{\ve}}f_0\CR\CH^\ve_{f_0}{f_R}\big)\big\rangle_1.
\end{align*}
Compared with the homogeneous case, the new difficulty is caused by
the singular term $\frac{1}{\ve}\nabla\cdot\big\langle\mm\mm\times\big(
f_0\CR\CH^\ve_{f_0}{f_R}\big)\big\rangle_1$.
To deal with it, it is natural to introduce the energy functional
\beno
\frac{1}{\ve}\big\langle{f}_R,\CH^\ve_{f_0}{f}_R\big\rangle
+\frac{Re}{1-\gamma}\langle{\vv}_R, {\vv}_R\rangle,
\eeno
since we have the following important observation:
\beno
\big\langle\mm\times(\nabla\vv_R)^T\cdot\mm{f_0},\CR{\CH}^\ve_{f_0} f_R\big\rangle+
\big\langle\langle\mm\mm\times{f_0\CR\CH^\ve_{f_0}{f_R}}\rangle_1,\nabla\vv_R\big\rangle=0.
\eeno
However, $\big\langle{f}_R,\CH_{f_0}^\ve{f}_R\big\rangle$ does not give a control for the part of $f_R$ inside the kernel.
To have a control for the part inside the kernel,
we need to introduce another functional $\big\langle{f}_R,\CA^{-1}_{f_0}f_R\big\rangle$ similar to
the homogeneous case. So, the suitable energy functional for the toy model should be
\beno
\big\langle{f}_R,\CA^{-1}_{f_0}f_R\big\rangle+\frac{1}{\ve}\big\langle{f}_R,\CH^\ve_{f_0}{f}_R\big\rangle
+\frac{Re}{1-\gamma}\langle{\vv}_R, {\vv}_R\rangle.
\eeno
However, if we take $\frac{1}{\ve}\f d {dt}\big\langle{f}_R,\CH^\ve_{f_0}{f}_R\big\rangle$, a new singular term
$\f1 \ve\big\langle{f_R}, \pa_t(\f1{f_0}){f_R}\big\rangle$ will appear.
Since there is no any decay in $\ve$ for the part of $f_R$ inside the kernel, this term seems to have the order
of $O(1/\ve)$(\textbf{Very singular!}). Surprisingly, by analyzing the nonlinear interactions for
$\big\langle{f_R}, \pa_t(\f1{f_0}){f_R}\big\rangle$ and using the lower bound inequality, we find that
it is bounded.

In order to control the nonlinear terms, we also need to introduce a
higher order analogous of the energy functional, whose choice is also very subtle.
In all, our energy functional takes the form
\begin{align*}
\Ee(t)\eqdefa&
\big\langle{f}_R,\CA^{-1}_{f_0}f_R\big\rangle+\frac{1}{\ve}\big\langle{f}_R,\CH^\ve_{f_0}{f}_R\big\rangle
+\frac{Re}{1-\gamma}\langle{\vv}_R, {\vv}_R\rangle\\
+&C_1\ve\big\langle\nabla{f}_R, \CA^{-1}_{f_0}\nabla{f}_R\big\rangle
+C_2\Big[\ve\big\langle\nabla{f}_R,\CH^\ve_{f_0}\nabla{f}_R\big\rangle
+\frac{Re}{1-\gamma}\ve^2\langle\nabla{\vv}_R,\nabla{\vv}_R\rangle\Big]\\
+&\ve^3\big\langle\Delta{f}_R, \CA^{-1}_{f_0}\Delta{f}_R\big\rangle
+C_3\Big[{\ve^3}\big\langle\Delta{f}_R,\CH^\ve_{f_0}\Delta{f}_R\big\rangle
+\frac{Re}{1-\gamma}\ve^4\big\langle\Delta{\vv}_R, \Delta{\vv}_R\big\rangle\Big],\\
\Fe(t)\eqdefa&\frac{1}{\ve}\big\langle{f}_R,\CH^\ve_{f_0}{f}_R\big\rangle+
\frac{1}{\ve^2}\big\langle\CH^\ve_{f_0}{f}_R,
\CA\CH^\ve_{f_0}{f}_R\big\rangle
+\frac{\gamma}{1-\gamma}\big\langle\nabla{\vv}_R, \nabla{\vv}_R\big\rangle\\
&+C_1\big\langle\nabla{f}_R, \CH^\ve_{f_0}\nabla{f}_R\big\rangle+
C_2\Big[\big\langle\CH^\ve_{f_0}\nabla{f}_R,
\CA\CH^\ve_{f_0}\nabla{f}_R\big\rangle
+\ve^2\frac{\gamma}{1-\gamma}\big\langle\nabla^2{\vv}_R, \nabla^2{\vv}_R\big\rangle\Big]\\
&+{\ve^2}\big\langle\Delta{f}_R,\CH^\ve_{f_0}\Delta{f}_R\big\rangle
+C_3\Big[\ve^2\big\langle\CH^\ve_{f_0}\Delta{f}_R,
\CA\CH^\ve_{f_0}\Delta{f}_R\big\rangle
+\ve^4\frac{\gamma}{1-\gamma}\big\langle\nabla\Delta{\vv}_R,
\nabla\Delta{\vv}_R\big\rangle\Big].
\end{align*}
Here the constants $C_1, C_2$ and $C_3$ bigger than one will
be determined later.

\begin{Proposition}\label{prop:error}
There exist $c_1>0$ and $\ve_0>0$  such that for any $\ve\in (0,\ve_0)$ and $t\in
[0,T]$, there holds
\begin{align*}
\frac{d}{d{t}}\Ee(t)+c_1\Fe(t)\le
 C\big(1+\Ee(t)+\ve^{1/4}\Ee(t)^{3/2}+\ve \Ee(t)^2\big)+C\ve^{3/2}\Ee(t)^{1/2}\Fe(t).
\end{align*}
where the constant $C$ depends on
$\|f_i\|_{L^\infty(0,T;H^3(\Om\times\BS))}$$(i=0, 1,2,3)$,
$\|\vv_i\|_{L^\infty(0,T;H^3(\Om))}(i=0,1,2)$.
\end{Proposition}

\no{\bf Proof.}\,From the definition of $\Ee(t),\Fe(t)$ and Lemma
\ref{lem:H-ve}, it is easy to see that
\begin{align*}
&\|f_R\|_{L^2}^2+C_1\|\ve^{1/2}{f}_R\|_{H^{0,1}}^2+\|\ve^{3/2}{f}_R\|_{H^{0,2}}^2\\
&\qquad\qquad+\|\vv_R\|_{L^2}^2+C_2\|\ve\vv_R\|_{H^1}^2+C_3\|\ve^2\vv_R\|_{H^2}^2\le {C}\Ee(t),\\
&\|\nabla\vv_R\|_{L^2}^2+C_2\|\ve\nabla\vv_R\|_{H^1}^2+C_3\|\ve^2\nabla\vv_R\|_{H^2}^2\le{C}\Fe(t).
\end{align*}
And it is easy to show that \beno
\|L_1\|_{H^{0,2}}+\|L_2\|_{H^2}\le C. \eeno These facts will be
repeatedly used in the following calculations. For the simplicity of notations,
we denote $\CA=\CA_{f_0}$ and $\CH_\ve=\CH^\ve_{f_0}$ in what follows.
\vspace{0.1cm}

\no{\bf Step 1.}\,$L^2$ energy estimate\vspace{0.1cm}

Making $L^2(\Om\times \BS)$ inner product to (\ref{eq:error-L-f})
with $\CA^{-1}_{f_0}f_R$, we get
\begin{align*}
&\big\langle\frac{\pa}{\pa{t}}{f}_R,
\CA^{-1}f_R\big\rangle+\frac{1}{\ve}\big\langle\CH_\ve {f_R},
f_R\big\rangle
=-\big\langle\tv\cdot\nabla{f_R}, \CA^{-1}f_R\big\rangle\nonumber\\
&\qquad+\big\langle\mm\times(\nabla\vv_R)^T\cdot\mm{f_0},
\CR\CA^{-1}f_R\big\rangle+\big\langle
F_R+{L_1},\CA^{-1}f_R\big\rangle.
\end{align*}
By Lemma \ref{lem:nonest-1} and Lemma
\ref{lem:convolution}, we have
\begin{align*}
&\big\langle F_4, \CA^{-1}f_R\big\rangle\le{C}
\ve^{1/2}\|\ve^2\nabla\vv_R\|_{H^1}\|\ve^{1/2}f_R\|_{H^{0,1}}\|f_R\|_{L^2}
\le{C}\ve^{1/2}\Ee^{3/2}(t),\\
&\big\langle F_5, \CA^{-1}f_R\big\rangle
\le{C}\ve^{5/4}\|f_R\|_{L^2}^3\le{C}\ve^{5/4}\Ee(t)^{3/2},\\
&\big\langle F_6,\CA^{-1}f_R\big\rangle
\le{C}\ve\|f_R\|_{L^2}^2\|\ve^2\vv_R\|_{H^2}\le{C}\ve\Ee(t)^{3/2},
\end{align*}
and the other terms can be estimated as follows
\begin{align*}
&\big\langle\tv\cdot\nabla{f_R},\CA^{-1}f_R\big\rangle+\big\langle
F_2,\CA^{-1}f_R\big\rangle
\le{C}\|f_R\|_{L^2}^2\le{C}\Ee(t),\\
&\big\langle F_1, \CA^{-1}f_R\big\rangle
\le {C}\|\vv_R\|_{L^2}\|f_R\|_{L^2}^2\le{C}\Ee(t),\\
&\big\langle F_3, \CA^{-1}f_R\big\rangle\le{C}\ve\|\na
\vv_R\|_{L^2}\|f_R\|_{L^2}\le{C}\Ee(t),
\end{align*}
So, we get
\begin{align}\label{eq:energyL-L2-1}
\big\langle\frac{\pa{f_R}}{\pa{t}}, \CA^{-1}f_R\big\rangle
+{\frac{1}{\ve}}\big\langle\CH_\ve {f_R}, f_R\big\rangle \le
C\big(1+\Ee(t)+\ve^{1/2}\Ee(t)^{3/2}\big)+\delta\Fe(t).
\end{align}

Make $L^2(\Om\times \BS)$ inner product to (\ref{eq:error-L-f})
with $\CH_\ve  f_R$ to obatin
\begin{eqnarray*}
&&\frac{1}{\ve}\big\langle\frac{\partial}{\partial{t}}{f}_R,
\CH_\ve {f}_R\big\rangle+\frac{1}{\ve^2}\big\langle\CR\CH_\ve {f_R},
f_0\CR\CH_\ve {f}_R\big\rangle
=-\frac{1}{\ve}\big\langle\tv\cdot\nabla{f_R}, \CH_\ve {f}_R\big\rangle\nonumber\\
&&\qquad+\frac{1}{\ve}\big\langle\mm\times(\nabla\vv_R)^T
\cdot\mm{f_0},\CR{\CH}_\ve f_R\big\rangle+\frac{1}{\ve}\big\langle
F_R+{L_1},\CH_\ve {f}_R\big\rangle.
\end{eqnarray*}
By Lemma \ref{lem:nonest-1} and Lemma \ref{lem:convolution}, we have
\begin{align*}
&\f 1\ve\big\langle F_4+F_6, {\CH}_\ve
f_R\big\rangle\le{C}\ve^{1/2}\|{\ve}^{1/2}{f}_R\|_{H^{0,1}}\|\ve^2\vv_R\|_{H^2}\|\frac{1}{\ve}\CR\CH_\ve {f}_R\|_{L^2}
\le{C}\ve^{1/2}\Ee(t)\Fe(t)^{1/2},\\
&\f 1\ve\big\langle F_5, \CH_\ve {f}_R\big\rangle\le
{C}\ve^{1/4}\|f_R\|_{L^2}^2\|\CR\CH_\ve {f}_R\|_{L^2}\le{C}\ve^{5/4}\Ee(t)\Fe(t)^{1/2}.
\end{align*}
Noting that \beno \int_{\BS}\nabla(f_0+\ve{\tf})\ud\mm=0,\quad
\int_{\BS}N_1\ud\mm=0, \eeno we get by Poincar\'{e} inequality that
\begin{align*}
&\frac{1}{\ve}\big\langle F_1,\CH_\ve {f}_R\big\rangle\le{C}
\|\vv_R\|_{L^2}\|\frac{1}{\ve}\CR\CH_\ve {f}_R\|_{L^2}\le{C}\Ee(t)+\delta\Fe(t),\\
&\frac{1}{\ve}\big\langle{N}_1, \CH_\ve {f}_R\big\rangle\le{C}
\|\frac{1}{\ve}\CR\CH_\ve {f}_R\|_{L^2}\le{C}+\delta\Fe(t).
\end{align*}
We infer from Lemma \ref{lem:refined-estimate} that \beno
\frac{1}{\ve}\big\langle\tv\cdot\nabla{f_R},\CH_\ve {f}_R\big\rangle
\le{C}\|\nabla{f}_R\|_{L^2}\|\frac{1}{\ve}\CR\CH_\ve {f}_R\|_{L^2}
\le \f {C_0} {\sqrt{C_1}}\Fe(t)+C\Ee(t)+\delta\Fe(t).
\eeno
Here and what follows $C_0$ denotes a constant independent of $\delta$.
The other terms  are estimated as follows
\begin{align*}
&\frac{1}{\ve}\big\langle F_2, \CH_\ve {f}_R\big\rangle
\le{C}\|f_R\|_{L^2}\|\frac{1}{\ve}\CR\CH_\ve {f}_R\|_{L^2}\le {C}\Ee(t)+\delta\Fe(t),\nonumber\\
&\frac{1}{\ve}\big\langle F_3, \CH_\ve {f}_R\big\rangle
\le{C}\ve\|\nabla\vv_R\|_{L^2}\|\frac{1}{\ve}\CR\CH_\ve {f}_R\|_{L^2}\le{C}\ve\Fe(t),\nonumber\\
\end{align*}
Hence, we obtain
\begin{align}\label{eq:energyL-L2-2}
&\frac{1}{\ve}\big\langle\frac{\pa{f_R}}{\pa{t}},
\CH_\ve {f}_R\big\rangle\non
+{\frac{1}{\ve^2}}\big\langle{f}_0\CR\CH_\ve {f_R},\CR\CH_\ve {f}_R\big\rangle\\
&\le C\big(1+\Ee(t)\big)+{C}\ve^{1/2}\Fe(t)^{1/2}\Ee(t)+\big(\delta+\ve+\f {C_0} {\sqrt{C_1}}\big)\Fe(t)\non\\
&\qquad+\frac{1}{\ve}\big\langle\mm\times(\nabla\vv_R)^T\cdot \mm
f_0, \CR\CH_\ve {f}_R\big\rangle.
\end{align}

Make $L^2(\Om)$ inner product to (\ref{eq:error-L-v}) with
$\vv_R$ to get
\begin{align*}
&\frac{1}{2}\frac{d}{d{t}}\langle\vv_R,\vv_R\rangle
+\frac{\gamma}{Re}\langle\nabla\vv_R,\nabla\vv_R\rangle
+\frac{1-\gamma}{2Re}\big\langle \big(\DD_R:\langle\mm\mm\mm\mm{{f_0}}\rangle_1\big), \na\vv_R\big\rangle\\
&=\frac{1-\gamma}{Re}\Big\langle\mm\mm\times\big({
\frac{1}{\ve}}f_0\CR\CH_\ve {f_R}\big),
\nabla\vv_R\Big\rangle+\big\langle G_R+{L_2},\vv_R\big\rangle.
\end{align*}
Obviously, $\big\langle G_6, \vv_R\big\rangle=0$.
We have by Lemma \ref{lem:nonest-1} and Lemma \ref{lem:convolution} that
\begin{align*}
&\big\langle G_5, \vv_R\big\rangle\le{C}
\ve^{1/2}\|\ve^{1/2}f_R\|_{H^{0,1}}\|\ve^2\vv_R\|_{H^2}\|\nabla\vv_R\|_{L^2}
\le{C}\ve^{1/2}\Ee(t)\Fe(t)^{1/2},\\
&\big\langle G_7, \vv_R\big\rangle\le{C}
\ve^{1/4}\|f_R\|_{L^2}^2\|\ve\nabla\vv_R\|_{L^2}\le{C}\ve^{1/4}\Ee(t)^{3/2},\\
&\big\langle G_8, \vv_R\big\rangle\le{C}\ve^{3/4}
\|f_R\|_{L^2}^2\|\vv_R\|_{L^2}\le{C}\ve^{3/4}\Ee(t)^{3/2}.
\end{align*}
While,
$\big\langle\DD_R:\langle\mm\mm\mm\mm{f_0}\rangle_1,\nabla\vv_R\big\rangle=
\big\langle\DD_R:\langle\mm\mm\mm\mm{f_0}\rangle_1,\DD_R\big\rangle\ge0$.
The other terms  are estimated as follows
\begin{align*}
&\langle G_1,\vv_R\rangle\le{C}\|\vv_R\|_{L^2}^2\le{C}\Ee(t),\\
&\big\langle G_2, \vv_R\big\rangle\le{C}\ve\|\nabla\vv_R\|_{L^2}^2\le {C}\Ee(t)+\delta\Fe(t),\\
&\big\langle G_3, \vv_R\big\rangle
\le{C}\|f_R\|_{L^2}\|\nabla\vv_R\|_{L^2}\le{C}\Ee(t)+\delta\Fe(t),\\
&\big\langle G_4,\vv_R\big\rangle\le C\f1\ve\|\CR\CH_\ve  f_R\|_{L^2}\|\vv_R\|_{L^2}+{C}\|f_R\|_{L^2}\|\vv_R\|_{L^2}
\le{C}\Ee(t)+\delta\Fe(t).
\end{align*}
Thus, we obtain
\begin{align}\label{eq:energyL-L2-3}
&\f12\f d {dt}\big\langle\vv_R,\vv_R\big\rangle
+\frac{\gamma}{Re}\langle\nabla\vv_R,\nabla\vv_R\rangle\nonumber\\
&\le C\big(1+\Ee(t)+\ve^{1/4}\Ee(t)^{3/2}\big)+{C}\ve^{1/2}\Ee(t)\Fe(t)^{1/2}+\delta\Fe(t)\non\\
&\qquad+\big\langle\mm\mm\times
{\frac{1}{\ve}}f_0\CR\CH_\ve {f_R},\nabla\vv_R\big\rangle.
\end{align}

\no{\bf Step 2.}\,$H^1$ energy estimate\vspace{0.1cm}

Taking the derivative to (\ref{eq:error-L-f}) with respect to
$\xx_i$, then making $L^2(\Om\times\BS)$ inner product with
$\ve\CA^{-1}\pa_if_R$, we get
\begin{align*}
&\ve\big\langle\frac{\partial}{\partial{t}}\pai{f_R},\CA^{-1}\pai{f_R}
\big\rangle+\big\langle\pai{f_R},\CH_\ve \pai{f_R}\big\rangle\\
&=-\big\langle\pai{f}_0\CR\CH_\ve {f}_R,\CR\CA^{-1}\pai{f}_R\big\rangle
-\big\langle\pai(\frac{1}{f_0})f_R,\pai{f}_R\big\rangle
-\ve\big\langle\pai(\tv\cdot\nabla{f_R}), \CA^{-1}\pai f_R\big\rangle\nonumber\\
&\quad+\ve\big\langle\pai(\mm\times(\nabla\vv_R)^T\cdot\mm{f_0}),
\CR\CA^{-1}f_R\big\rangle +\ve\big\langle\pai F_R+\pai
{L_1},\CA^{-1}\pai f_R\big\rangle.
\end{align*}
By Lemma \ref{lem:nonest-1} and Lemma
\ref{lem:convolution}, we get
\begin{align*}
&\ve\big\langle \pai F_4, \CA^{-1}{\pai}f_R\big\rangle
\le{C}\ve\|\ve\na\vv_R\|_{H^1}\|\ve^{1/2}{f}_R\|_{H^{0,1}}\|\ve^{3/2}\pai{f}_R\|_{H^{0,1}}
\le{C}\ve\Fe(t)^{1/2}\Ee(t),\\
&\ve\big\langle\pai F_5, \CA^{-1}\pai{f}_R\big\rangle
\le C\ve^{5/4}\|f_R\|_{L^2}\|\ve^{1/2}f_R\|_{H^{0,1}}^2\le{C}\ve^{5/4}\Ee(t)^{3/2},\\
&\ve\big\langle\pai F_6,\CA^{-1}{\pai}f_R\big\rangle
\le{C}\ve\|\ve\na\vv_R\|_{H^1}\|\ve^{1/2}{f}_R\|_{H^{0,12}}\|\ve^{3/2}
\pai{f}_R\|_{H^{0,1}}\le{C}\ve\Fe(t)^{1/2}\Ee(t).
\end{align*}
It follows from Lemma \ref{lem:refined-estimate} that
\begin{align*}
&\ve\big\langle\pai\big(\mm\times(\nabla\vv_R)^T\cdot\mm{f_0}\big),
\CR\CA^{-1}\pai{f}_R\big\rangle+\ve\big\langle\pai F_3, \CA^{-1}\pai{f}_R\big\rangle\\
&\le{C}\|\ve\na\vv_R\|_{H^1}\|\pai{f}_R\|_{L^2}\le \f {C_0} {\sqrt{C_2}}\Fe(t)+\delta\Fe(t)+C\Ee(t).
\end{align*}
The other terms are estimated as follows
\begin{align*}
&\big\langle\pai{f}_0\CR\CH_\ve {f}_R,\CR\CA^{-1}\pai{f}_R\big\rangle
+\big\langle\pai(\frac{1}{f_0})f_R,\pai{f}_R\big\rangle
\le{C}\|f_R\|_{L^2}\|\pai{f}_R\|_{L^2}\le{C}\Ee(t)+\delta\Fe(t),\\
&-\ve\big\langle\pai\big(\tv\cdot\nabla{f_R}\big),\CA^{-1}{\pai}f_R\big\rangle+
\ve\big\langle\pai F_2, \CA^{-1}{\pai}f_R\big\rangle
\le{C}\ve||{f}_R\|_{H^{0,1}}^2\le{C}\Ee(t),\\
&\ve\big\langle\pai F_1,\CA^{-1}\pai{f}_R
\big\rangle\le{C}\ve^{1/2}\|\vv_R\|_{H^1}\|\ve^{1/2}\pai{f}_R\|_{L^2}\le
\delta\Fe(t)+C\Ee(t).
\end{align*}
So, we get
\begin{align}\label{eq:energyL-H1-1}
&\big\langle\frac{\pa }{\pa{t}}\pai{f_R}, \CA^{-1}\pai
f_R\big\rangle
+{\frac{1}{\ve}}\big\langle\CH_\ve {\pai f_R}, \pai f_R\big\rangle\non\\
&\le
C\big(1+\Ee(t)+\ve^{5/4}\Ee(t)^{3/2}\big)+{C}\ve\Fe(t)^{1/2}\Ee(t)+\big(\delta+\f {C_0} {\sqrt{C_2}}\big)\Fe(t).
\end{align}

Taking the derivative to (\ref{eq:error-L-f}) with respect to
$\xx_i$, then making $L^2(\Om\times\BS)$ inner product with
$\ve\CH\pa_if_R$, we get
\begin{align*}
&\ve\big\langle\frac{\pa}{\pa{t}}{\partial_if_R},\CH_\ve {\pai}f_R\big\rangle+
\big\langle{f}_0\CR\CH_\ve \pai{f_R},\CR\CH_\ve {\pai}f_R\big\rangle\\
&=-\big\langle\pai{f}_0\CR\CU_\ve{f}_R+{f}_R\CR\pai(\CU_0{f_0}),\CR\CH_\ve \pai{f}_R\big\rangle
-\ve\big\langle\pai(\tv\cdot\nabla{f_R}),{\CH}_\ve\pai{f_R}\big\rangle\nonumber\\
&\quad+\ve\big\langle\pai(\mm\times(\nabla\vv_R)^T\cdot\mm{f_0}),\CR{\CH}_\ve\pai{f_R}\big\rangle
+\ve\big\langle\pai F_R+\pai{L_1},{\CH}_\ve\pai{f_R}\big\rangle.
\end{align*}
The first term on the right hand side is bounded by
\beno
\|f_R\|_{L^2}^2\|\CR\CH_\ve \pai{f}_R\|_{L^2}^2\le {C}\Ee(t)+\delta\Fe(t).
\eeno
By Lemma \ref{lem:nonest-1} and Lemma
\ref{lem:convolution}, we get
\begin{align*}
&\ve\big\langle \pai F_4, \CH_\ve {\pai}f_R\big\rangle
\le{C}\ve^{1/2}\|\ve^2\vv_R\|_{H^2}\|\ve^{3/2}{f}_R\|_{H^{0,2}}\|\CR\CH_\ve \pai{f}_R\|_{L^2}
\le{C}\ve^{1/2}\Ee(t)\Fe(t)^{1/2},\\
&\ve\big\langle\ \pai F_5, \CR\CH_\ve \pai{f}_R\big\rangle\le C\ve^{7/4}\|\ve^{1/2}{f}_R\|_{H^{0,1}}\|f_R\|_{L^2} \|\CR\CH_\ve \pai{f}_R\|_{L^2}
\le{C}\ve^{7/4}\Ee(t)\Fe(t)^{1/2},\\
&\ve\big\langle \pai F_6,\CH_\ve {\pai}f_R\big\rangle
\le{C}\ve^{1/2}\|\ve^2\vv_R\|_{H^2}\|\ve^{3/2}{f}_R\|_{H^{0,2}}\|\CR\CH_\ve \pai{f}_R\|_{L^2}
\le{C}\ve^{1/2}\Ee(t)\Fe(t)^{1/2}.
\end{align*}
and by Lemma \ref{lem:refined-estimate},
\begin{align*}
&\ve\big\langle\pai\big(\tv\cdot\nabla{f_R}\big),\CH_\ve {\pai}f_R\big\rangle\\
&=\ve\big\langle\pai\tv\cdot\nabla{f}_R,
\frac{\pai{f}_R}{f_0}\big\rangle
-\ve\big\langle\tv\cdot\nabla(\frac{1}{f_0})\pai{f}_R,\pai{f}_R\big\rangle
-\ve\big\langle\tv\cdot\nabla{f}_R,\pai(\CU_{\ve}{\pai}f_R)\big\rangle\\
&\le{C}\ve^{1/2}\|\nabla{f}_R\|_{L^2}\|\pai{f}_R\|_{L^2}
\le C\Ee(t)+\delta\Fe(t).
\end{align*}
The other terms are estimated as follows
\begin{align*}
&\ve\big\langle\pai F_1,\CH_\ve \pai{f}_R
\big\rangle\le{C}\|\ve\nabla\vv_R\|_{L^2}\|\CR\CH_\ve \pai{f}_R\|_{L^2}\le{C}\Ee(t)+\delta\Fe(t),\\
&\ve\big\langle\pai F_2, \CH_\ve {\pai}f_R\big\rangle\le{C}\ve^{1/2}\|\ve^\f12{f}_R\|_{H^{0,1}}\|\CR\CH_\ve \pai{f}_R\|_{L^2}^2
\le{C}\Ee(t)+\delta\Fe(t),\\
&\ve\big\langle\pai F_3, \CH_\ve \pai{f}_R\big\rangle \le
C\|\ve^2\vv_R\|_{H^2}\|\CR\CH_\ve \pai{f}_R\|_{L^2}\le{C}\Ee(t)+\delta\Fe(t),\\
&\ve\big\langle\mm\times(\nabla\vv_R)^T\cdot\mm{\pai{f}_0},\CR\CH_\ve {\pai}f_R\big\rangle
\le{C}\|\ve\vv_R\|_{H^1}\|\CR\CH_\ve {\pai}f_R\|_{L^2}\le{C}\Ee(t)+\delta\Fe(t).
\end{align*}
Thus, we obtain
\begin{align}
&\ve\big\langle\frac{\pa}{\pa{t}}{\partial_if_R},\CH_\ve {\pai}f_R\big\rangle+
\big\langle{f}_0\CR\CH_\ve \pai{f_R},\CR\CH_\ve {\pai}f_R\big\rangle\non\\
&\le {C}\big(1+\Ee(t)\big)+{C}\ve^{1/2}\Ee(t)\Fe(t)^{1/2}+\delta\Fe(t)\non\\
&\qquad\quad+\frac{1}{\ve}\big\langle\mm\times(\nabla\pai\vv_R)^T\cdot \mm f_0,
\CR\CH_\ve \pai{f}_R\big\rangle.
\end{align}

Making $L^2(\Om)$ inner product to (\ref{eq:error-L-v}) with
$\ve^2\pa_i^2\vv_R$, we get
\begin{align*}
&\frac{\ve^2}{2}\frac{d}{d{t}}\big\langle\pai\vv_R,\pai\vv_R\big\rangle
+\frac{\gamma}{Re}\ve^2\big\langle\nabla\pai\vv_R,\nabla\pai\vv_R\big\rangle+
\frac{1-\gamma}{2Re}\ve^2\big\langle\pai \big(\DD_R:\langle\mm\mm\mm\mm{{f_0}}\rangle_1\big), \na\pai\vv_R\big\rangle\\
&=\frac{1-\gamma}{Re}\ve\big\langle\mm\mm\times\pai\big(f_0\CR\CH_\ve {f_R}\big), \nabla\pai\vv_R\big\rangle
+\ve^2\big\langle \pai G_R+\pai{L_2},\pai\vv_R\big\rangle.
\end{align*}
First of all, we know that $\big\langle\big(\pai\DD_R:\langle\mm\mm\mm\mm{{f_0}}\rangle_1\big), \na\pai\vv_R\big\rangle\ge 0$ and
\beno
\ve^2\big\langle\big(\DD_R:\pai\langle\mm\mm\mm\mm{{f_0}}\rangle_1\big), \na\pai\vv_R\big\rangle
\le C\|\ve\vv_R\|_{H^1}\|\ve\na\pai \vv_R\|_{L^2}\le C\Ee(t)+\delta\Fe(t).
\eeno
By Lemma \ref{lem:nonest-1} and Lemma
\ref{lem:convolution}, we get
\begin{align*}
&\ve^2\big\langle\pai G_5, \pai\vv_R\big\rangle\le
C\ve^{1/2}\|\ve^{1/2}f_R\|_{H^{0,1}}\|\ve^2\vv_R\|_{H^2}\|\ve^2\na\vv_R\|_{H^2}\le C\ve^\f12\Ee(t)\Fe(t)^{1/2},\\
&\ve^2\big\langle\pai G_6, \pai\vv_R\big\rangle\le C\ve\|\ve^2\vv_R\|_{H^2}^2\|\na\vv_R\|_{L^2}\le C\ve\Ee(t)\Fe(t)^{1/2},\\
&\ve^2\big\langle\pai G_7, \pai\vv_R\big\rangle
\le{C}\ve^{3/4}\|\ve^{1/2}{f}_R\|_{H^{0,1}}\|f_R\|_{L^2}\|\ve^2\vv_R\|_{H^2}
\le{C}\ve^{3/4}\Ee(t)^{3/2},\nonumber\\
&\ve^2\big\langle\pai G_8,\pai\vv_R\big\rangle
\le{C}\ve\|{f}_R\|_{L^2}\|\ve^{1/2}{f}_R\|_{H^{0,1}}\|\ve^2\vv_R\|_{H^2}\le
C\ve^{3/4}\Ee(t)^{3/2}.
\end{align*}
And the other terms are estimates as
\begin{align*}
&\ve^2\big\langle\pai G_1,\pai\vv_R\big\rangle
\le{C}\|\ve\vv_R\|^2_{H^1}\le{C}\Ee(t),\\
&\ve^2\big\langle\pai G_2, \pai\vv_R\big\rangle
\le{C}\|\ve^2\vv_R\|_{H^2}\|\ve\Delta \vv_R\|_{L^2}\le {C}\ve\Ee(t)+\delta\Fe(t),,\\
&\ve^2\big\langle\pai G_3, \pai\vv_R\big\rangle
\le{C}\ve^{1/2}\|\ve^{1/2}{f}_R\|_{H^{0,1}}\|\ve\nabla\pai\vv_R\|_{L^2}\le {C}\Ee(t)+\delta\Fe(t),\\
&\ve^2\big\langle\pai G_4, \pai\vv_R\big\rangle
\le {C}\|f_R\|_{L^2}\|\ve\Delta\vv_R\|_{L^2}^2\le C\Ee(t)+\delta\Fe(t).\nonumber
\end{align*}
So, we get
\begin{align}
&\ve^2\big\langle\frac{\pa}{\pa{t}}{\pai\vv_R},\pai\vv_R\big\rangle
+\ve^2\frac{Re}{1-\gamma}\langle\nabla\pai\vv_R,\nabla\pai\vv_R\rangle\non\\
&\le C\big(1+\Ee(t)+\ve^{3/4}\Ee(t)^{3/2}\big)+C\ve^\f12\Ee(t)\Fe(t)^{1/2}+\delta\Fe(t)\non\\
&\qquad+\big\langle\mm\mm\times {\frac{1}{\ve}}f_0\CR\CH_\ve {\pai
f_R},\nabla\pai\vv_R\big\rangle.
\end{align}

\no{\bf Step 3.}\,$H^2$ energy estimate\vspace{0.1cm}

Taking $\Delta$ to (\ref{eq:error-L-f}), then making
$L^2(\Om\times\BS)$ inner product with $\ve^3\CA^{-1}\Delta f_R$, we
get
\begin{align*}
&\ve^3\big\langle\frac{\partial}{\partial{t}}\Delta{f_R},\CA^{-1}
\Delta{f_R}\big\rangle+\ve^2\big\langle\Delta{f_R},\CH_\ve \Delta{f_R}\big\rangle\\
&=\ve^2\langle{f_R}\CR\Delta{\CU_0{f_0}}+\Delta{f_0}\CR\CU_\ve{f_R}+2\pai{f_0}\CR\CU_\ve\pai{f_R}
+2\pai{f_R}\CR\CU_0\pai{f_0},\CR\CA^{-1}\Delta{f_R}\rangle\nonumber\\
&\quad-\ve^3\big\langle\Delta(\tv\cdot\nabla{f_R}), \CA^{-1}\Delta f_R\big\rangle
+\ve^3\big\langle\Delta(\mm\times(\nabla\vv_R)^T\cdot\mm{f_0}), \CR\CA^{-1}\Delta f_R\big\rangle\\
&\quad+\ve^3\big\langle \Delta F_R+\Delta {L_1},\CA^{-1}\Delta f_R\big\rangle.
\end{align*}
The first term on the right hand side is bounded by
\beno
\|\ve^{3/2}\CR\CA^{-1}\Delta{f_R}\|_{L^2}^2+\|\ve^{1/2}\pai{f}_R\|_{L^2}^2+\ve\|f_R\|_{L^2}^2
\le{C}\Ee(t).
\eeno
By Lemma \ref{lem:nonest-1} and Lemma \ref{lem:convolution}, we have
\begin{align*}
&\ve^3\big\langle\Delta F_4, \CA^{-1}{\Delta}f_R\big\rangle
\le{C}\ve\|\ve^2\nabla\vv_R\|_{H^2}\|\ve^{3/2}\Delta{f}_R\|_{L^2}\|\ve^{3/2}{f}_R\|_{H^{0,2}}
\le{C}\ve\Ee(t)\Fe(t)^{1/2},\\
&\ve^3\langle\Delta F_5,
\CA^{-1}\Delta{f}_R\rangle
\le{C}\ve^{5/4}\|f_R\|_{L^2}\|f_R\|_{H^{0,2}}^2\le{C}\ve^{5/4}\Ee(t)^{3/2},\\
&\ve^3\big\langle \Delta F_6,\CA^{-1}{\Delta}f_R\big\rangle
\le{C}\ve\|\ve^2\nabla\vv_R\|_{H^2}\|\ve^{3/2}\Delta{f}_R\|_{L^2}\|\ve^{3/2}{f}_R\|_{H^{0,2}}\le{C}\ve\Ee(t)\Fe(t)^{1/2}.
\end{align*}
And by Lemma \ref{lem:refined-estimate}, the term $\ve^3\big\langle\Delta(F_1+F_2+F_3), \CA^{-1}\Delta{f_R}\big\rangle$ is bounded by
\beno
&&\ve^{1/2}\|\ve\vv_R\|_{H^2}\|\ve^{3/2}\Delta{f}_R\|_{L^2}+\|\ve^{3/2}{f}_R\|_{H^{0,2}}^2
+\|\ve^2\nabla\vv_R\|_{H^2}\|\ve\Delta{f}_R\|_{L^2}\\
&&\le C\Ee(t)+\big(\f {C_0} {\sqrt{C_3}}+\delta\big)\Fe(t).
\eeno
So, we get
\begin{align}\label{eq:energyL-H1-1}
&\big\langle\frac{\pa}{\pa{t}}\Delta{f_R}, \CA^{-1}\Delta
f_R\big\rangle
+{\frac{1}{\ve}}\big\langle\CH_\ve {\Delta f_R}, \Delta f_R\big\rangle\non\\
&\le
C\big(1+\Ee(t)+\ve^{5/4}\Ee(t)^{3/2}\big)+{C}\ve\Fe(t)^{1/2}\Ee(t)+\big(\delta+\f {C_0} {\sqrt{C_3}}\big)\Fe(t).
\end{align}

Taking $\Delta$ to (\ref{eq:error-L-f}), then making
$L^2(\Om\times\BS)$ inner product with $\ve^3\CH\Delta f_R$, we get
\begin{align*}
&\ve^3\langle\frac{\pa}{\pa{t}}{\partial_if_R},\CH_\ve {\Delta}f_R\rangle+
\ve^2\langle{f}_0\CR\CH_\ve \Delta{f_R},\CR\CH_\ve {\Delta}f_R\rangle\non\\
&=\ve^2\langle{f_R}\CR\Delta{\CU_0{f_0}}+\Delta{f_0}\CR\CU_\ve{f_R}+2\pai{f_0}\CR\CU_\ve\pai{f_R}
+2\pai{f_R}\CR\CU_0\pai{f_0},\CR\CH_\ve \Delta{f_R}\rangle\nonumber\\
&\quad+\ve^3\big\langle\Delta(\mm\times(\nabla\vv_R)^T\cdot\mm{f_0}),\CR{\CH}_\ve\Delta{f_R}\big\rangle
+\ve^3\big\langle\Delta F_R+\Delta{L_1},{\CH}_\ve\pai{f_R}\big\rangle.
\end{align*}
The first term on the right hand side is bounded by
\beno
\|\ve\CR\CH_\ve \Delta{f_R}\|_{L^2}\big(\ve^{1/2}\|\ve^{1/2}\pai{f}_R\|_{L^2}+\ve\|f_R\|_{L^2}\big)
\le \delta\Fe(t)+C\Ee(t)).
\eeno
By Lemma \ref{lem:nonest-1} and Lemma
\ref{lem:convolution}, we get
\begin{align*}
\ve^3\langle\Delta F_4,
\CH_\ve {\Delta}f_R\rangle
&\le{C}\ve^{3/2}\|\ve^2\nabla\vv_R\|_{H^2}\|\ve^{3/2}{f}_R\|_{H^{0,2}}
\|\ve\CR\CH_\ve \Delta{f}_R\|_{L^2}\le{C}\ve^{3/2}\Fe(t)\Ee(t)^{1/2},\\
\ve^3\langle \Delta F_5,
\CR\CH_\ve \Delta{f}_R\rangle&\le
C\ve^{7/4}\|f_R\|_{L^2}\|\ve^{3/2}{f}_R\|_{H^{0,2}}\|\ve\CR\CH_\ve \Delta{f}_R\|_{L^2}
\le{C}\ve^{7/4}\Ee(t)\Fe(t)^{1/2},
\end{align*}
and for $F_6$, we have
\begin{align*}
&\ve^3\big\langle \Delta F_6,\CH_\ve {\Delta}f_R\big\rangle
=\ve^6\big\langle\Delta\vv_R\cdot\nabla{f}_R,
\frac{\Delta{f}_R}{f_0}\big\rangle
+2\ve^6\big\langle\pai\vv_R\cdot\nabla\pai{f}_R,\frac{\Delta{f}_R}{f_0}\big\rangle\\
&\qquad\qquad+\ve^6\big\langle\vv_R\cdot\nabla{f}_R,\Delta(\CU_{\ve}{\Delta}f_R)\big\rangle
-{\f{\ve^6}2\big\langle\vv_R\cdot\nabla(\frac{1}{f_0})\Delta{f}_R,\Delta{f}_R\big\rangle}\\
&\quad\le{C}\ve\big(\|\ve^2\na\vv_R\|_{H^2}\|\ve^{3/2}\nabla{f}_R\|_{H^{0,1}}
\|\ve^{3/2}\Delta{f}_R\|_{L^2}\\
&\qquad\qquad+\ve^{7/4}\|\vv_R\|_{L^2}\|\ve^{1/2}\na{f}_R\|_{L^2}
\|\ve^{3/2}\Delta{f}_R\|_{L^2}+\|\ve^2\vv_R\|_{H^2}\|\ve^{3/2}\Delta{f}_R\|_{L^2}^2\big)\\
&\quad\le C\ve\Fe(t)^{1/2}\Ee(t)+{C}\ve\Ee(t)^{3/2}.
\end{align*}
And the term $\ve^3\big\langle\Delta(F_1+F_2+F_3), \CH_\ve\Delta{f_R}\big\rangle$ is bounded by
\beno
&&\|\ve^2\vv_R\|_{H^2}\|\ve\CR\CH_\ve \Delta{f}_R\|_{L^2}+
\ve^{1/2}\|\ve^{3/2}f_R\|_{H^{0,2}}\|\ve\CR\CH_\ve \Delta{f}_R\|_{L^2}\\
&&\qquad+\ve\|\ve^2\na\vv_R\|_{H^2}\|\ve\CR\CH_\ve \Delta{f}_R\|_{L^2}
\le C\Ee(t)+(\ve+\delta)\Fe(t),
\eeno
and by Lemma \ref{lem:convolution},
\begin{align*}
&\ve^3\big\langle\Delta\big(\tv\cdot\nabla{f_R}\big),\CH_\ve {\Delta}f_R\big\rangle
=\ve^3\big\langle\Delta\tv\cdot\nabla{f}_R,
\frac{\Delta{f}_R}{f_0}\big\rangle
+2\ve^3\big\langle\pai\tv\cdot\nabla\pai{f}_R,\frac{\Delta{f}_R}{f_0}\big\rangle\\
&\qquad\qquad+\ve^3\big\langle\tv\cdot\nabla{f}_R,\Delta(\CU_{\ve}{\Delta}f_R)\big\rangle
-\f{\ve^3}2\big\langle\tv\cdot\nabla(\frac{1}{f_0})\Delta{f}_R,\Delta{f}_R\big\rangle\\
&\le C\|\ve^{1/2}\nabla{f}_R\|_{L^2}\|\ve^{3/2}\Delta{f}_R\|^2_{L^2}+C\|\ve^{3/2}\Delta{f}_R\|^2_{L^2} \le C\Ee(t).
\end{align*}
So, we get
\begin{align}\label{eq:energy-H2-2}
&\ve^3\big\langle\frac{\partial}{\partial{t}}\Delta{f_R},\CH_\ve
\Delta{f_R}\big\rangle+\ve^2\big\langle\CR\CH_\ve\Delta{f_R},f_0\CR\CH_\ve\Delta{f_R}\big\rangle\non\\
&\quad\le C\big(1+\Ee(t)+\ve\Ee(t)^{3/2}\big)+C\ve\Ee(t)\Fe(t)^{1/2}+C\big(\delta+\ve^{3/2}\Ee(t)^{1/2}\big)\Fe(t)\non\\
&\qquad\qquad+\ve^3\big\langle\mm\times(\nabla\Delta\vv_R)^T\cdot\mm{f_0},
\CR\CH_\ve \Delta{f_R})\big\rangle.
\end{align}

Making $L^2(\Om)$ inner product to (\ref{eq:error-L-v}) with
$\ve^4\Delta^2\vv_R$, we get
\begin{align*}
&\frac{\ve^4}{2}\frac{d}{d{t}}\big\langle\Delta\vv_R,\Delta\vv_R\big\rangle
+\frac{\gamma}{Re}\ve^4\big\langle\nabla\Delta\vv_R,\nabla\Delta\vv_R\big\rangle\\
&=-\frac{1-\gamma}{2Re}\ve^4\big\langle\Delta \big(\DD_R:\langle\mm\mm\mm\mm{{f_0}}\rangle_1\big), \na\Delta\vv_R\big\rangle\\
&\quad+\frac{1-\gamma}{Re}\ve^3\big\langle\mm\mm\times\Delta\big(f_0\CR\CH_\ve {f_R}\big), \nabla\Delta\vv_R\big\rangle
+\ve^4\big\langle \Delta G_R+\Delta{L_2},\Delta\vv_R\big\rangle.
\end{align*}
Again, $\big\langle\big(\Delta \DD_R:\langle\mm\mm\mm\mm{{f_0}}\rangle_1\big), \na\Delta\vv_R\big\rangle\ge 0$.
The other part of the first term is bounded by
\beno
\|\ve^2\vv_R\|_{H^2}\|\ve^2\na\Delta\vv_R\|_{L^2}\le C\Ee(t)+\delta\Fe(t).
\eeno
By Lemma \ref{lem:nonest-1} and Lemma
\ref{lem:convolution}, we get
\begin{align*}
&\ve^4\big\langle\Delta G_5,
\Delta\vv_R\big\rangle\le
C\ve^{3/2}\|\ve^{3/2}f_R\|_{H^{0,2}}\|\ve^2\nabla\Delta\vv_R\|_{L^2}\|\ve^2\nabla\vv_R\|_{H^2}
\le{C}\ve^{3/2}\Ee(t)^{1/2}\Fe(t),\\
&\ve^4\big\langle\Delta G_6,
\Delta\vv_R\big\rangle\le C\ve\|\ve^2\vv_R\|_{H^2}^2\|\ve^2\na\Delta\vv_R\|_{L^2}\le C\ve\Ee(t)\Fe(t)^{1/2},\\
&\ve^4\big\langle\Delta G_7, \Delta\vv_R\big\rangle
\le{C}\ve^{7/4}\|f_R\|_{L^2}\|\ve^{3/2}{f}_R\|_{H^{0,2}}\|\ve^2\nabla\Delta\vv_R\|_{L^2}
\le{C}\ve^{7/4}\Ee(t)\Fe(t)^{1/2},\nonumber\\
&\ve^4\big\langle\Delta G_8,\Delta\vv_R\big\rangle
\le{C}\ve^{9/4}\|\ve^{1/2}{f}_R\|_{H^{0,1}}\|f_R\|_{L^2}\|\ve^2\nabla\Delta\vv_R\|_{L^2}\le \ve\Ee(t)\Fe(t)^{1/2}.
\end{align*}
And the term $\ve^4\big\langle\Delta(G_1+G_2+G_3+G_4), \Delta{\vv_R}\big\rangle$ is bounded by
\beno
\|\ve^2\vv_R\|_{H^2}^2+\ve\|\ve^2\na \vv_R\|_{H^2}^2+
\ve^{1/2}\|\ve^{3/2}{f}_R\|_{H^{0,2}}\|\ve^2\nabla\Delta\vv_R\|_{L^2}
\le C\Ee(t)+(\delta+\ve)\Fe(t).
\eeno
Then we get
\begin{align}\label{eq:energy-H2-3}
&\ve^4\langle\frac{\pa{\Delta\vv_R}}{\pa{t}},\Delta\vv_R\rangle
+\f \ga {Re}\ve^4\langle\nabla\Delta\vv_R,\nabla\Delta\vv_R\rangle\nonumber\\
&\le C\big(1+\Ee(t)\big)+(\delta+\ve)\Fe(t)+C\ve\Ee(t)\Fe(t)^{1/2}\non\\
&\qquad+C\ve^{3/2}\Ee(t)^{1/2}\Fe(t)+\big\langle\mm\mm\times
{\frac{1}{\ve}}f_0\CR\CH_\ve {\Delta
f_R},\nabla\Delta\vv_R\big\rangle.
\end{align}

\no{\bf Step 4.} The closing of the energy estimates \vspace{0.1cm}

Noting that
\beno
\frac{1}{\ve}\big\langle\mm\times(\nabla\vv_R)^T\cdot\mm{f_0},\CR{\CH}_\ve f_R\big\rangle+
\frac{1}{\ve}\big\langle\langle\mm\mm\times{f_0\CR\CH_\ve{f_R}}\rangle_1,\nabla\vv_R\big\rangle=0,
\eeno
and then summing up (\ref{eq:energyL-L2-1})-(\ref{eq:energy-H2-3}), and taking $C_1$ big enough,
and then $C_2, C_3$ big enough, and finally taking $\delta$ small enough,
we infer that there exist $\ve>0$ and $c_1>0$ such that for any $\ve\in (0,\ve_0)$, there holds
\beno
&&\big\langle\frac{\pa}{\pa{t}}{f}_R, \CA^{-1}f_R\big\rangle+\frac{1}{\ve}\big\langle\frac{\partial}{\partial{t}}{f}_R, \CH_\ve{f}_R\big\rangle
+\f {Re} {1-\ga}\big\langle\frac{\pa}{\pa{t}}{\vv}_R, \vv_R\big\rangle\\
&&+C_1\ve\big\langle\frac{\partial}{\partial{t}}\na{f_R},\CA^{-1}\na{f_R}
\big\rangle+C_2\ve\big\langle\frac{\partial}{\partial{t}}\na{f_R},{\CH}_\ve\na{f_R}\big\rangle+
C_2\ve^2\f {Re} {1-\ga}\big\langle\frac{\pa}{\pa{t}}\na{\vv}_R, \na\vv_R\big\rangle\\
&&+\ve^3\big\langle\frac{\partial}{\partial{t}}\Delta{f_R},\CA^{-1}
\Delta{f_R}\big\rangle+C_3\ve^3\big\langle\frac{\partial}{\partial{t}}\Delta{f_R},\CH_\ve
\Delta{f_R}\big\rangle+C_3\ve^4
\f {Re} {1-\ga}\big\langle\frac{\pa}{\pa{t}}\Delta{\vv}_R, \Delta\vv_R\big\rangle+c_1\Fe(t)\\
&&\le C\big(1+\Ee(t)+\ve^{1/4}\Ee(t)^{3/2}+\ve \Ee(t)^2\big)+C\ve^{3/2}\Ee(t)^{1/2}\Fe(t).
\eeno
Now Proposition \ref{prop:nonest-key} implies that
\begin{align*}
&\frac{1}{2}\frac{d}{d{t}}\langle{f}_R,\CH_\ve {f}_R\rangle
\le\langle\frac{\partial{f}_R}{\partial{t}},\CH_\ve {f}_R\rangle
+\delta\Fe(t)+C\Ee(t),\\
&\frac{1}{2}\frac{d}{d{t}}\langle\pai{f}_R,\CH_\ve \pai{f}_R\rangle
\le\langle\frac{\partial}{\partial{t}}\pai{f}_R,\CH_\ve \pai{f}_R\rangle
+\delta\Fe(t)+C\Ee(t),\\
&\frac{1}{2}\frac{d}{d{t}}\langle\Delta{f}_R,\CH_\ve \Delta{f}_R\rangle
\le\langle\frac{\partial\Delta{f}_R}{\partial{t}},\CH_\ve \Delta{f}_R\rangle
+\delta\Fe(t)+C\Ee(t),
\end{align*}
and we have the trivial inequality
\beno
\f12\f d {dt}\big\langle f, \CA^{-1}f\rangle\le \big\langle \f \pa {\pa t}f, \CA^{-1}f\rangle
+C\|\CR\CA^{-1}f\|_{L^2}^2.
\eeno
Thus, we can deduce that
\beno
\frac{d}{d{t}}\Ee(t)+c_1\Fe(t)\le C\big(1+\Ee(t)+\ve^{1/4}\Ee(t)^{3/2}+\ve \Ee(t)^2\big)+C\ve^{3/2}\Ee(t)^{1/2}\Fe(t).
\eeno
This completes the proof of Proposition \ref{prop:error}.\ef

Now we are ready to prove Theorem \ref{thm:deborah-inhom-L}.
Given the initial data $\big(f_0^\ve,\vv_0^\ve\big)$, we can show by the energy method \cite{ZZ-SIAM} that
there exists $T_\ve>0$ and a unique solution $\big(f^{\ve}(\xx,\mm,t), \vv^\ve(\xx,t)\big)$ on $[0,T_\ve]$
to (\ref{eq:LCP-nonL-f})-(\ref{eq:LCP-nonL-v}) such that
\beno
f^\ve(t)-1\in C\big([0,T_\ve]; H^2(\Om\times \BS)\big),
\quad \vv^\ve(t)\in C\big([0,T_\ve]; H^2(\Om)\big)\cap L^2(0,T_\ve;H^3(\Om)).
\eeno
While, Proposition \ref{prop:error} tells us that
\beno
\frac{d}{d{t}}\Ee(t)+c_1\Fe(t)\le C\big(1+\Ee(t)+\ve^{1/4}\Ee(t)^{3/2}+\ve \Ee(t)^2\big)+C\ve^{3/2}\Ee(t)^{1/2}\Fe(t),
\eeno
for any $t\in [0,T_\ve]$. Due to the assumptions of Theorem \ref{thm:deborah-inhom-L}, we know that
$\Ee(0)\le C$. Thus, there exist $\ve_0>0$ depending on $T$ such that
for any $\ve\in (0,\ve_0)$ and $t\in [0,\min(T,T_\ve)]$, there holds
\beno
\Ee(t)+c_1\int_0^t\Fe(s)\ud s\le C.
\eeno
This in turn implies $T_\ve\ge T$ by a continuous argument.
Then Theorem \ref{thm:deborah-inhom-L} follows.\ef

\section{The dissipation of the Ericksen-Leslie energy}

Recall that the Ericksen-Leslie equation has the following energy law
\begin{eqnarray}
&&-\frac{\ud}{\ud{t}}\Big(\int_{\Omega}\frac{Re}{2(1-\gamma)}|\vv|^2\ud\xx+E_F\Big)\non\\
&&\quad=\int_{\Omega}\Big(\frac{\gamma}{1-\gamma}|\nabla\vv|^2+(\alpha_1+\frac{\gamma_2^2}{\gamma_1})|\DD:\nn\nn|^2
+\alpha_4\DD:\DD\nonumber\\
&&\qquad\qquad+(\alpha_5+\alpha_6-\frac{\gamma_2^2}{\gamma_1})|\DD\cdot\nn|^2
+\frac{1}{\gamma_1}|\nn\times\hh|^2\Big)\ud\xx.\quad\label{EL:energy}
\end{eqnarray}
Because the relations between six Leslie coefficients are unclear in Physics, whether the energy is dissipated
remains open. In \cite{Lin1}, Lin and Liu present some constrains on the Leslie coefficients to ensure
that the energy is dissipated. We will show that the energy (\ref{EL:energy}) is dissipated for the Ericksen-Leslie equation derived
from the Doi-Onsager equation. More precisely,

\begin{Theorem}\label{thm:energy}
If the Leslie coefficients are determined by (\ref{Leslie cofficients1}) and (\ref{Leslie cofficients2}), then there holds
\beno
(\alpha_1+\frac{\gamma_2^2}{\gamma_1})|\DD:\nn\nn|^2
+\alpha_4\DD:\DD+\big(\al_5+\al_6-\frac{\gamma_2^2}{\gamma_1}\big)|\DD\cdot\nn|^2\ge 0
\eeno
for any symmetric matrix $\DD$ and $\nn\in\BS$.
\end{Theorem}
\begin{Remark}
Recall that $\ga_1=S_2/\lambda$. By taking $\uu=\uu'$ in (\ref{eq:kernel-product}), we see that $\lambda>0$,
thus $\ga_1>0$.
\end{Remark}

Throughout this section, we denote by $f_0=h_{\eta,\nn}$ a critical point of $A[f]$.

\subsection{Some useful identities}
Recall that $S_2=\langle{P}_2(\mm\cdot\nn)\rangle_{f_0}$ and
$S_4=\langle{P}_4(\mm\cdot\nn)\rangle_{f_0}$, where
$P_k(x)$ is the $k$-th Legendre polynomial. We define
\beno
\MM^{(2)}=\langle\mm\mm\rangle_{f_0},\quad \MM^{(4)}=\langle\mm\mm\mm\mm\rangle_{f_0}.
\eeno

\begin{Lemma}\label{lem:M-tensor}
It holds that
\begin{align*}
&\MM^{(2)}=S_2\nn\nn+\frac{1-S_2}{3}\II,\\
&\MM^{(4)}_{\alpha\beta\gamma\mu}=S_4n_{\alpha}n_{\beta}n_{\gamma}n_{\mu}
+\frac{S_2-S_4}{7}\big(n_{\alpha}n_{\beta}\delta_{\gamma\mu}+
n_{\gamma}n_{\mu}\delta_{\alpha\beta}+n_{\alpha}n_{\gamma}\delta_{\beta\mu}
+n_{\beta}n_{\mu}\delta_{\alpha\gamma}\big)\nonumber\\
&\qquad+n_{\alpha}n_{\mu}\delta_{\beta\gamma}+
n_{\beta}n_{\gamma}\delta_{\alpha\mu})+\big(\frac{S_4}{35}-\frac{2S_2}{21}
+\frac{1}{15}\big)\big(\delta_{\alpha\beta}
\delta_{\gamma\mu}+\delta_{\alpha\gamma}\delta_{\beta\mu}
+\delta_{\alpha\mu}\delta_{\beta\gamma}\big).
\end{align*}
\end{Lemma}

The lemma is a direct consequence of Lemma \ref{lem:Q-tensor}. Especially, the lemma implies that

\begin{Lemma}\label{lem:tensor-product}
For any symmetric matrix $\DD$, there hold
\begin{align*}
\MM^{(2)}\cdot\DD=&S_2\nn(\DD\cdot\nn),\quad\DD\cdot\MM^{(2)}=S_2(\DD\cdot\nn)\nn;\\
\MM^{(4)}:\DD=&S_4\nn\nn(\DD:\nn\nn)+\frac{2(S_2-S_4)}{7}
\big((\DD\cdot\nn)\nn+\nn(\DD\cdot\nn))\nonumber\\
&+2\big(\frac{S_4}{35}-\frac{2S_2}{21}
+\frac{1}{15}\big)\DD+\frac{S_2-S_4}{7}\II(\DD:\nn\nn).
\end{align*}
\end{Lemma}

\begin{Lemma}\label{lem:identity}
For any symmetric constant matrix $\DD$, there holds
\begin{align*}
\big\langle\CR\cdot\big(\mm\times\DD\cdot\mm{f_0}\big),f\big\rangle
=\frac{1}{2}\DD:\int_\BS(\mm\mm-\frac{1}{3}\II)\CR\cdot(f_0\CR{f})\ud\mm.
\end{align*}
\end{Lemma}

\no{\bf Proof.}\,It is easy to show that for any vector field $\vv$ defined on $\BS$,
\begin{align}
\langle(\mm\mm-\frac{1}{3}\II)\CR\cdot(f\vv)\rangle_1=\langle(\mm\times\vv)\mm+\mm(\mm\times\vv)\rangle_f.\non
\end{align}
Applying it with  $\vv=\CR{g}$ and $\vv=\mm\times(\kappa\cdot\mm)$,
we deduce that
\begin{align}
\int_{\BS}(\mm\mm-\frac{1}{3}\II)\CR\cdot(f\CR{g})\ud\mm=\langle\mm\times\CR{g}\mm+\mm\mm\times\CR{g}\rangle_f.\non
\end{align}
Thus, we have
\begin{align}
&\big\langle\CR\cdot\big(\mm\times\DD\cdot\mm{f_0}\big),f\big\rangle
=\DD:\langle\mm(\mm\times\CR{f})\rangle_{f_0}\nonumber\\
&=\frac{1}{2}\DD:\big(\langle\mm(\mm\times\CR{f})\rangle_{f_0}+\langle(\mm\times\CR{f})\mm\rangle_{f_0}\big)\nonumber\\
&=\frac{1}{2}\DD:\int_\BS(\mm\mm-\frac{1}{3}\II)\CR\cdot(f_0\CR{f})\ud\mm.\non
\end{align}
The lemma follows. \ef

\begin{Lemma}\label{Lem:antisymmetric}
For any antisymmetric constant matrix $\BOm$, we have
\begin{align*}
\CR\cdot\big(\mm\times(\BOm\cdot\mm){f}_0\big)-(\nn\times(\BOm\cdot\nn))\cdot\CR{f_0}=0,
\end{align*}
\end{Lemma}
\no{\bf Proof.}\,
The lemma is a direct consequence of the following identities
\begin{align*}
&\CR\cdot\big(\mm\times(\BOm\cdot\mm)\big)=\CR_i(\epsilon^{ijk}m_j\Omega_{kl}m_l)=(\II-3\mm\mm):\BOm=0,\\
&(\mm\times(\BOm\cdot\mm))\cdot\CR{f_0}=(\mm\times(\BOm\cdot\mm))\cdot(\mm\times\nn)f_0'\\
&\quad=(\nn\times(\BOm\cdot\nn))\cdot(\mm\times\nn)f_0'=(\nn\times(\BOm\cdot\nn))\cdot\CR{f_0}.
\end{align*}
The proof is finished. \ef

\subsection{Projection operator and properties}

We denote by $\Pin$ the projection operator from $\CP_0(\BS)$ to $\mathrm{Ker}\CG_{f_0}$, and  denote by $\Pout$
the projection operator from $\CP_0(\BS)$ to $(\mathrm{Ker}\CG^*_{f_0})^{\bot}$.
Since $\mathrm{Ker}\CG_{f_0}$ is orthogonal to $(\mathrm{Ker}\CG^*_{f_0})^{\bot}$
under the inner product $\langle\cdot,\CA_{f_0}^{-1}(\cdot)\rangle$, we have
\begin{align*}
\langle{f},\CA^{-1}_{f_0}f\rangle=\langle\Pin{f},\CA^{-1}_{f_0}\Pin{f}\rangle+\langle\Pout{f},\CA^{-1}_{f_0}\Pout{f}\rangle.
\end{align*}
For any constant matrix $\kappa$, we define
\begin{align*}
\mathcal{K}(\kappa)=\Pin\big[\CR\cdot\big(\mm\times(\kappa\cdot\mm){f}_0\big)\big],
\quad\mathcal{L}(\kappa)=\Pout\big[\CR\cdot\big(\mm\times(\kappa\cdot\mm){f}_0\big)\big].
\end{align*}

\begin{Lemma}\label{Lem:ker-proj}
It holds that
\begin{align*}
\mathcal{K}(\kappa)=\big(\nn\times(\lambda\DD\cdot\nn-\BOm\cdot\nn)\big)\cdot\CR{f_0}.
\end{align*}
Here $\DD=\f12(\kappa+\kappa^T), \BOm=\f12(\kappa^T-\kappa)$.
\end{Lemma}

\no{\bf Proof.}\, By Theorem \ref{thm:G-kernel}, we may assume that
$$\Pin[\CR\cdot(\mm\times(\kappa\cdot\mm){f_0})]=\ww\cdot\CR{f}_0$$
for some vector $\ww$ with $\ww\bot\nn.$
Thus for all $\mathbf{\Theta}\cdot\CR{f}_0\in\mathrm{Ker}\CG_{f_0}$,
$$\big\langle\CR\cdot(\mm\times\kappa\cdot\mm{f_0}), \CA^{-1}_{f_0}(\mathbf{\Theta}\cdot\CR{f}_0)\big\rangle
=\big\langle\ww\cdot\CR{f}_0,\CA^{-1}_{f_0}(\mathbf{\Theta}\cdot\CR{f}_0)\big\rangle.$$
First we claim that
\begin{align}\label{eq:kernel-claim}
\langle\ww\cdot\CR{f}_0,\CA^{-1}_{f_0}(\mathbf{\Theta}\cdot\CR{f}_0)\rangle
=\mathbf{\Theta}\cdot\big(\nn\times\big({S}_2\DD\cdot\nn-\frac{S_2}{\lambda}\BOm\cdot\nn\big)\big).
\end{align}
Let $\uu$ and $\uu'$ be any vectors. By Proposition \ref{prop:kernel of Gstar}, we may write
\beno
\CA^{-1}_{f_0}(\uu'\cdot\CR{f}_0)=(u'_1\sin\phi-u'_2\cos\phi)g_0(\theta).
\eeno
Then we get by a direct computation that
\begin{align}
&\big\langle\uu\cdot\CR{f}_0,\CA^{-1}_{f_0}(\uu'\cdot\CR{f}_0)\big\rangle\nonumber\\
&=\int_\BS2\eta\sin\theta\cos\theta{f}_0(u_1\sin\phi-u_2\cos\phi)
(u'_1\sin\phi-u'_2\cos\phi)g(\theta)\ud\mm\nonumber\\\label{eq:kernel-product}
&=\frac{1}{2}(\uu\times\nn)(\uu'\times\nn)\int_\BS{f}_0\frac{\ud{u_0}}{\ud\theta}g(\theta)\ud\mm
=\frac{S_2}{\lambda}\big(\uu-(\uu\cdot\nn)\nn\big)\cdot\uu'.
\end{align}
Therefore, $\ww=\nn\times(\lambda\DD\cdot\nn-\BOm\cdot\nn)$.

Now, we prove (\ref{eq:kernel-claim}). By Lemma \ref{lem:identity}, Lemma \ref{Lem:antisymmetric} and (\ref{eq:kernel-product}), we have
\begin{align*}
&\big\langle\CR\cdot(\mm\times\kappa\cdot\mm{f_0}),\CA^{-1}(\mathbf{\Theta}\cdot\CR{f}_0)\big\rangle\\
&=\big\langle\CR\cdot(\mm\times\DD\cdot\mm{f_0}),\CA^{-1}(\mathbf{\Theta}\cdot\CR{f}_0)\big\rangle
-\big\langle\CR\cdot(\mm\times\BOm\cdot\mm{f_0}),\CA^{-1}(\mathbf{\Theta}\cdot\CR{f}_0)\big\rangle\\
&=-\frac{1}{2}\langle(\mm\mm:\DD), \mathbf{\Theta}\cdot\CR{f}_0\rangle
-\langle(\nn\times(\BOm\cdot\nn))\cdot\CR{f_0},\CA^{-1}(\mathbf{\Theta}\cdot\CR{f}_0)\rangle\\
&=\langle\mathbf{\Theta}\cdot(\mm\times(\DD\cdot\mm))\rangle_{f_0}
-\frac{S_2}{\lambda}\mathbf{\Theta}\cdot(\nn\times(\BOm\cdot\nn))\\
&=S_2\mathbf{\Theta}\cdot(\nn\times(\DD\cdot\nn))
-\frac{S_2}{\lambda}\mathbf{\Theta}\cdot(\nn\times(\BOm\cdot\nn)).
\end{align*}
The claim follows. \ef

\begin{Lemma}\label{Lem:ker-out-proj}
$\mathcal{L}(\BOm)=0$ for any antisymmetric matrix $\BOm$.
\end{Lemma}

\no{\bf Proof} This is equivalent to prove $\mathcal{K}(\BOm)=\CR\cdot\big(\mm\times(\BOm\cdot\mm){f}_0\big)$,
which is a consequence of Lemma \ref{Lem:antisymmetric} and Lemma \ref{Lem:ker-proj}.\ef

\begin{Lemma}\label{Lem:ker-out-positive}
For any symmetric matrix $\DD$, there holds
\begin{align*}
\int_\BS\mathcal{L}(\DD)\CA^{-1}_{f_0}\mathcal{L}(\DD)\ud\mm=&\big(\frac{3S_2+4S_4}{7}-\lambda{S_2}\big)|\DD\cdot\nn|^2\\
&-2\big(\frac{S_4}{35}-\frac{2S_2}{21}+\frac{1}{15}\big)\DD:\DD+(\lambda{S_2}-S_4)(\DD:\nn\nn)^2.
\end{align*}
\end{Lemma}

\no{\bf Proof}.\,
Applying Lemma \ref{lem:identity} with $f=\CA^{-1}_{f_0}\big(\CR\cdot(\mm\times\DD\cdot\mm{f_0})\big)$ and Lemma \ref{lem:tensor-product}, we get
\begin{align*}
&\big\langle\CR\cdot(\mm\times(\DD\cdot\mm){f_0}),\CA^{-1}_{f_0}\CR\cdot(\mm\times(\DD\cdot\mm){f_0})\big\rangle\\
&=-\frac{1}{2}\DD:\int_\BS
(\mm\mm-\frac{1}{3})\CR\cdot(\mm\times(\DD\cdot\mm){f_0})\ud\mm\\
&=-\frac{1}{2}\DD:\big(2\DD:\langle\mm\mm\mm\mm\rangle_{f_0}-\DD\cdot\langle\mm\mm\rangle_{f_0}-\langle\mm\mm\rangle_{f_0}\cdot\DD\big)\\
&=-\DD:\Big(S_4\nn\nn(\DD:\nn\nn)+\frac{2(S_2-S_4)}{7}
\big((\DD\cdot\nn)\nn+\nn(\DD\cdot\nn))\\
&\quad+2\big(\frac{S_4}{35}-\frac{2S_2}{21}
+\frac{1}{15}\big)\DD+\frac{S_2-S_4}{7}\II(\DD:\nn\nn)\Big)+S_2(\DD\cdot\nn)^2\\
&=\frac{3S_2+4S_4}{7}(\DD\cdot\nn)^2-2\big(\frac{S_4}{35}-\frac{2S_2}{21}
+\frac{1}{15}\big)\DD:\DD-S_4(\DD:\nn\nn)^2,
\end{align*}
which along with  Lemma \ref{Lem:ker-proj} gives
\begin{align*}
&\langle\mathcal{L}(\DD),\CA^{-1}_{f_0}\mathcal{L}(\DD)\rangle=\big\langle \Pout\big(\CR\cdot(\mm\times(\DD\cdot\mm){f_0})),
\CA^{-1}_{f_0}\Pout\big(\CR\cdot(\mm\times(\DD\cdot\mm){f_0})\big)\big\rangle\\
&=\big\langle\CR\cdot(\mm\times(\DD\cdot\mm){f_0}),\CA^{-1}_{f_0}\CR\cdot(\mm\times(\DD\cdot\mm){f_0})\big\rangle
-\lambda{S_2}|\nn\times(\DD:\nn)|^2\\
&=\big(\frac{3S_2+4S_4}{7}-\lambda{S}_2\big)(\DD\cdot\nn)^2
-2\big(\frac{S_4}{35}-\frac{2S_2}{21}
+\frac{1}{15}\big)\DD:\DD+(\lambda{S}_2-S_4)(\DD:\nn\nn)^2.
\end{align*}
The proof is finished.\ef\vspace{0.1cm}

\subsection{Proof of Theorem \ref{thm:energy} and application} Let us first prove Theorem \ref{thm:energy}.
By (\ref{Leslie cofficients1})-(\ref{Leslie cofficients2}) and Lemma \ref{Lem:ker-out-positive},
we find that
\begin{align}
&(\alpha_1+\frac{\gamma_2^2}{\gamma_1})|\DD:\nn\nn|^2
+\alpha_4\DD:\DD+(\gamma_3-\frac{\gamma_2^2}{\gamma_1})|\DD\cdot\nn|^2\nonumber\\
&=(-\frac{S_4}{2}+\frac{\gamma_2^2}{\gamma_1})|\DD:\nn\nn|^2
+(-\frac{S_4}{35}-\frac{5S_2}{21}+\frac{4}{15})\DD:\DD
+(\frac{5S_2+2S_4}{7}-\frac{\gamma_2^2}{\gamma_1})|\DD\cdot\nn|^2\nonumber\\
&=\big(\frac{3S_2+4S_4}{7}-\frac{\gamma_2^2}{\gamma_1}\big)(\DD\cdot\nn)^2
-2\big(\frac{S_4}{35}-\frac{2S_2}{21}+\frac{1}{15}\big)\DD:\DD
+\big(-S_4+\frac{\gamma_2^2}{\gamma_1}\big)(\DD:\nn\nn)^2\nonumber\\
&\quad+\frac{2(S_2-S_4)}{7}(\DD\cdot\nn)^2+\frac{S_4}{2}(\DD:\nn\nn)^2
+\big(\frac{S_4}{35}-\frac{3S_2}{7}+\frac{2}{5}\big)\DD:\DD\label{eq:corspc-energy}\\
&\ge\big(\frac{3S_2+4S_4}{7}-\frac{\gamma_2^2}{\gamma_1}\big)(\DD\cdot\nn)^2
-2\big(\frac{S_4}{35}-\frac{2S_2}{21}+\frac{1}{15}\big)\DD^2
+\big(-S_4+\frac{\gamma_2^2}{\gamma_1}\big)(\DD:\nn\nn)^2\ge0,\nonumber
\end{align}
since all the coefficients in the line (\ref{eq:corspc-energy}) are positive. Indeed, we have
\begin{align*}
S_2&=\langle\frac{1}{2}(3(\mm\cdot\nn)^2-1)\rangle_{f_0}=\frac{3A_2-A_0}{2A_0},\\
S_4&=\big\langle\frac{1}{8}(35(\mm\cdot\nn)^4-30(\mm\cdot\nn)^2+3)\big\rangle_{f_0}=\frac{35A_4-30A_2+3A_0}{8A_0}\\
&=\frac{1}{8A_0(2\eta)^2}(A_8-2A_6+A_4)>0.
\end{align*}
Hence,
\begin{align*}
&S_2-S_4=\frac{7}{8A_0}(6A_2-5A_4-A_0)=\frac{7}{16A_0\eta}(A_6-2A_4+A_2)>0.\\
&\frac{S_4}{35}-\frac{3S_2}{7}+\frac{2}{5}=\frac{1}{8A_0}(5A_0-6A_2+A_4)>0.
\end{align*}
This complete the proof of Theorem \ref{thm:energy}. \ef\vspace{0.1cm}

As a byproduct, we get the following dissipation law,
which has been used in the existence of the Hilbert expansion.

\begin{Lemma}\label{Lem:stress-dissipation}
For any matrix $\kappa$, there holds
\begin{align*}
\big\langle(\mm\mm-\frac{1}{3}\II)\mathcal{L}(\kappa)\big\rangle_1:\kappa \le0.
\end{align*}
\end{Lemma}

\no{\bf Proof}.\,Lemma \ref{lem:identity} implies that
\begin{align*}
\CA_{f_0}\big(\DD:(\mm\mm-\f13\II)\big)=-2\CR\cdot\big(\mm\times\DD\cdot\mm{f_0}\big)
=-2(\mathcal{K}(\DD)+\mathcal{L}(\DD)).
\end{align*}
Here  $\DD=\f12(\kappa+\kappa^T)$.
From Lemma \ref{Lem:ker-out-proj}, we know that $\mathcal{L}(\kappa)=\mathcal{L}(\DD)$. Hence,
\begin{align*}
&\big\langle(\mm\mm-\frac{1}{3}\II)\mathcal{L}(\kappa)\big\rangle_1:\kappa
=\int_\BS(\mm\mm-\frac{1}{3}\II):\DD\mathcal{L}(\DD)\ud\mm\\
&\quad=-2\big\langle\CA^{-1}_{f_0}\mathcal{L}(\DD),~\mathcal{K}(\DD)+\mathcal{L}(\DD)\big\rangle
=-2\big\langle\CA^{-1}_{f_0}\mathcal{L}(\DD),~\mathcal{L}(\DD)\big\rangle\le0.
\end{align*}
The proof is finished.\ef

\end{document}